\documentclass[a4paper,12pt]{amsart}
\usepackage{amssymb,amscd}
\usepackage[cp1251]{inputenc}

\usepackage{amsmath}
\usepackage{amsthm}
\usepackage{mathtext}
\usepackage{euscript}

\usepackage[matrix, arrow, curve]{xy}

\input xypic

\usepackage{graphicx}
\graphicspath{}
\DeclareGraphicsExtensions{.pdf,.png,.jpg,.jpeg}

\usepackage[english]{babel}

\theoremstyle{plain}
\newtheorem{theorem}{Theorem}[section]
\newtheorem{propos}[theorem]{Proposition}
\newtheorem{coro}[theorem]{Corollary}
\newtheorem{lemma}[theorem]{Lemma}
\newtheorem{prob}[theorem]{Problem}
\theoremstyle{definition}
\newtheorem{definition}[theorem]{Definition}
\newtheorem{exam}[theorem]{Example}
\newtheorem{constr}[theorem]{Construction}

\newtheorem{remark}[theorem]{Remark}
\numberwithin{equation}{section}

\newcommand{\zp}{\mathcal Z_P}
\newcommand{\rp}{\mathcal R_P}
\newcommand{\zk}{\mathcal Z_K}
\newcommand{\rk}{\mathcal R_K}
\DeclareMathOperator{\Dj}{\mbox{\textit{DJ}}}
\DeclareMathOperator{\sk}{\mbox{\textit{sk}}}

\DeclareMathOperator{\Tor}{Tor}
\DeclareMathOperator{\cone}{cone}

\DeclareMathOperator{\ko}{k}

\DeclareMathOperator{\conv}{conv}
\DeclareMathOperator{\pr}{pr}
\DeclareMathOperator{\cc}{cc}
\DeclareMathOperator{\vc}{vc}
\DeclareMathOperator{\fc}{fc}
\DeclareMathOperator{\MF}{MF}
\DeclareMathOperator{\Link}{Link}
\DeclareMathOperator{\im}{im}
\newcommand{\mb}[1]{{\textbf {\textit#1}}}
\newcommand{\field}[1]{\mathbb{#1}}
\def\C{\field{C}}
\def\k{\mathbf{k}}
\newcommand{\R}{\field{R}}
\newcommand{\Q}{\field{Q}}
\newcommand{\Z}{\field{Z}}

\begin{document}

\title[Massey products]{Massey products, toric topology\\and combinatorics of polytopes}
\author{Victor Buchstaber}
\address{Steklov Mathematical Institute of the Russian Academy of Sciences}
\email{buchstab@mi-ras.ru}
\author{Ivan Limonchenko}
\address{National Research University ``Higher School of Economics''}
\email{ilimonchenko@hse.ru}

\thanks{The research of the first author was supported by the Russian Foundation for Basic Research (grants nn. 17-01-00671 and 18-51-50005). The research of the second author was carried out within the University Basic Research Program of the Higher School of Economics and was funded by the Russian Academic Excellence Project 5-100.}




\maketitle

\begin{abstract}
In this paper we introduce a direct family of simple polytopes $P^{0}\subset P^{1}\subset\ldots$ such that for any $2\leq k\leq n$ there are non-trivial strictly defined Massey products of order $k$ in the cohomology rings of their moment-angle manifolds $\mathcal Z_{P^n}$. We prove that the direct sequence of manifolds $\ast\subset S^{3}\hookrightarrow\ldots\hookrightarrow\mathcal Z_{P^n}\hookrightarrow\mathcal Z_{P^{n+1}}\hookrightarrow\ldots$ has the following properties: every manifold $\mathcal Z_{P^n}$ is a retract of $\mathcal Z_{P^{n+1}}$, and one has inverse sequences in cohomology (over $n$ and $k$, where $k\to\infty$ as $n\to\infty$) of the Massey products constructed.
As an application we get that there are non-trivial differentials $d_k$, for arbitrarily large $k$ as $n\to\infty$ in Eilenberg--Moore spectral sequence connecting the rings $H^*(\Omega X)$ and $H^*(X)$ with coefficients in a field, where $X=\mathcal Z_{P^n}$.
\end{abstract}


\emph{Dedicated to the memory of Professor Evgeniy S. Golod (1935-2018)}

\section{Introduction}

In the paper~\cite{Go} of Golod (1962), Massey products were used to formulate a necessary and sufficient condition for the Serre inequality for Poincar\'e series of local rings to become an equality. This motivated the notions of a Golod ring and a Golod simplicial complex that currently attract much interest in combinatorial commutative algebra and toric topology. 
A number of questions arising in algebraic topology, homological algebra, and symplectic geometry lead us to the problem of nontriviality of Massey products.

Massey products were defined in the 1950s, and since that time they have been studied intensively in algebraic topology~\cite{M-U,Mass,M,Kr,N}. Some important problems of rational homotopy theory, homological algebra, group theory, knot theory and other areas of mathematics can be stated in terms of them.


For example, the existence of a non-trivial (matric) Massey product of order $k\geq 3$ in the (singular) cohomology ring $H^*(X;\mathbb Q)$ is an obstruction to the rational formality of a topological space $X$, see~\cite[Theorem 4]{BaTa}. A well-known classical result~\cite{D-G-M-S} states that any compact K\"ahler manifold is formal (over a field). 

Applying Gromov's embedding theorem for symplectic manifolds, Babenko and Taimanov~\cite{BaTa}, using the result that nontriviality of Massey products is preserved under symplectic blow-up procedure, obtained examples of simply connected symplectic manifolds with non-trivial Massey products in cohomology.
Non-formal manifolds having non-trivial Massey products in cohomology were also studied in the framework of symplectic geometry and topology~\cite{TrOp,RuTr}, homotopy theory of coordinate subspace arrangements and toric topology~\cite{BaskM,DS}.

The paper of Davis and Januszkiewicz~\cite{DJ} was devoted to introducing a topological analogue of algebraic nonsingular toric varieties. They suggested axiomatics, which allowed them to construct a model being a $2n$-dimensional smooth manifold with a compact torus action satisfying these axioms. 
Toric topology emerged from the papers of Buchstaber and Panov~\cite{bu-pa00-2}, based on the seminal paper~\cite{DJ}.
 
In toric topology, with any (abstract) simplicial complex $K$ on the set of vertices $[m]=\{1,2,\ldots,m\}$ one can associate its moment-angle complex $\zk$. The key object of the approach by Buchstaber and Panov appeared to be a moment-angle complex. Their crucial observation was that to {\emph{any}} simplicial complex $K$ we can associate a moment-angle complex $\zk$ in such a way that this correspondence is functorial. In this approach, a topological counterpart of a toric variety arised as a quotient space of a moment-angle manifold over a maximal freely acting compact torus subgroup. These manifolds were later called quasitoric and they satisfied the Davis and Januszkiewicz axiomatics. It was essential for the Buchstaber and Panov approach that when $K$ is a polytopal sphere its moment-angle complex $\zk$ is a manifold and the above functorial correspondence $K\to\zk$ respects the Poincar\'e duality. 

The idea of a moment-angle complex turned out to be closely related to the ideas of many constructions from homotopy theory and lead to the notion of a polyhedral product. The development of the general theory of polyhedral products was initiated in the paper of Bahri, Bendersky, Cohen, and Gitler~\cite{BBCG10}. 

The main topological and algebraic invariants of $\zk$ are its homotopy type, the cohomology ring $H^*(\zk;\ko)$, and the bigraded (or, algebraic) Betti numbers $\beta^{-i,2j}(\ko[K])$, $0\leq i,j\leq m$, of the corresponding Stanley--Reisner ring (the face ring) $\ko[K]$, where $\ko$ is a ring with identity. The ring $H^*(\zk;\ko)$ was described by Buchstaber and Panov~\cite{TT}, Theorem 4.5.4, as the (multigraded) Tor-algebra $\Tor_{\ko[v_{1},\ldots,v_{m}]}^{*,*}(\ko[K],\ko)$ of the Stanley--Reisner ring $\ko[K]=\ko[v_{1},\ldots,v_{m}]/I$ of the corresponding simplicial complex $K$. Bigraded Betti numbers are the ranks of the bigraded components of the Tor-algebra mentioned above. Their computation for various classes of face rings $\ko[K]$ is an important problem of combinatorial commutative algebra (monomial ring theory)~\cite{HH}. 

Categorical methods based on the use of homotopy colimits in the model categories of topological monoids and differential algebras~\cite{p-r-v04,P-R} enabled Panov and Ray~\cite{P-R} to obtain explicit homotopy decompositions, see~\cite[Proposition 8.1.1]{TT}, for moment-angle complexes and other toric spaces. Applying Sullivan -- de Rham functor of rational PL-forms $A_{PL}$, this implies toric and quasitoric manifolds are formal in the sense of rational homotopy theory. 

Bahri, Bendersky, Cohen, and Gitler~\cite{BBCG10} showed that after one suspension a wide class of polyhedral products over a simplicial complex $K$, including (real) moment-angle complexes, have a homotopy decomposition into wedges of suspensions over full subcomplexes in $K$. In combination with the constructions and technics of~\cite{GT,IK,G-P-T-W} this constitutes a powerful method to obtain explicit homotopy decompositions for polyhedral products.

Baskakov~\cite{BaskM} constructed the first example of a non-formal moment-angle complex. He introduced a class of triangulated spheres $K$ for which the corresponding moment-angle manifolds $\zk$ have non-trivial triple Massey products in cohomology. Generalizing Baskakov's results, Denham and Suciu~\cite{DS} asserted that a non-trivial triple Massey product of 3-dimensional classes in the cohomology ring $H^{*}(\zk)$ exists if and only if the 1-skeleton $\sk^{1}(K)$ contains an induced subgraph on 6 vertices isomorphic to one of the five graphs explicitly described in their paper. Later, Grbi\'c and Linton~\cite{GL1} corrected this result: they found one additional graph, which also realizes a non-trivial strictly defined triple Massey product of 3-dimensional classes. Therefore, there are six obstruction graphs that give one the possibility to realize such a Massey product in cohomology of a moment-angle complex. 

Using this criterion, Limonchenko~\cite{L2,L3} obtained necessary and sufficient conditions for a non-trivial triple Massey product to exist in $H^*(\zp)$, where $P=P_{\Gamma}$ is a graph-associahedron. It was proved by Zhuravleva~\cite{Zh} that there is a non-trivial triple Massey product in $H^*(\zp)$ for any Pogorelov polytope $P$. The latter class consists of all 3-dimensional flag polytopes without 4-belts of facets (in particular, it contains all fullerenes, see~\cite{BE2017}) and coincides with the class of polytopes having a (unique, up to isometry) right-angled realization in the Lobachevsky space $\mathbb L^3$~\cite{Pog,And,BE2017,BEMPP}. 

Grbi\'c and Linton~\cite{GL1} constructed the first example of a simplicial complex $K$ such that there exists a non-trivial triple Massey product of 3-dimensional classes with a non-zero indeterminacy in $H^{*}(\zk)$. Moreover, the main result of~\cite{GL1} is a criterion of existence of a non-trivial triple Massey product of 3-dimensional cohomology classes for a moment-angle complex. Note that in the case of graphs of~\cite{DS} there is a single valued Massey product defined on certain triples of elements of $H^3(\zk)$ and taking its value in $H^8(\zk)$. Grbi\'c and Linton found two new graphs for which there is a multivalued Massey product $\langle\alpha_{1},\alpha_{2},\alpha_{3}\rangle$ defined on certain triples of elements of $H^3(\zk)$ and taking its values in $H^8(\zk)$; it is defined modulo sums of products including either $\alpha_{1}$, or $\alpha_{3}$.

The first examples of polyhedral products with non-trivial $n$-fold Massey products in cohomology for any $n\geq 4$ were constructed by Limonchenko in~\cite{L1,L2}. Note that all the moment-angle manifolds appearing in~\cite{DS,L1,L2,Zh} in constructing non-trivial Massey products in cohomology, were 2-connected. Later, in~\cite{L3}, sufficient conditions for a strictly defined $n$-fold Massey product for any $n\geq 3$ to exist were formulated in terms of the vanishing of certain multigraded Betti numbers. These conditions hold, in particular, for any of the non-trivial triple Massey products of 3-dimensional classes described in~\cite{DS}, Theorem 6.1.1. The proof that the $n$-fold Massey product is strictly defined uses a multigraded structure (\cite{L3}, Lemma 3.3). The second author used this result to introduce in~\cite{L3} a family of polytopes such that their moment-angle manifolds may have arbitrarily large connectedness, and non-trivial strictly defined Massey products of any prescribed order may exist in their cohomology.

Recently, Grbi\'c and Linton~\cite{GL2} showed that for any given $n$-tuple of simplicial complexes there exists a simplicial complex K, obtained by performing a sequence of stellar subdivisions on their join, such that there is a nontrivial $n$-fold Massey product in $H^*(\zk)$. Their construction is a natural generalization of both the constructions due to Baskakov~\cite{BaskM} and Limonchenko~\cite{L1}. On the other hand, it was also proved in~\cite{GL2} that if a series of certain combinatorial operations performed on a simplicial complex $K$ (namely, edge contractions for edges that satisfy the so called ``link condition'') leads to a simplicial complex $\tilde K$ with a property that $H^*(\mathcal Z_{\tilde K})$ has a nontrivial $n$-fold Massey product of cohomology classes defined on induced subcomplexes in $\tilde K$ (see Theorem~\ref{zkcoh}), then the initial simplicial complex $K$ had the same property itself.

Note that in our constructions below all non-trivial higher Massey products are strictly defined. Based on our results and following the constructions of~\cite{GL2} one can obtain new families of polytopes with non-trivial Massey products or prove that a given simplicial complex $K$ has a moment-angle complex with a non-trivial triple or higher Massey product in cohomology.
We decided to keep the Problem~\ref{probMasseyFamilies}, which was formulated before~\cite{GL2} appeared. Now, we can say that using the results of~\cite[Proposition 4.2]{GL2} we immediately obtain the solution of Problem~\ref{probMasseyFamilies}.  
 

One of the important applications of Massey products in algebra is the theory of Golod rings and their generalizations, which is now actively developing within the framework of homological algebra, toric topology, and combinatorial commutative algebra.

Let $A$ be a Noetherian commutative local ring and let $\k=A/I$, where $I$ is the unique maximal ideal of $A$.
The Poincar\'e series $P(A;t)$ of $A$ is defined to be the Hilbert series of the graded $\k$-module $\Tor_{A}(\k,\k)$. In the early 1960s the following inequality was discussed in the literature:
$$
P(A;t)\leq\frac{(1+t)^m}{1-\sum\limits_{i\geq 0}c_{i}t^{i+1}},
$$
Here $m$ is equal to the number of elements in a minimal system of generators of the maximal ideal of the ring $A$, and the $c_{i}$ are the ranks of the homology groups of its Koszul complex. In~\cite{Go} it was refered to as Serre inequality.


The main result of the paper of Golod~\cite{Go} is a necessary and sufficient condition for the Serre inequality to become an equality. This condition says that the multiplication and all the Massey products, triple and higher, vanish in the Koszul homology of $A$.
Such rings were named \emph{Golod ring} in the monograph of Gulliksen and Levin~\cite{GL}.

The Koszul homology of the Stanley--Reisner ring $\ko[K]$ coincides with its Tor-algebra $\Tor_{\ko[v_{1},\ldots,v_{m}]}^{*,*}(\ko[K],\ko)$. Thus, thanks to a theorem of Buchstaber and Panov, the Golod property acquired a topological interpretation in toric topology: it is equivalent to the triviality of the multiplication and all the Massey products, triple and higher, in $H^{+}(\zk;\ko)$. In toric topology, simplicial complexes whose Stanley--Reisner rings are Golod over any field, are said to be \emph{Golod}, cf.~\cite[Definition 4.9.5]{TT}. 

In the framework of toric topology, Grbi\'c and Theriault~\cite{GT,G-T13,GT2}, Iriye and Kishimoto~\cite{IK,IK2014}, Grbi\'c, Theriault, Panov, and Wu~\cite{G-P-T-W}, and others considered a wide class of simplicial complexes well-known in combinatorics (shifted, shellable, totally fillable, dual sequentially Cohen-Macaulay). They proved the moment-angle complexes of the simplicial complexes in the above classes to be homotopy equivalent to wedges of spheres, and therefore, these simplicial complexes are Golod. 
As a counterexample to an earlier claim in~\cite{BJ}, Katth\"an~\cite{Kat} constructed the first example of a non-Golod Stanley--Reisner ring with trivial multiplication in its Tor-algebra.

In~\cite{BJ}, Berglund and J\"ollenbeck introduced the notion of a \emph{minimally non-Golod} Stanley--Reisner ring (resp., simplicial complex). This means that $K$ is not a Golod complex itself, but deleting any vertex from $K$ results in a Golod complex. 

Bosio--Meersemann~\cite{bo-me06}, Lopez de Medrano~\cite{LdM} and Gitler--Lopez de Medrano~\cite{G-LdM} considered a wide class of moment-angle manifolds being diffeomorphic to connected sums of manifolds each of which is a product of two spheres. 
For the corresponding classes of simple polytopes $P$, Limonchenko found necessary and sufficient conditions for minimal non-Golodness of their Stanley--Reisner rings $\ko[P]$ (over any field $\ko$). This includes the cases when $P$ is a generalized truncation polytope~\cite[Theorem 3.2]{L2014}, a dual even-dimensional neighbourly polytope~\cite[Proposition 3.6]{L2014}, or an $n$-dimensional simple polytope with $m\leq n+3$ facets~\cite[Theorem 3.5(c)]{L2015}.

As pointed out by Buchstaber and Panov, when $K$ is a polytopal sphere the moment-angle complex $\zk$ has a canonical smooth structure as a complete (nonsingular) intersection of Hermitian quadrics in $\C^n$. Bosio and Meersemann~\cite{bo-me06} introduced a notion of a 2-link and showed that $\zp$ or $\zp\times S^1$ has a structure of a complex manifold.
Therefore, thanks to the results on nontriviality of Massey products in cohomology of moment-angle manifolds, complex geometry obtained a new huge class of complex non-K\"ahler manifolds with zero integral 2-dimensional cohomology groups.  


Finally, we would like to draw attention to a relation between the theory of Golod rings and the theory of Lusternik--Schnirelmann category ($cat(\cdot)$) and theory of cohomology length ($cup(\cdot)$) for topological spaces, see~\cite{G-P-T-W,BG}.      

Significant progress has recently been achieved in the problem of finding algebraic and combinatorial conditions under which the moment-angle complex $\zk$ is homotopy equivalent to a wedge of spheres (condition (I)), or to a connected sum of manifolds, each of which is a product of two spheres (condition (II)). In particular, it was shown in~\cite{G-P-T-W} that for {\emph{flag}} complexes $K$ condition (I) is equivalent to the Golodness of the face ring $\ko[K]$ and also to the chordality of the graph $\sk^1(K)$, but condition (II) is equivalent to minimal non-Golodness of $\ko[K]$ and also to $K$ being a boundary of an $m$-gon with $m\geq 4$. 

The results described above convince us that, in general, they should be replaced by topological conditions: (I) by the condition $cat(\zk)=1$, and (II) by the condition $cat(\zk)=2$. Our results on the invariant $cat(X)$ and its relation to the Milnor spectral sequence are given in \S6.

This paper is an extended version of~\cite{BL}. For our preprints in this direction, see~\cite{BL1} and~\cite{BL2}. We construct a direct family $\mathcal P_{Mas}$: $P^{0}\subset P^{1}\subset\ldots$ of simple polytopes, where $P^n$ is a facet of $P^{n+1}$, and prove that for every $k$, $2\leq k\leq n$, non-trivial strictly defined $k$-fold Massey products exist in the cohomology rings of their moment-angle-manifolds $M_{n}=\mathcal Z_{P^n}$. It follows from general results of toric topology that $\mathcal P_{Mas}$ provides us with a direct sequence of smooth manifolds $\ast\hookrightarrow S^{3}\hookrightarrow\ldots\hookrightarrow M_{n}\hookrightarrow M_{n+1}\hookrightarrow\ldots$, where $M_{n}$ is a smooth submanifold in $M_{n+1}$. We show that every $M_{n}$ is a retract of $M_{n+1}$ and, in the inverse sequence of the cohomology rings of the manifolds $M_{n}$, one has inverse systems (over $n$ and $k$, where $k\to\infty$ as $n\to\infty$) of the constructed Massey products. 

Note that a retraction $\pi_{n+1}\colon M_{n+1}\to M_{n}$ can not be turned into a smooth map. This follows from results of differential topology, Theorem~\ref{zkcoh}, see Proposition~\ref{nofib}.   

It follows immediately that, for an infite-dimensional manifold $M=\varinjlim M_{n}$ there are non-trivial strictly defined Massey products of every order in  $\varprojlim\,H^{*}(M_n)$. We plan to study the properties of the filtration $\mathcal Z_{P^0}\subset\mathcal Z_{P^1}\subset\ldots\subset Z_{P^n}\subset\ldots$ and topology of the direct limit, $M=\mathcal Z_{\mathcal P}$, in subsequent publications. 

In the first part of this paper we discuss known constructions and methods of toric topology and polytope theory. We also introduce the new ones, which we need in the second part. 
We develop the theory of the differential ring of polytopes~\cite{B}, compute the action of the boundary operator on the polytopes in our families (\S4) and introduce differential equations for the 2-parameter generating series of these families. 

We prove that polytopes in $\mathcal P_{Mas}$ are flag nestohedra, where $P^0$ is a point, $P^1$ is a segment, $P^2$ is a square, and $P^n$, $n>2$, admits a non-trivial strictly defined Massey product $\langle\alpha_{1},\ldots,\alpha_{k}\rangle$ in $H^*(\mathcal Z_{P^n})$ with $\dim\alpha_{i}=3$, where $2\leq k\leq n$. We introduce new notions in the theory of polytope families: complexity of a family, algebraic and geometric direct families of polytopes, as well as a special direct family of polytopes with non-trivial Massey products. In these new terms our main results says that the family $\mathcal P_{Mas}$ is a special geometric direct family of polytopes of complexity 4 with non-trivial Massey products.

In~\cite{BV}, a family of 2-truncated cubes was introduced and proved to contain a number of famous polytope families arising in different areas of Mathematics. Our  sequence $\mathcal P_{Mas}$ belongs to this family (\S5).

We introduce a family of smooth closed manifolds $\{M_n\}$, with a single manifold for every $n\geq 2$, such that $cup(M_n)=cat(M_n)=n$. For this sequence of manifolds we show: 

(1) the differentials $d_{r},1\leq r\leq n-1$ in the Eilenberg--Moore spectral sequence for $M_n$ are non-trivial, 

(2) the differentials $d_{r}, r\geq n+1$ in the Milnor spectral sequence for $M_n$ vanish. 

We give necessary details on the Milnor spectral sequence in \S6. 

In the last section we discuss results closely related to problems in which Massey products play an important role, as well as problems in the theory of Massey products arising in topology, algebra and elsewhere.

We are grateful to Anthony Bahri, Nickolay Erokhovets, Jelena Grbi\'c and Taras Panov for many fruitful discussions during the course of this research.

\section{Moment-angle manifolds: definitions and constructions}

We start this section with the next basic notion.

\begin{definition}
An \emph{(abstract) simplicial complex} $K$ on the vertex set $[m]=\{1,\ldots,m\}$ is a subset in $2^{[m]}$ such that if $\sigma\in K$ and $\tau\subseteq\sigma$, then $\tau\in K$.

Elements of $K$ are called its \emph{simplices} and the maximal dimension of a simplex $\sigma\in K$ is the \emph{dimension} of $K$ and denoted by $\dim K=n-1$, where $\max|\sigma|=n$. Finally, for any subset of vertices $I\subset [m]$ the subset of $K$ that is equal to $K\cap 2^{I}$ is itself a simplicial complex (on the vertex set $I$). This complex is called a \emph{full subcomplex} in $K$ (on $I$) and is denoted by $K_{I}$.
\end{definition}

We assume in what follows that there are no \emph{ghost vertices} in $K$, that is, $\{i\}\in K$ holds for any $1\leq i\leq m$. A simplicial complex $K$ is a partially ordered set with the natural ordering by inclusion. Therefore $K$ (and every full subcomplex $K_{I}, I\subset [m]$) is determined by its maximal (by inclusion) simplices and $\dim K$ is equal to the dimension of one of them. If all the maximal simplices in $K$ have the same dimension, then $K$ is called a \emph{pure} simplicial complex. 

We use the following definition of a simple (convex) polytope. 

\begin{definition}\label{Simplepolytopes}
A \emph{convex $n$-dimensional polytope} $P$ in Euclidean space $\R^n$ with inner product $\langle\;,\:\rangle$ is a bounded intersection of $m$ closed halfspaces, that is,
 $$
  P=\bigl\{\mb x\in\R^n\colon\langle\mb a_i,\mb
  x\rangle+b_i\ge0\quad\text{for }
  i=1,\ldots,m\bigr\},\eqno (1)
$$
where $\mb a_i\in\R^n$, $b_i\in\R$. 
Here we assume that \emph{facets} of $P$, 
$$
  F_i=\bigl\{\mb x\in P\colon\langle\mb a_i,\mb
  x\rangle+b_i=0\bigr\},\quad\text{for } i=1,\ldots,m.
$$
are in general position and that there are no redundant inequalities above, that is, none of the inequalities can be removed without changing the set of points~$P$ (equivalently, there are no ghost vertices in $K_P=\partial P^*$). In particular, the latter means that the polytope $P$ is {\emph{simple}}, that is, every vertex is the intersection of exactly $n$ facets.
\end{definition}

Let $P$ be an $n$-dimensional convex simple polytope with $m$ facets $F_{1},\ldots,F_{m}$. In this paper, we are only interested in the structure of the partially ordered set of faces (that is, the combinatorial equivalence class), rather than in a particular embedding in the ambient Euclidean space $\R^n$. In what follows we write $P=Q$ when $P$ and $Q$ are combinatorially equivalent. 

\begin{exam}
Consider the nerve of the (closed) covering of the boundary $\partial P$ of a polytope $P$ by all its facets $F_{i}, 1\leq i\leq m$. The resulting simplicial complex on the vertex set $[m]$ is called the \emph{nerve complex} of $P$ and denoted by $K_P$. Note that the geometric realization of $K_P$ coincides with the boundary of the simplicial polytope $P^*$ combinatorially dual to $P$. The nerve complex $K_P$ of a simple $n$-dimensional polytope $P$ is obviously a pure simplicial complex of dimension $n-1$.  
\end{exam}

It is easy to see that a simplicial complex $K$ can either be determined by all its maximal simplices or by all its minimal non-faces, that is, the subsets $I\in 2^{[m]}\backslash K$ minimal by inclusion. It is obvious that $I$ is a minimal non-face of $K$ if and only if $I\notin K$ and each of its proper subsets is a simplex in $K$. We denote the set of minimal non-faces of $K$ by $\MF(K)$. Thus, $I\in\MF(K)$ if and only if $K_{I}=\partial\Delta_I$. Using $\MF(K)$, we can easily introduce on the set of simplicial complexes and simple polytopes the next combinatorial operation originally used in the framework of toric topology by Bahri, Bendersky, Cohen, and Gitler~\cite{BBCG15}. 

\begin{constr}\label{simpmultwedge}
Let $J=(j_{1},\ldots,j_{m})$ be an ordered set of $m$ positive integers. Consider the following set of vertices:
$$
m(J)=\{11,\ldots,1j_{1},\ldots,m1,\ldots,mj_{m}\}.
$$
In order to define a \emph{simplicial multiwedge (or, a $J$-construction)} of $K$, which is a simplicial complex $K(J)$ on $m(J)$, we say that its set of minimal non-faces $\MF(K(J))\subset 2^{[m(J)]}$ consists of the subsets of $m(J)$ of the following type 
$$
I(J)=\{i_{1}1,\ldots,i_{1}j_{i_1},\ldots,i_{k}1,\ldots,i_{k}j_{i_k}\},
$$
where $I=\{i_{1},\ldots,i_{k}\}\in\MF(K)$. Note that if $J=(1,\ldots,1)$, then $K(J)=K$. 
\end{constr}

\begin{exam}
Suppose $K=K_P$ is a polytopal sphere on $m$ vertices, that is, a nerve complex of an $n$-dimensional polytope with $m$ facets. According to~\cite[Theorem 2.4]{BBCG15}, its simplicial multiwedge $K(J)$ is then a polytopal sphere for any vector $J$, and so it is also a nerve complex of a certain simple polytope $Q=P(J)$, and therefore $K_{P(J)}=K_{P}(J)$. If we denote $d(J)=j_{1}+\ldots+j_{m}$, then it is easy to see that $m(P(J))=d(J)$ and $n(P(J))=d(J)+n-m$. Thus, the value $m(P(J))-n(P(J))=m(P)-n(P)$ remains the same after applying a $J$-construction. 
\end{exam}

For the geometric description of the multiwedge operation on polytopes and its applications we refer to~\cite{G-LdM}.  

In what follows we shall denote by $\ko$ a field of zero characteristic of the ring of integers. Let $\ko[v_{1},\ldots,v_{m}]$ be a graded polynomial algebra with $m$ generators and $\deg(v_{i})=2$. 

\begin{definition}
By a \emph{Stanley--Reisner ring}, or \emph{face ring} of a simplicial complex  $K$ (over $\ko$) we mean the quotient ring
$$
   \ko[K]=\ko[v_{1},\ldots,v_{m}]/I(K),
$$
where $I(K)$ is an ideal generated by square-free monomials of the type $v_{i_{1}}\cdots{v_{i_{k}}}$ such that $\{i_{1},\ldots,i_{k}\}\in\MF(K)$. The monomial ideal $I(K)$ is called the \emph{Stanley--Reisner ideal} of~$K$. Observe that the natural projection gives $\ko[K]$ a structure of a $\ko$-algebra and that of a module over polynomial algebra $\ko[v_{1},\ldots,{v_{m}}]$. 
\end{definition}

By a result of Bruns and Gubeladze~\cite{br-gu96}, two simplicial complexes $K_1$ and $K_2$ are simplicially isomorphic if and only if their Stanley--Reisner rings are isomorphic. Thus, $\ko[K]$ is a complete combinatorial invariant of a simplicial complex $K$.

\begin{definition}
(1) A simplicial complex $K$ is called \emph{flag} if it coincides with the
\emph{clique complex} $\Delta(\Gamma)$ of its {\emph{graph}} $\Gamma=\sk^1(K)$, that is, for any subset of vertices $I\subset [m]$ which are pairwisely joint by edges in $\Gamma$, the corresponding full subcomplex $K_{I}=\Delta_{I}$. A simple polytope $P$ is called \emph{flag}, if its nerve complex $K_P$ is flag. (2) A simplicial complex $K$ is called $q$-\emph{connected} ($q\geq 1$) if $K_{I}=\Delta_{I}$ for any subset of vertices $I\subset [m]$ in $K$ with $|I|=q$ elements. Note that 1-connectedness is equivalent to having no ghost vertices.  
\end{definition}

\begin{remark}
(1) $K$ is a flag simplicial complex if and only if for any $I\in\MF(K)$ one has $|I|=2$, that is, $I(K)$ is generated by monomials of degree 4. (2) $K$ is a $q$-connected simplicial complex if and only if for any $I\in\MF(K)$ one has $|I|\geq q+1$.
\end{remark}

The next class of toric spaces is a main object of study in toric topology. Let $({\bf{X}},{\bf{A}})=\{(X_i,A_i)\}_{i=1}^{m}$ be an ordered set of topological pairs. The case $X_{i}=X, A_{i}=A$ was firstly introduced in~\cite{bu-pa00-2} as a \emph{K-power} and was then intensively studied and generalized in a series of more recent papers, among which are~\cite{BBCG10,G-T13,IK}. 

\begin{definition}
A {\emph{polyhedral product}} over a simplicial complex $K$ on the vertex set $[m]$ is a topological space
$$
({\bf{X}},{\bf{A}})^K=\bigcup\limits_{I\in K}({\bf{X}},{\bf{A}})^I\subseteq\prod\limits_{i=1}^{m}\,X_{i},
$$
where $({\bf{X}},{\bf{A}})^I=\prod\limits_{i=1}^{m} Y_{i}$ and $Y_{i}=X_{i}$ if $i\in I$, and $Y_{i}=A_{i}$ if $i\notin I$.
\end{definition}
The term ``polyhedral product'' was suggested by Browder (cf.~\cite{BBCG10}).

\begin{exam}
Suppose $X_{i}=X$ and $A_{i}=A$ for all $1\leq i\leq m$. Then the next spaces are particular cases of the above construction of a polyhedral product $({\bf{X}},{\bf{A}})^{K}=(X,A)^K$.
\begin{itemize}
\item[(1)] {\emph{Moment-angle-complex}} 
$\zk=(\mathbb{D}^2,\mathbb{S}^1)^K$;
\item[(2)] {\emph{Real moment-angle-complex}} 
$\mathcal R_K=([-1;1],\{-1,1\})^K$;
\item[(3)] \emph{Davis--Januszkiewicz space} $\Dj(K)\simeq(\mathbb{C}P^{\infty},*)^K$;
\item[(5)] Cubical complex $\cc(K)=(I^1,1)^K$ in the $m$-dimensional cube $I^{m}=[0,1]^m$, which is PL-homeomorphic to a cone over a barycentric subdivision of $K$.
\end{itemize} 
\end{exam}

The next construction appeared firstly in the work of Davis and Januszkiewicz~\cite{DJ}.

\begin{constr}(Moment-angle manifolds I)\label{mamfdDJ}
\begin{itemize}
\item[(1)] Suppose $P^n$ is a simple polytope with a set of facets $\mathcal F(P^n)=\{F_{1},\ldots,F_{m}\}$. Denote by $T^{F_{i}}$ the corresponding 1-dimensional coordinate subgroup in the $m$-dimensional torus $T^{\mathcal F}\cong T^{m}$ for each $1\leq i\leq m$ and denote by $T^{G}=\prod\,T^{F_i}\subset T^{\mathcal F}$, where $G=\cap\,F_{i}$ in the polytope $P^n$. Then a \emph{moment-angle manifold} of~$P$ is defined to be a quotient space
$$
\zp=T^{\mathcal F}\times P^{n}/\sim,
$$
where $(t_{1},p)\sim (t_{2},q)$ if and only is $p=q\in P$, $t_{1}t_{2}^{-1}\in T^{G(p)}$, and $G(p)$ is a minimal, with respect to inclusion relation, face of $P$ containing the point $p=q$;
\item[(2)] Suppose $P^n$ is a simple polytope with a set of facets $\mathcal F(P^n)=\{F_{1},\ldots,F_{m}\}$. Denote by $\mathbb{Z}_{2}^{F_{i}}$ the corresponding 1-dimensional coordinate subgroup in the $m$-dimensional real torus $\mathbb{Z}_{2}^{\mathcal F}\cong\mathbb{Z}_{2}^{m}$ for each $1\leq i\leq m$ and denote by $\mathbb{Z}_{2}^{G}=\prod\,\mathbb{Z}_{2}^{F_i}\subset\mathbb{Z}_{2}^{\mathcal F}$, where $G=\cap\,F_{i}$ in the polytope $P^n$. Then a \emph{real moment-angle manifold} of~$P$ is defined to be a quotient space
$$
\rp=\mathbb{Z}_{2}^{\mathcal F}\times P^{n}/\sim,
$$
where $(t_{1},p)\sim (t_{2},q)$ if and only if $p=q\in P$, $t_{1}-t_{2}\in\mathbb{Z}_{2}^{G(p)}$, and $G(p)$ is a minimal, with respect to inclusion relation, face of $P$, containing the point $p=q$.
\end{itemize}
\end{constr}

\begin{remark}
By means of Construction~\ref{mamfdDJ} one can prove that if $P_{1}$ and $P_{2}$ are \emph{combinatorially equivalent}, i.e, their face posets are isomorphic (or, equivalently, $K_{P_1}$ and $K_{P_2}$ are simplicially isomorphic), then their moment-angle manifolds $\mathcal Z_{P_{1}}$ and $\mathcal Z_{P_{2}}$ are homeomorphic, as well as their real moment-angle manifolds $\mathcal R_{P_1}$ and $\mathcal R_{P_2}$. The opposite statement fails to be true already for the truncation polytopes $P=\vc^{k}(\Delta^n)$ for $k\geq 3$, cf.~\cite{bo-me06}.
\end{remark}

Next we shall consider another, equivalent, construction which also leads to moment-angle manifolds.

Let $A_P$ be an $m\times n$-matrix with rows $\mb a_i\in\mathbb{R}^n$, and let $\mb b_P\in\mathbb{R}^m$ be a vector-column of numbers $b_i\in\R$. Then we can rewrite the inequalities in $(1)$ 
in the following way
\[
  P=\bigl\{\mb x\in\R^n\colon A_P\mb x+\mb b_P\ge\mathbf 0\},
\]
and consider an affine map: 
\[
  i_P\colon \R^n\to\R^m,\quad i_P(\mb x)=A_P\mb x+\mb b_P.
\]
This map, obviously, embeds the polytope $P$ into the nonnegative ortant
\[
  \R^m_\ge=\{\mb y\in\R^m\colon y_i\ge0\quad\text{for }
  i=1,\ldots,m\}.
\]

Moment-angle manifolds and complexes were studied in a series of works by Buchstaber and Panov by means of algebraic topology, combinatorial commutative algebra, and polytope theory. These investigations has led to the foundation of a new area of geometry and topology, Toric Topology, see~\cite{BP04,TT}. In particular, they introduced the following constrcution and proved it to be equivalent to the Constrcution~\ref{mamfdDJ}.

\begin{constr}(Moment-angle manifolds II)\label{mamfdBP}
\begin{itemize}
\item[(1)] Define a \emph{moment-angle manifold} $\mathcal Z_P$ of a polytope $P$ as a pullback from the following commutative diagram
$$\begin{CD}
  \mathcal Z_P @>i_Z>>\C^m\\
  @VVV\hspace{-0.2em} @VV\mu V @.\\
  P @>i_P>> \R^m_\ge
\end{CD}\eqno 
$$
where $\mu(z_1,\ldots,z_m)=(|z_1|^2,\ldots,|z_m|^2)$. The projection map $\zp\rightarrow P$ in the above diagram is a quotient map of the canonical compact torus $\mathbb{T}^m$ action on $\zp$, induced by the standard coordinatewise action
\[
  \mathbb T^m=\{\mb z\in\C^m\colon|z_i|=1\quad\text{for }i=1,\ldots,m\}
\]
on~$\C^m$. Therefore, $\mathbb T^m$ acts on $\zp$ with an orbit space $P$, and $i_Z$ is a $\mathbb T^m$-equivariant embedding;
\item[(2)] Define a \emph{real moment-angle manifold} $\mathcal R_P$ of a polytope $P$ as a pullback from the following commutative diagram
$$\begin{CD}
  \mathcal R_P @>i_R>>\R^m\\
  @VVV\hspace{-0.2em} @VV\mu V @.\\
  P @>i_P>> \R^m_\ge
\end{CD}\eqno 
$$
where $\mu(x_1,\ldots,x_m)=(|x_1|^2,\ldots,|x_m|^2)$. The projection map $\rp\rightarrow P$ in the above diagram is a quotient map of the canonical action of $\mathbb{Z}_{2}^m$ on $\rp$, induced by the standard coordinatewise action
\[
  \mathbb Z_{2}^{m}=\{\mb x\in\R^m\colon|x_i|=1\quad\text{for }i=1,\ldots,m\}
\]
on~$\R^m$. Therefore, $\mathbb Z_{2}^{m}$ acts on $\rp$ with an orbit space $P$, and $i_R$ is a $\mathbb Z_{2}^{m}$-equivariant embedding.
\end{itemize}
\end{constr}

It follows immediately from the Construction~\ref{mamfdBP}, see~\cite[\S3]{BR}, that $\zp$ is a complete intersection of Hermitian quadrics in $\C^m$, and $\rp$ is a nondegenrate intersection of real quadrics in $\R^m$. Thus, $\zp$ and $\rp$ both acquire canonical equivariant smooth structures. For any simple $n$-dimensional polytope $P$ with $m$ facets, its moment-angle manifold $\zp$ will be a 2-connected closed $(m+n)$-dimensional smooth manifold, and its real moment-angle manifold $\rp$ will be (non simply connected) closed orientable $n$-dimensional manifold.

\begin{remark}
The following classical problem is well known: how many $K$-invariant smooth structures are there on $\zp$, where $K\cong T^{m-n}$ is a maximal subgroup in $\mathbb{T}^m$, freely acting on $\zp$. This problem remains open. In the case $P=\Delta^3$ and $K\cong S^1$ (i.e., $\zp=S^7$), it was solved in the work of Bogomolov~\cite{Bg}.
\end{remark}

The next properties of the previously defined polyhedral products hold.
\begin{itemize}
\item Suppose $K=K_P$ is a nerve complex of an $n$-dimensional simple polytope $P$ with $m$ facets. In~\cite{bu-pa00-2} canonical equivariant homemorphisms $h_{P}:\zp\cong\mathcal Z_{K_P}$ and $h^{\R}_{P}:\rp\cong\mathcal R_{K_P}$ were constructed;
\item If $K$ is $q$-connected, then $\zk$ is a $2q$-connected CW-complex;
\item One has a commutative diagram
$$\begin{CD}
  \zk @>>>(\mathbb{D}^2)^m\\
  @VVrV\hspace{-0.2em} @VV\rho V @.\\
  \cc(K) @>i_c>> I^m
\end{CD}\eqno 
$$
where $i_{c}:\,\cc(K)\hookrightarrow I^{m}=(I^1,I^1)^{[m]}$ is an embedding of a cubical subcomplex, induced by an embedding of pairs: $(I^1,1)\subset (I^1,I^1)$, the maps $r$ and $\rho$ are projection maps on the orbit spaces of a torus $\mathbb{T}^m$-action, induced by coordinatewise action of $\mathbb{T}^m$ on the complex unitary polydisk $(\mathbb{D}^2)^m$ in $\C^m$;
\item One has a commutative diagram
$$\begin{CD}
  \rk @>>>[-1;1]^m\\
  @VVrV\hspace{-0.2em} @VV\rho V @.\\
  \cc(K) @>i_c>> I^m
\end{CD}\eqno 
$$
where $i_{c}:\,\cc(K)\hookrightarrow I^{m}=(I^1,I^1)^{[m]}$ is an embedding of a cubical subcomplex, induced by an embedding of pairs: $(I^1,1)\subset (I^1,I^1)$, the maps $r$ and $\rho$ are projection maps on the orbit spaces of a group $\mathbb{Z}_{2}^m$-action, induced by a coordinatewise action of $\mathbb{Z}_{2}^m$ on the ``big'' cube $[-1;1]^m$;

\item There is a homotopy fibration
$$\begin{CD}
  X @>>> ET^m\\
  @VVV\hspace{-0.2em} @V\pi VV @.\\
  \Dj(K) @>p>> BT^{m}
\end{CD}\eqno 
$$
where $\pi$ denotes the universal $\mathbb{T}^m$-bundle, and the map $p:\,\Dj(K)\rightarrow BT^{m}$ is induced by the embedding of pairs $(\mathbb{C}P^{\infty},*)\subset (\mathbb{C}P^{\infty},\mathbb{C}P^{\infty})$. 
Furthermore, its homotopy fiber $X\simeq\zk$ and the Davis--Januszkiewicz space $\Dj(K)=ET^{m}\times_{\mathbb{T}^m}\zk$ is homotopy equivalent to the polyhedral product $(\mathbb{C}P^{\infty},*)^K$. Therefore, for the equivariant cohomology of moment-angle-complex $\zk$ one has:
$$
H^{*}_{\mathbb{T}^m}(\zk)=H^{*}(\Dj(K))\cong\mathbb{Z}[K].
$$
The details can be found in~\cite[Ch.~4,8]{TT}.
\end{itemize}

Using the homotopy fibration defined above, Buchstaber and Panov~\cite{bu-pa00-2} showed that its Eilenberg--Moore spectral sequence, which converges to $H^{*}(\zk,\ko)$ over a field $\ko$, collapses in the $E_2$-term. Moreover, in their joint paper with Baskakov~\cite{BBP} they obtained the next result on the structure of the cohomology algebra of $\zk$ with coefficients in an arbitrary ring with unit $\ko$. 

\begin{theorem}[{\cite[Theorem 4.5.4]{TT}}]\label{zkcoh}
The cohomology algebra of a moment-angle complex $\zk$ is given by isomorphisms
\[
\begin{aligned}
  H^{*,*}(\mathcal Z_K;\ko)&\cong\Tor_{\ko[v_1,\ldots,v_m]}^{*,*}(\ko[K],\ko)\\
  &\cong H^*\bigl[R(K),d\bigr]\\
  &\cong \bigoplus\limits_{I\subset [m]}\widetilde{H}^{*}(K_{I};\ko),
\end{aligned}
\]
in which bigrading and differential in cohomology of a differential (bi)graded algebra are defined in the following way:
\[
  \mathop{\mathrm{bideg}} u_i=(-1,2),\;\mathop{\mathrm{bideg}} v_i=(0,2);\quad
  du_i=v_i,\;dv_i=0
\]
and $R^*(K)=\Lambda[u_{1},\ldots,u_{m}]\otimes\ko[K]/(v_{i}^{2}=u_{i}v_{i}=0,1\leq i\leq m)$.
Here by $\widetilde{H}^*(K_{I})$ (in what follows we drop $\ko$ in the notations) we denoted reduced simplicial cohomology of a simplicial complex $K_{I}$. The last isomorphism above is a sum of isomorphisms 
$$
H^{p}(\zk)\cong\sum\limits_{I\subset [m]}\widetilde{H}^{p-|I|-1}(K_{I}).
$$
In order to determine a product of two cohomology classes $\alpha=[a]\in\tilde{H}^{p}(K_{I_1})$ and $\beta=[b]\in\tilde{H}^{q}(K_{I_2})$, let us introduce an embedding of subsets $i:\,K_{I_{1}\sqcup I_{2}}\rightarrow K_{I_1}*K_{I_2}$ and a canonical $\ko$-module isomorphism of cochains:
$$
j:\,\tilde{C}^{p}(K_{I_1})\otimes\tilde{C}^{q}(K_{I_2})\rightarrow\tilde{C}^{p+q+1}(K_{I_{1}}*K_{I_2}).
$$
Then a product of the classes $\alpha$ and $\beta$ is given by:
$$
\alpha\cdot\beta=\begin{cases}
0,&\text{if $I_{1}\cap I_{2}\neq\varnothing$;}\\
i^{*}[j(a\otimes b)]\in\tilde{H}^{p+q+1}(K_{I_{1}\sqcup I_{2}}),&\text{if $I_{1}\cap I_{2}=\varnothing$.}
\end{cases}
$$
\end{theorem}

Additive structure of the Tor-algebra described in Theorem~\ref{zkcoh}, can be computed using the free Koszul resolution for $\ko$ viewed as a $\ko[m]=\ko[v_1,\ldots,v_m]$-module. Multiplicative structure can be obtained using the free Taylor resolution for $\ko[K]$, viewed as a $\ko[m]$-module. In general, the latter resolution will no longer be minimal, however it can often be useful in combinatorial proofs. By Hochster theorem~\cite{Hoch} we get that the $\ko$-module structure on $\Tor_{\ko[v_1,\ldots,v_m]}^{*,*}(\ko[K],\ko)\cong H^*(\zk;\ko)$ is determined by all the reduced simplicial cohomology groups of all full subcomplexes in $K$ (including $\varnothing$ and the $K$ itself). 
 
We shall make use of the following multigraded version of Theorem~\ref{zkcoh}. Indeed, differential graded algebra $R^*(K)$ has a natural multigrading and the next statement holds.

\begin{theorem}[{\cite[Construction 3.2.8, Theorem 3.2.9]{TT}}]\label{mgrad}
For any simplicial complex $K$ on the vertex set $[m]$ one has
$$
\Tor^{-i,2(j_{1},\ldots,j_{m})}_{\ko[v_{1},\ldots,v_{m}]}(\ko[K],\ko)\cong\widetilde{H}^{|J|-i-1}(K_{J};\ko),
$$
where for a subset $J\subset [m]$ we set $j_{k}=1$, if $k\in J$, and $j_{k}=0$, if $k\notin J$.
$\Tor^{-i,2{\bf{a}}}_{\ko[v_{1},\ldots,v_{m}]}(\ko[K],\ko)=0$, when ${\bf{a}}$ is not a $(0,1)$-vector. Furthermore, 
$$
\Tor^{-i,2{\bf{a}}}_{\ko[v_{1},\ldots,v_{m}]}(\ko[K],\ko)\cong H^{-i,2{\bf{a}}}[R(K),d].
$$
\end{theorem}


\begin{constr}\label{fullsubcompretract}
Consider a subset of vertices $S\subset [m]$ of a simplicial complex $K$. There exists a natural simplicial embedding $j_{S}:\,K_{S}\to K$. We shall construct here a retraction to the induced embedding of moment-angle-complexes $\hat{j}_{S}:\mathcal Z_{K_S}\to\mathcal\zk$. To do this, consider a projection map
$$
p_{\sigma}:\,\prod\limits_{i\in\sigma}\,D^{2}\times\prod\limits_{i\notin\sigma}\,S^{1}\rightarrow\prod\limits_{i\in\sigma\cap S}D^{2}\times\prod\limits_{S\backslash\sigma}S^{1}
$$
for each $\sigma\in K$. Note that the image $p_{\sigma}$ belongs to $\mathcal Z_{K_S}$ for any $\sigma\in K$, since $K_{S}=\{\sigma\cap S|\,\sigma\in K\}$. It is easy to see that $r_{S}:=\cup_{\sigma\in K}\,p_{\sigma}:\,\zk\to\mathcal Z_{K_S}$ is the desired retraction. A similar formula in th case $\rk=(D^{1},S^{0})^K$ yields a retraction for the induced embedding of real moment-angle-complexes $\hat{j}_{S}:\mathcal R_{K_S}\to\mathcal\rk$; the general result can be found in~\cite[Proposition 2.2]{PV}
\end{constr}

\begin{coro}\label{splitepi}
Let $j_{S}:\,K_{S}\hookrightarrow K$ be an embedding of a full subcomplex on the vertex set $S\subseteq [m]$ in $K$. Then the embedding of the corresponding moment-angle complexes $\hat{j}_{S}:\,\mathcal Z_{K_S}\to\zk$ has a retraction, and the induced cohomology ring homomorphism $j^*_{S}:\,H^*(\zk)\rightarrow H^*(\mathcal Z_{K_S})$ is a split epimorphism of rings. Similarly, the induced embedding of real moment-angle-complexes $\hat{j}_{S}:\,\mathcal R_{K_S}\to\rk$ has a retraction, and the induced cohomology ring homomorphism $j^*_{S}:\,H^*(\rk)\rightarrow H^*(\mathcal R_{K_S})$ is a split epimorphism of rings.
\end{coro}
\begin{proof}
The statement about embedding of moment-angle-complexes follows directly from their definition and Construction~\ref{fullsubcompretract} (cf.~\cite[Exercise 4.2.13]{TT}). The rest follows from that homomorphism in cohomology, induced by a retraction, is a split epimorphism of rings: $(ri)^{*}=i^{*}r^{*}=1_{A}^{*}$ for a retract $A=\mathcal Z_{K_S}\subset\zk$.
\end{proof}

The main result of this section is the following: if $P$ is a flag polytope and $F\subset P$ is its face, then the maps $\hat{i}_{F,P}:\mathcal Z_{F}\to\mathcal Z_P$ and $\hat{j}_{F,P}:\mathcal Z_{K_F}\to\mathcal Z_{K_P}$, induced by an affine embedding of a face $i_{F,P}:\,F\to P$ and an embedding of the corresponding full subcomplex $j_{F,P}:K_F\to K_P$, and linked together by canonical equivariant homeomorphisms $h_{P}:\mathcal Z_P\to\mathcal Z_{K_P}$ and $h_{F}:\mathcal Z_{F}\to\mathcal Z_{K_F}$ in a commutative diagram. Later on, in section 5, we shall show that a similar result holds for all nestohedra (not necessarily flag ones), see Definition~\ref{nest}. That is, we shall prove that if $P=P_B$ is a nestohedron, then a face embedding $i_{P_{B|_S}}$ and an embedding of the corresponding full subcomplex $j_{P_{B|_S}}$ induce the maps of moment-angle manifolds and complexes, which form a commutative diagram with the equivariant homeomorphisms described above, see Proposition~\ref{phiPsi}.

Thus, we now need to discuss various ways of constructing an equivariant embedding of a moment-angle manifold $\mathcal Z_{F^r}$ to a moment-angle manifold $\mathcal Z_{P^n}$, induced by a face emebedding $i^{n}_{r}:\,F^{r}\to P^{n}$. Note that, although $\mathcal Z_{F^r}$ will always ne a submanifold in $\mathcal Z_{P^n}$, a retraction $\mathcal Z_{P^n}\to\mathcal Z_{F^r}$ may not exist, in general, see Example~\ref{prism} below.

\begin{constr}(Mappings of moment-angle manifolds I)\label{mapmfds}
Let $F^{r}=F_{i_{1}}\cap\ldots\cap F_{i_{n-r}}$ be an $r$-dimensional face in $P^n$ and $i^{n}_{r}:\,F^{r}\hookrightarrow\partial P^{n}\subset P^n$ be its affine embedding into $P^n$. Then $F^r$ is itself a convex polytope with $m(F)$ facets $G_{i},1\leq i\leq m(F)$. Therefore, each facet $G_{\alpha}$ of $F^r$ can be uniquely represented in the form 
$$
G_{\alpha}=(F_{i_{1}}\cap\ldots\cap F_{i_{n-r}})\cap F_{j}
$$
for some facet $F_{j}$ of $P^n$.

Then a map $\phi^{n}_{r}:\,[m(F)]\rightarrow [m(P)]$ is defined such that $\phi^{n}_{r}(\alpha)=j$.

Now, applying Construction~\ref{mamfdDJ}, we are going to construct a map $\hat{\phi}^{n}_{r}:\,\mathcal Z_{F^r}\rightarrow\mathcal Z_{P^n}$, induced by the mappings $\phi^{n}_{r}$ and $i^{n}_{r}$.

First, let us consider the following map
$$
\tilde{\phi}^{n}_{r}:\,T^{\mathcal F(F^r)}\rightarrow T^{\mathcal F(P^n)},
$$
where $\mathcal F(F^r)=\{G_{1},\ldots,G_{m(F)}\}$ and $\mathcal F(P^n)=\{F_{1},\ldots,F_{m(P)}\}$ denote the sets of facets of the polytopes $F^r$ and $P^n$, respectively, and
$$
\tilde{\phi}^{n}_{r}(t_{1},\ldots,t_{m(F)})=(\tau_{1},\ldots,\tau_{m(P)}),
$$
where
$$
\tau_{i}=\begin{cases}
t_{(\phi^{n}_{r})^{-1}(i)},&\text{if $i\in\im\phi_{r}^{n}$;}\\
1,&\text{otherwise.}
\end{cases}
$$

It is easy to see that $\tilde{\phi}^{n}_{r}:\,T^{m(F)}\rightarrow T^{m(P)}$ is a group homomorphism.

Finally, we are able to define a map
$$
\hat{\phi}^{n}_{r}:\,\mathcal Z_{F^r}=(T^{\mathcal F(F^r)}\times F^{r})/\sim\rightarrow\mathcal Z_{P^n}=(T^{\mathcal F(P^n)}\times P^{n})/\sim
$$
by formula:
$$
\hat{\phi}^{n}_{r}([t,p])=[\tilde{\phi}^{n}_{r}(t),i_{r}^{n}(p)].
$$

Let us prove that the above definition is correct. Due to Construction~\ref{mamfdDJ}, it suffices to prove that if $T_{1}^{-1}T_{2}\in T^{G_{F^r}(p=q)}$, then $\tilde{\phi}^{n}_{r}(T_{1}^{-1})\tilde{\phi}^{n}_{r}(T_{2})\in T^{G_{P^n}(p=q)}$.
Since $\tilde{\phi}^{n}_{r}$ is a group homomorphism, we need to show that
$$
\tilde{\phi}^{n}_{r}(T_{1}^{-1}T_{2})\in T^{G_{P^n}(p=q)},
$$
when
$$
T_{1}^{-1}T_{2}\in T^{G_{F^r}(p=q)}.
$$

Suppose $T_{1}^{-1}T_{2}=(t_{1},\ldots,t_{m(F)})$ and $\tilde{\phi}^{n}_{r}(T_{1}^{-1}T_{2})=(\tau_{1},\ldots,\tau_{m(P)})$.
Observe that $G:=G_{F^r}(p=q)=G_{P^n}(p=q)$ is a uniquely determied face in $F^r\subset P^n$, for which the point $p=q$ belongs to its interior.

Set $G=G_{\alpha_{1}}\cap\ldots\cap G_{\alpha_{k}}=F^{r}\cap F_{j_{1}}\cap\ldots\cap F_{j_k}$. Then, by definition of $\phi^{n}_{r}$, we have:
$\phi^{n}_{r}(\alpha_{p})=j_p$ for any $1\leq p\leq k$.

Now it suffices to show that if $\tau_{i}\neq 1$, then $i\in\{j_{1},\ldots,j_{k}\}$.
Since $\tau_{i}\neq 1$, then, due to definition of $\tilde{\phi}^{n}_{r}$, we have: $\tau_{i}=t_{(\phi^{n}_{r})^{-1}(i)}\neq 1$ and $i\in\im\phi^{n}_{r}$.
Since $t_{(\phi^{n}_{r})^{-1}(i)}\neq 1$ and $(t_{1},\ldots,t_{m(F)})\in T^{G}=T^{G_{\alpha_{1}}}\times\ldots\times T^{G_{\alpha_{k}}}$, we obtain: $t_{(\phi^{n}_{r})^{-1}(i)}\in T^{G_{\alpha_{s}}}$ for some $s\in [k]$. Again, by definition of the map $\tilde{\phi}^{n}_{r}$, the latter means that $(\phi^{n}_{r})^{-1}(i)=\alpha_{s}$, or equivalently, $\phi^{n}_{r}(\alpha_{s})=i$. As we previously had $\phi^{n}_{r}(\alpha_{s})=j_{s}$ for the $G$ defined above, it follows that $i=\phi^{n}_{r}(\alpha_{s})=j_{s}\in\{j_{1},\ldots,j_k\}$, and the proof of the definition of the map $\hat{\phi}^{n}_{r}$ be correct is done.

To prove the map $\hat{\phi}^{n}_{r}$ is continuous, observe that, by definition of $\hat{\phi}^{n}_{r}$, the next commutative diagram holds:
$$\begin{CD}
  T^{m(F)}\times F^{r} @>\tilde{\phi}^{n}_{r}\times i_{r}^{n}>> T^{m(P)}\times P^{n}\\
  @VV\pr_{r} V\hspace{-0.2em} @VV\pr_{n} V @.\\
  \mathcal Z_{F^r} @>\hat{\phi}^{n}_{r}>>\mathcal Z_{P^n}
\end{CD}\eqno 
$$
Since $\mathcal Z_{F^r}$ and $\mathcal Z_{P^n}$ both have quotient topologies, determied by the canonical projection maps $pr_{r}$ and $pr_{n}$, respectively (cf. Construction~\ref{mamfdDJ}), we deduce that the map $\hat{\phi}^{n}_{r}$ is continuous provided the composition map $\hat{\phi}^{n}_{r}\pr_{r}=\pr_{n}\;(\tilde{\phi}^{n}_{r}\times i_{r}^{n})$ is continuous. The latter mapping is continuous as a composition of continuous maps, which finishes the proof.
\end{constr}

We are going now to construct a special simplicial map $\Phi^{n}_{r}:\,K_{F^r}\rightarrow K_{P^n}$, determined by the function $\phi_{r}^{n}$ such that it induces a continuous map of moment-angle complexes $\hat{\Phi}^{n}_{r}:\,\mathcal Z_{K_{F^r}}\rightarrow\mathcal Z_{K_{P^n}}$.

\begin{constr}\label{mapmcxs}
First, observe that $\phi^{n}_{r}$ induces an injective map of sets of facets $\overline{\phi}^{n}_{r}:\,\mathcal F(F^r)\rightarrow\mathcal F(P^n)$, where $\overline{\phi}^{n}_{r}(G_{\alpha})=F_{\phi^{n}_{r}(\alpha)}$ and $G_{\alpha}=F^r\cap F_{\phi^{n}_{r}(\alpha)}$.

Then the simplex $\sigma=(\alpha_{1},\ldots,\alpha_k)\in K_{F^r}$ is in one-to-one correspondence with a nonempty intersection of facets of a face of $F^r$:
$$
G_{\alpha_1}\cap\ldots\cap G_{\alpha_k}=F^{r}\cap F_{\phi^{n}_{r}(\alpha_{1})}\cap\ldots\cap F_{\phi^{n}_{r}(\alpha_k)}\neq\varnothing.
$$
It follows that $F_{\phi^{n}_{r}(\alpha_{1})}\cap\ldots\cap F_{\phi^{n}_{r}(\alpha_k)}\neq\varnothing$, that is, we are able to define a nondegenerate injective simplicial mapping $\Phi^{n}_{r}:\,K_{F^r}\rightarrow K_{P^n}$ by formula
$$
\Phi^{n}_{r}(\sigma)=(\phi_{r}^{n}(\alpha_{1}),\ldots,\phi^{n}_{r}(\alpha_k))\in K_{P^n}.
$$
Note that:
\begin{itemize}
\item[(1)] If $F^{r}=F_{i_{1}}\cap\ldots\cap F_{i_{n-r}}$ in $P$, then there is a simplex $\Delta({F})=\{i_{1},\ldots,i_{n-r}\}\in K_P$;
\item[(2)] By definition of link, $\Phi_{r}^{n}(K_{F^r})=\Link_{K_{P^n}}\Delta(F)$. 
\end{itemize} 

It follows directly from the definition of a moment-angle-complex that $\Phi^{n}_{r}$ determines a continuous map of moment-angle-complexes:
$$
\hat{\Phi}^{n}_{r}:\,\mathcal Z_{K_{F^r}}\rightarrow\mathcal Z_{K_{P^n}},
$$
induced by homeomorphisms:
$$
(D^2,S^1)^{\sigma}\cong (D^2,S^1)^{\Phi^{n}_{r}(\sigma)}.
$$
\end{constr}

\begin{remark}
Note that both mappings $\hat{\phi}^{n}_{r}$ and $\hat{\Phi}^{n}_{r}$ are weakly equivariant with respect to the compact torus $\mathbb{T}^{m(F)}$-action, induced by the map $\tilde{\phi}^{n}_{r}:\,\mathbb{T}^{m(F)}\rightarrow\mathbb{T}^{m(P)}$.
\end{remark}

Now, applying Construction~\ref{mamfdBP}, we obtain a description of a map between moment-angle manifolds $\mathcal Z_{F^r}\rightarrow\mathcal Z_{P^n}$, induced by the face embedding $i_{r}^{n}:\,F^r\to P^n$, which is equivalent to that considered in Construction~\ref{mapmfds}.

\begin{constr}(Mappings of moment-angle manifolds II)\label{mapmfds2}
Let $F^{r}$ be an $r$-dimensional face in $P^n$ with the set of facets $\mathcal F(P^n)=\{F_{1},\ldots,F_{m(P)}\}$.

Let us show that there exists an induced embedding $\hat{i}_{r}^{n}:\,\mathcal Z_{F^r}\hookrightarrow\mathcal Z_{P^n}$.
Consider an affine embedding $f:\,\R^n \to\R^{m(P)}$ such that its restriction to the polytope $P^n$ coincides with the map $i_{P^n}$, see Construction~\ref{mamfdBP}, and an embedding $g:\,\R^r \to\R^n$ with its restriction on the face $F^r$ coinciding with the map $i_{r}^{n}$.
Now consider composition of the mappings defined above $f\,g:\,\R^r \to \R^n \to R^{m(P)}$ and a continuous section $s_{n}: \R^{m(P)}_\ge \to \C^{m(P)}$.
We remind here, that a tori embedding $\tilde{\phi}_{r}^{n}:\,T^{m(F)}\to T^{m(P)}$, described in Construction~\ref{mapmfds}, gives rise to a group $T^{m(F)}$ action on the complex Euclidean space $\C^{m(P)}$. Thus, we get an action of the compact torus $T^{m(F)}$ on the image (under the composition map) $s_{n}\,f\,g(F^r)$ of the face $F^r$ in $\C^{m(P)}$, where $s_{n}$ is a continuous section of the moment map $\mu_{n}:\,\C^{m(P)}\to\R^{m(P)}_{\ge}$. Denote the corresponding moment-angle manifold (cf. Constrcution~\ref{mamfdDJ}) by $W$. Then, due to Construction~\ref{mamfdBP}, $W$ is embedded into $\mathcal Z_{P^n}\hookrightarrow\C^{m(P)}$ and an induced mapping  $\hat{i}_{r}^{n}:\,\mathcal Z_{F^r}\to\mathcal Z_{P^n}$ is defined. Its image coincides with $W$.

Explicit formulae for the induced embedding $\hat{i}^{n}_{r}$ described above were given in~\cite[Construction 2.10]{BL1}. 
\end{constr}

\begin{remark}
Obviously, one has real analogues of the Constructions 2.11, 2.12, and 2.13 for an induced embedding of real moment-angle manifolds: $\hat{i}_{r}^{n}:\,\mathcal R_{F^r}\to\mathcal R_{P^n}$ and real moment-angle-complexes: $\hat{\Phi}^{n}_{r}:\,\mathcal R_{K_{F^r}}\rightarrow\mathcal R_{K_{P^n}}$.
\end{remark}


In terms and notations from the above constructions let us introduce the following simplicial complexes: $K^{r-1}=K_{F^r}$, $K^{n-1}=K_{P^n}$, and $K_{n,r}=K^{n-1}_{\phi_{r}^{n}[m(F)]}$. Then, due to Construction~\ref{mapmcxs}, a nondegenerate injective simplicial mapping $\Phi_{r}^{n}:\,K^{r-1}\to K^{n-1}$ gives rise to an embedding
$$
\Phi^{n}_{r}(K^{r-1})\subseteq K_{n,r}
$$
and is a composition of simplicial maps:
$K^{r-1}\rightarrow K_{n,r}\rightarrow K^{n-1}$, where the first map is induced by $\Phi^{n}_{r}$, and the second map is a natural embedding of a full subcomplex into its ambient simplicial complex.

The next general result holds.

\begin{propos}\label{GenCase}
There exists a commutative diagram
$$
\xymatrix{
& \mathcal Z_F\ar[rr]^{\hat{i}^{n}_{r}} \ar[dl]_{h_F} && \mathcal Z_P\ar[dr]^{h_P} \\
\mathcal Z_{K^{r-1}} \ar[drr] &&&& \mathcal Z_{K^{n-1}},\\
&& \mathcal Z_{K_{n,r}} \ar[urr]
}
$$
where $h_{F}$ and $h_{P}$ are canonical equivariant homeomorphisms, and the composition in the lower row of arrows equals $\hat{\Phi}^{n}_{r}$. A similar statement holds for induced mappings of real moment-angle manifolds and complexes.
\end{propos}
\begin{proof}
The above diagram commutes. For, it suffices to apply Construction~\ref{mapmfds} to the equivariant homeomorphism $h_{P}:\,\mathcal Z_{P}\cong\mathcal Z_{K_{P}}$ constructed in~\cite[Lemma 3.1.6, formula (37)]{bu-pa00-2}, using the cubical subdivision $\mathcal C(P)\subset I^{m(P)}$. Indeed, we first observe that for any face $F^{r}\subset P^n$ the corresponding cubical subdivision $\mathcal C(F^r)$ is a cubical subcomplex in $\mathcal C(P^n)$, and the next diagram commutes
$$
\begin{CD}
  \mathcal Z_{F^r} @>\hat{i}_{r}^{n}>> \mathcal Z_{P^n} @>h_{P}>> \mathcal Z_{K^{n-1}} @>>> (\mathbb{D}^2)^m\\
  @VVr V\hspace{-0.2em} @VV\rho V\hspace{-0.2em} @VV\rho V\hspace{-0.2em} @VV\rho V @.\\
  F^{r} @>i_{r}^{n}>>P^n @>j_P>> \cc(K^{n-1}) @>>> I^{m(P)},
\end{CD}\eqno 
$$
where $r$ denotes a projection map onto the orbit space of the canonical  $\mathbb{T}^{m(F)}$-action, $\rho$ denotes the projection map onto the orbit space of the canonical $\mathbb{T}^{m(P)}$-action, and $j_P$ is an embedding of a cubical subdivision $\mathcal C(P)$ of the polytope $P$ into $I^{m(P)}$, with its image being $\cc(K^{n-1})$, see~\cite{bu-pa00-2}.

Note that the composition of mappings in the lower row of arrows is equal to $j_F$, since $j_{P}(\mathcal C(F))=\cc(K^{r-1})\subset\cc(K^{n-1})$. Therefore, the composition in the upper row of arrows equals $\hat{\Phi}^{n}_{r}h_{F}$, which finishes the proof.
\end{proof}

In fact, we are able to prove the next statement, which is more general than the above one in the case of a flag polytope.

\begin{lemma}\label{FlagCommuteEmbed}
Let $P^n$ be a flag polytope and $F^r\subset P^n$ be its $r$-dimensional face. Then $\Phi^{n}_{r}(K^{r-1})=K_{n,r}$.
Furthermore, the following diagram commutes
$$\begin{CD}
  \mathcal Z_{F^r} @>\hat{i}_{r}^{n}>> \mathcal Z_{P^n}\\
  @VVH_{1} V\hspace{-0.2em} @VVH_{2} V @.\\
  \mathcal Z_{K^{r-1}} @>\hat{\Phi}_{r}^{n}>>\mathcal Z_{K^{n-1}},
\end{CD}\eqno 
$$
where $H_{1}, H_{2}$ are homeomorohisms and $\hat{i}_{r}^{n}$ induces a split epimorphism in cohomology. A similar result holds for the induced mappings of real moment-angle manifolds and complexes.
\end{lemma}
\begin{proof}
Let us prove the first part of the statemet. By Construction~\ref{mapmcxs}, for any simple polytope $P^n$ one has the next embedding of simplicial complexes, in which the both are subcomplexes on the vertex set $\phi^{n}_{r}[m(F)]$ of the simplicial complex $K^{n-1}$:
$$
\Phi^{n}_{r}(K^{r-1})\subseteq K_{n,r}.
$$
We must prove that the opposite inclusion holds. Assuming the converse, suppose that a set of indices $\{j_{1},\ldots,j_k\}\subset\phi_{r}^{n}[m(F)]$ is such that:
$$
\sigma=(j_{1},\ldots,j_k)\in K_{n,r},\,\sigma\notin\Phi^{n}_{r}(K^{r-1}).
$$
These formulae are equivalent to the following relations in the face poset of the polytope $P^n$:
$$
F_{j_{1}}\cap\ldots\cap F_{j_k}\neq\varnothing,\,F^{r}\cap F_{j_1}\cap\ldots\cap F_{j_k}=\varnothing.
$$
Thus, $F_{j_s}\cap F_{j_t}\neq\varnothing$ for any pair of indices $s,t\in [k]$, and
$$
F_{i_1}\cap\ldots\cap F_{i_{n-r}}\cap F_{j_1}\cap\ldots\cap F_{j_k}=\varnothing.
$$
On the other hand, by definition of the map $\phi_{r}^{n}$, we have the next relations in the face poset of the polytope $F^r$
$$
F_{j_s}\cap F_{i_{1}}\cap\ldots\cap F_{i_{n-r}}=G_{\alpha_s},
$$
in which $\phi_{r}^{n}(\alpha_{s})=j_s$ for $s\in [k]$. It follows that $F_{j_s}\cap F_{i_t}\neq\varnothing$ for any indices $s\in [k], t\in [n-r]$.
Moreover, $F_{i_q}\cap F_{i_p}\neq\varnothing$ for any pair of indices $p,q\in [n-r]$, since $F_{i_{1}}\cap\ldots\cap F_{i_{n-r}}=F^r\neq\varnothing$.

Now recall that the polytope $P^n$ is flag, and we have already proved that each pair of its facets from the set
$$
F_{i_1},\ldots,F_{i_{n-r}},F_{j_1},\ldots,F_{j_k}
$$
has a nonempty intersection. It follows that
$$
F_{i_1}\cap\ldots\cap F_{i_{n-r}}\cap F_{j_1}\cap\ldots\cap F_{j_k}\neq\varnothing,
$$
and we got a contradiction.

To prove that second part of the statement we first note that the diagram
$$
\xymatrix{
& \mathcal Z_F\ar[rr]^{\hat{i}^{n}_{r}} \ar[dl]_{h_F} && \mathcal Z_P\ar[dr]^{h_P} \\
Z_{K^{r-1}} \ar[drr] &&&& \mathcal Z_{K^{n-1}}\ar[dll]_{r_{\phi_{r}^{n}[m(F)]}}\\
&& Z_{K_{n,r}}
}
$$
commutes, according to Proposition~\ref{GenCase}, where $r_{\phi_{r}^{n}[m(F)]}$ is a retraction for the induced embedding $\hat{j}_{\phi^{n}_{r}[m(F)]}$ of moment-angle-complexes. Observe that we have the following equality: $\Phi_{r}^{n}(K^{r-1})=K_{n,r}$, which we proved above. Therefore, in the flag case, $\mathcal Z_F$ is always a retract of $\mathcal Z_P$. The rest of the statement follows from the fact that a homomorphism, equivalent to a spit epimorphism, is a split epimorphism itself, cf. Corollary~\ref{splitepi}.
\end{proof}

\begin{propos}\label{FaceInduceSplit}
Let $F^r$ be an $r$-dimensional face of a polytope $P^n$. Then the next conditions are equivalent:
\begin{itemize}
\item[(1)] For any set $\{j_{1},\ldots,j_{k}\}\subset\phi_{r}^{n}[m(F)]$, if $F_{j_{1}}\cap\ldots\cap F_{j_k}\neq\varnothing$, then $F^{r}\cap F_{j_1}\cap\ldots\cap F_{j_k}\neq\varnothing$;
\item[(2)] $\Phi_{r}^{n}(K^{r-1})=\Link_{K^{n-1}}\Delta(F^r)=K_{n,r}$;
\item[(3)] $\hat{i}_{r}^{n}:\,\mathcal Z_{F^r}\rightarrow\mathcal Z_{P^n}$ is an embedding of a submanifold, having a retraction;
\item[(4)] $(i_{r}^{n})^{*}:\,H^*(\mathcal Z_{P^n})\rightarrow H^*(\mathcal Z_{F^r})$ is a split epimorphism of cohomology rings.
\end{itemize}
\end{propos}
\begin{proof}
Conditions (1) and (2) are equivalent by Construction~\ref{mapmcxs} and the proof of Lemma~\ref{FlagCommuteEmbed}.
The implications $(2)\Rightarrow (3)\Rightarrow (4)$ follow directly from Corollary~\ref{splitepi}. Finally, (4) implies (2), since, if $\Phi_{r}^{n}(K^{r-1})$ is not a full subcomplex, then there exists a simplex $\sigma=(j_{1},\ldots,j_k)\in K_{n,r}$, $\sigma\notin\Phi_{r}^{n}(K^{r-1})$, and $\partial\sigma\subset\Phi_{r}^{n}(K^{r-1})$, which, in view of (4) and Theorem~\ref{mgrad}, gives: $H^*(\partial\sigma)$ is a direct summand in $H^*(\mathcal Z_{K^{n-1}})$ (and, in particular, $\beta^{-i,2J}(K^{n-1})>0$, where $J=(j_{1},\ldots,j_{m})$ is an $m$-tuple of 0s and 1s such that $j_{t}=1$ if and only if $t\in\sigma$; $|J|-i-1=k-2$), which is a contradiction. 
\end{proof}

Therefore, when a polytope $P^n$ is flag, the mapping $\Phi_{r}^{n}$ from Construction~\ref{mapmcxs} provides a simplicial isomorphism between $K_{F^r}$ and the full subcomplex in $K_{P^n}$ on the vertex set $\phi_{r}^{n}[m(F)]$. So, in what follows, in the flag case, we shall identify those simplicial complexes and consider the corresponding simplicial embedding $j_{r}^{n}:\,K_{F^r}\rightarrow K_{P^n}$.

The next example clarifies the situation in the nonflag case.

\begin{exam}\label{prism}
Consider a triangular prism $P^3$, $m(P^3)=5$, and denote its triangular facets by $F_1$ and $F_5$. Consider its quadrangular facet $F_2$. Then $\phi_{2}^{3}[m(F_2)]=\{1,3,4,5\}\subset [m(P)]=[5]$, a nerve complex of the face $F_2$ is a boundary of a square, and the full subcomplex in $K_{P^3}$ on the vertex set $\{1,3,4,5\}$ is a boundary of a square alongside with the edge $\{3,4\}$. This corresponds to the well known fact that there exists no retraction $S^{3}\times S^{5}$ onto $S^{3}\times S^3$. On the other hand, for any triangular facet, in particular, for $F_{1}$, we have: $\Phi_{2}^{3}(K_{F_1})=K_{3,2}$; there exists a retraction $S^{3}\times S^{5}\to S^{5}$, and the corresponding split epimorphism in cohomology.
\end{exam}

\begin{remark}
For any face embedding $i_{r}^{n}:\,F^{r}\to P^n$ there exists an induced embedding of quasitoric manifolds $M^{2r}(F,\Lambda_{F})\to M^{2n}(P,\Lambda)$, which can be realized as a composition map of the embeddings of characteristic submanifolds, induced by face emebeddings:
$$
F^{r}=F_{i_1}\cap\ldots\cap F_{i_{n-r}}\subset F_{i_2}\cap\ldots\cap F_{i_{n-r}}\subset\ldots\subset F_{i_{n-r}}\subset P.
$$
An alternative description of the induced mapping of charactersitic submanifolds can be obtained similar to that determined in Construction~\ref{mapmfds}, using the definition of a quasitoric manifold $M(P,\Lambda)$ as a quotient space of a moment-angle manifold $\zp$ by a freely acting subtorus of dimension $(m-n)$.

However, the induced mappings of quasitoric manifolds, in general, do not have retractions, even for a flag polytope $P^n$, as the next example shows.
\end{remark}

\begin{exam}
Consider a pentagon $P_5^2$, embedded into $\mathbb{R}^2$ as shown in the figure below, with facets (edges) $F_1,\ldots,F_5$, denoted here simply by $\{1,2,3,4,5\}$.


\begin{center}
\begin{picture}(90,60)(0,0)
{\thicklines
  \put(25,40){\line(1,0){20}}
  \put(65,0){\line(0,1){20}}
  \put(65.1,0){\line(0,1){20}}
  \put(65,20){\line(-1,1){20}}
  }
   \put(20,0){\vector(1,0){60}}
   \put(25,-5){\vector(0,1){60}}
  \put(77,-5){$x_1$}
  \put(18,53){$x_2$}
   \put(20,21){$1$}
   \put(45,-5){$2$}
   \put(67,10){$3$}
   \put(58,29){$4$}
   \put(35,41){$5$}
\end{picture}
\end{center}
\vskip 1.0cm

Let us take $P^2=P_5^2 = \{ x\in \mathbb{R}^2\;: Ax+b\geqslant 0 \}$, where
\[ A^\top =
  \begin{pmatrix}
    1 & 0 & -1 & -1 & 0 \\
    0 & 1 & 0 & -1 & -1 \\
  \end{pmatrix},\qquad
 b^\top = (0,0,2,3,2).
\]
and $\top$ denotes matrix transposition operation.

By Davis--Januszkiewicz theorem, cf.~\cite{DJ}, we have the next description of cohomology ring of the quasitoric manifold $M=M_{P_5}$
$$
H^*(M(P_{5}^{2},A^{\top}))\cong\mathbb{Z}[v_{1},v_{2},v_{3},v_{4},v_{5}]/I,
$$
where
$$
I=(v_{1}-v_{3}-v_{4},v_{2}-v_{4}-v_{5},v_{1}v_{3},v_{1}v_{4},v_{2}v_{4},v_{2}v_{5},v_{3}v_{5}),
$$
from which it follows that $v_{1}^{2}=v_{2}^{2}=0$ and $v_{3}^{2}=v_{4}^{2}=v_{5}^{2}\neq 0$ in $H^*(M)$ (since, if all the squares of 2-dimensional generators equal zero, then all the monomials of degree greater than 2 in $H^*(M)$ are also equal to zero, which contradicts the fact that $H^{4}(M)\cong\mathbb{Z}$). We state that there is no retraction for the embedding of quasitoric manifolds, induced by face embedding $i_{1}^{3}:\,F_{3}\to P_5$. Otherwise one would have existence of a split epimorphism in cohomology: $H^{2}(M)\to H^{2}(M_{F_3})$, which would contradict $v_{3}^{2}\neq 0$ in $H^*(M)$ (recall that $M_{F_3}\cong\mathbb{C}P^1$, and so for its 2-dimensional cohomology generator one has: $v_{3}^{2}=0$).
\end{exam}

Now we would like to show that the opposite implication to that of Lemma~\ref{FlagCommuteEmbed} also holds. To do that, we need the next combinatorial criterion of flagness of a simplicial complex.

\begin{lemma}\label{FlagCriterionLemma}
A simplicial complex $K$ is flag if and only if $\Link_{K}(v)$ is a full subcomplex in $K$ for any vertex $v\in K$.
\end{lemma}
\begin{proof}
If $K$ is a flag complex, then suppose that the converse is true. Then there exists a vertex $v\in K$ and a simplex $\sigma\in K$ such that $|\sigma|\geq 2$, $\partial\sigma\subseteq\Link_{K}(v)$ and $\sigma\notin\Link_{K}(v)$. It implies that $v\notin\sigma$, $v\cup\sigma\notin K$ and $v\cup\sigma_{i}\in K$, for any $i$, where $\sigma_{i}$ denotes a facet of a simplex $\sigma$. The latter means that $v\cup\sigma\in MF(K)$ with $|v\cup\sigma|\geq 3$ elements. We got a contradiction with flagness of $K$.

Suppose now that $K$ is not flag. Then there exists a minimal nonface $\{i_{1},\ldots,i_{p}\}$ of $K$ on $p\geq 3$ vertices. Then $\partial\Delta_{\{i_{1},\ldots,i_k\}}$ is a full subcomplex in $K$ on the vertex set $\{i_{1},\ldots,i_{p}\}$. Observe that $\Link_{K}(i_1)$ contains all vertices from the set $\{i_{2},\ldots,i_k\}$, but does not contain the simplex $(i_{2},\ldots,i_k)\in K$. That is, $\Link_{K}(i_1)$ is not a full subcomplex in $K$ on its set of vertices, which finishes the proof.
\end{proof}

Finally, we are able to get the next result.

\begin{theorem}\label{FlagCriterion}
The following conditions are equivalent:
\begin{itemize}
\item[(1)] $P^n$ is a flag polytope;
\item[(2)] $\Phi_{r}^{n}(K^{r-1})=K_{n,r}$ for any face $F^r\subset P^n$;
\item[(3)] $\hat{i}_{r}^{n}:\,\mathcal Z_{F^r}\rightarrow\mathcal Z_{P^n}$ is a submanifold embedding, having retraction for any face $F^r\subset P^n$;
\item[(4)] $(i_{r}^{n})^{*}:\,H^*(\mathcal Z_{P^n})\rightarrow H^*(\mathcal Z_{F^r})$ is a split epimorphism of cohomology rings, for any face $F^r\subset P^n$.
\end{itemize}
\end{theorem}
\begin{proof}
Conditions (2), (3), and (4) are equivalent to each other by Proposition~\ref{FaceInduceSplit}. 

Suppose condition (2) holds. Then, in particular, for any facet $F^{n-1}\subset P^{n}$ one has: $\Phi_{n-1}^{n}(K^{n-2})=K_{n,n-1}$. On the other hand, $\Phi_{n-1}^{n}(K^{n-2})=\Link_{K_P}(v)$, where $v\in K_P$ corresponds to the facet $F^{n-1}\subset P^n$, cf. Construction~\ref{mapmcxs}. The latter means that $\Link_{K_P}(v)$ is a full subcomplex in $K_P$ for any vertex $v\in K_P$, therefore, by Lemma~\ref{FlagCriterionLemma}, $K_P$ is a flag simplicial complex, which implies condition (1). 

Finally, (1) implies (2) by Lemma~\ref{FlagCommuteEmbed}, which finishes the proof.
\end{proof}


In the final part of this section we shall introduce the definition of a Massey product in cohomology of a differential graded algebra. First of all, we need a notation of a defining system for an ardered set of cohomology classes. In our exposition we follow that in~\cite{Kr} and in~\cite[Appendix $\Gamma$]{BP04}, where one can find all the properties of Massey product, necessary in what follows.

Let $(A,d)$ be a differential graded algebra, $\alpha_{i}=[a_{i}]\in H^{*}[A,d]$ and $a_{i}\in A^{n_{i}}$ for $1\leq i\leq k$.
Then a \emph{defining system} for the set $(\alpha_{1},\ldots,\alpha_{k})$ is a $(k+1)\times (k+1)$-matrix $C$ such that the following conditions hold:
\begin{itemize}
\item[{(1)}] $c_{i,j}=0$ for $i\geq j$,
\item[{(2)}] $c_{i,i+1}=a_{i}$,
\item[{(3)}] $a\cdot E_{1,k+1}=dC-\bar{C}\cdot C$ for a certain element $a=a(C)\in A$, where $\bar{c}_{i,j}=(-1)^{\deg(c_{i,j})}\cdot c_{i,j}$ and $E_{1,k+1}$ denotes a $(k+1)\times(k+1)$-matrix with a unit in the position $(1,k+1)$ and zeros in all other positions.
\end{itemize} 

A direct computation shows that $d(a)=0$ and $a\in A^{m}$, $m=n_{1}+\ldots+n_{k}-k+2$. Thus, $a=a(C)$ is a cocycle for any defining system $C$, and a cohomology class $\alpha=[a]\in H^{m}[A,d]$ is defined correctly.

\begin{definition}\label{defMassey}
We say that a \emph{$k$-fold Massey product} $\langle\alpha_{1},\ldots,\alpha_{k}\rangle$ is \emph{defined}, if for the corresponding ordered set of cohomology classes there exists a defining system $C$. In this case, by definition, the Massey product $\langle\alpha_{1},\ldots,\alpha_{k}\rangle$ contains of all cohomology classes $\alpha=[a(C)]$, for $C$ -- a defining system. A defined Massey product $\langle\alpha_{1},\ldots,\alpha_{k}\rangle$ is said to be 
\begin{itemize}
\item \emph{trivial} (or, \emph{vanishing}), if $[a(C)]=0$ for some defining system $C$;
\item \emph{decomposable}, if $[a(C)]\in H^{+}(A)\cdot H^{+}(A)$ for some defining system $C$;
\item \emph{strictly defined}, if $\langle\alpha_{1},\ldots,\alpha_{k}\rangle=\{[a(C)]\}$ for some (and, therefore, for any) defining system $C$.
\end{itemize}
\end{definition}

\begin{remark}
By~\cite[Theorem 3]{Kr}, the set of cohomology classes $\langle\alpha_{1},\ldots,\alpha_{k}\rangle$ depends only on the cohomology classes $\alpha_{1},\ldots,\alpha_{k}$, rather than on the choice of particular representing cocycles $a_{1},\ldots,a_{k}$. Furthermore, if $[a_{i}]=0$ for some $1\leq i\leq k$ and $\langle\alpha_{1},\ldots,\alpha_{k}\rangle$ is defined, then $0\in\langle\alpha_{1},\ldots,\alpha_{k}\rangle$, and the $k$-fold Massey product is trivial.
\end{remark}

Let us introduce explicit formulae for the relations, see condition (3) above, on the elements of a defining system $C$ in the cases of 2-, 3-, and 4-fold defined Massey product.

\begin{exam}
Suppose $k=2$. If a Massey product $\langle\alpha_{1},\alpha_{2}\rangle$ is defined, then:
$$
a=d(c_{1,3})-\bar{a}_1\cdot a_{2};
$$
Therefore, our definition gives the usual cup-product in cohomology of topological spaces (up to sign).
\end{exam}

\begin{exam} 
Suppose $k=3$. If a Massey product $\langle\alpha_{1},\alpha_{2},\alpha_{3}\rangle$ is defined, then:
$$
a=d(c_{1,4})-\bar{a}_1\cdot c_{2,4}-\bar{c}_{1,3}\cdot a_{3},
$$
\[
\begin{aligned}
d(c_{1,3})&=&\bar{a}_1\cdot a_{2},\\
d(c_{2,4})&=&\bar{a}_2\cdot a_{3};    
\end{aligned}
\]
\end{exam}

\begin{exam} 
Suppose $k=4$.
If a Massey product $\langle\alpha_{1},\alpha_{2},\alpha_{3},\alpha_{4}\rangle$ is defined, then:
$$
a=d(c_{1,5})-\bar{a}_1\cdot c_{2,5}-\bar{c}_{1,3}\cdot c_{3,5}-\bar{c}_{1,4}\cdot a_{4},
$$
\[
\begin{aligned}
d(c_{1,3})&=&\bar{a}_1\cdot a_{2},\\
d(c_{1,4})&=&\bar{a}_1\cdot c_{2,4}+\bar{c}_{1,3}\cdot a_{3},\\
d(c_{2,4})&=&\bar{a}_2\cdot a_{3},\\
d(c_{2,5})&=&\bar{a}_2\cdot c_{3,5}+\bar{c}_{2,4}\cdot a_{4},\\
d(c_{3,5})&=&\bar{a}_3\cdot a_{4}. 
\end{aligned}
\]
\end{exam}

Obviously, each defined $k$-fold Massey product in $H^*(\zp)$ is trivial for $k>2$, when $P=I^n$, and is already trivial for $k>1$, when $P=\Delta^n$. In section 4 we construct a family of polytopes $\mathcal P=\{P^n|\,n\geq 0\}$ such that there exist nontrivial strictly defined Massey products of all orders $2\leq k\leq n$ in $H^*(\mathcal Z_{P^n})$.


\section{The theory of families of polytopes: methods and constructions}

The foundation of the theory of families of simple polytopes was made in the work of Buchstaber~\cite{B}. We start this section with giving a definition of the notion of a family of polytopes, following~\cite{B}. Below we consider combinatorially equivalent polytopes as equal.

\begin{definition}\label{Family}
By a \emph{family of polytopes} we mean a set $\mathcal F$ of combinatorial polytopes such that for any $n\geq 0$ a subset $\mathcal F_n\subset\mathcal F$ of its $n$-dimensional elements is nonempty and finite. The family $\mathcal F$ will be called \emph{linear}, if each subset $\mathcal F_i$ consists of a unique element.
\end{definition}

There is a number of sources, dedicated to the study of remarkable families of polytopes, see~\cite{Post05, PRW06, BEMPP,BE2017}. In~\cite{B} a notion of a bigraded differential ring of simple polytopes $(\mathcal P,d)$ was introduced, where the differential $d$ maps a polytope $P$ into the disjoint union of its facets. We shall give all the necessary details from~\cite{B} and introduce an oriented version of this construction below.

\begin{definition}
A {\emph{bigraded ring of polytopes}} is a free Abelian group $\mathcal P=\oplus\mathcal P_{n}$, where $\mathcal P_{n}=\oplus\mathcal P_{n,k}$. Here, we denoted by $\mathcal P_{n,k}$ a finitely generated free Abelian group with generators corresponding to simple $n$-dimensional polytopes $P^n$ with $m$ facets, and $k=m-n$. In the group $\mathcal P$, apart from addition, which corresponds to disjoint union of polytopes, there is also multiplication, which corresponds to a Cartesian product of polytopes. As a zero element in $\mathcal P$ it is natural to take the empty set $P^{-1}$, and as a unit -- the point $P^0$. The ring $\mathcal P$ just described turns out to be a bigraded one, since, if $a\in\mathcal P_{n,k}$ and $b\in\mathcal P_{n',k'}$, then $ab\in\mathcal P_{n+n',k+k'}$.

Recall, that in algebra, a \emph{differential ring} is defined to be a graded ring with an operator, satsfying Leibniz rule. The ring $\mathcal P$ is a differential ring (with respect to the first grading) with respect to the operator $d:\mathcal P\rightarrow\mathcal P$ that maps a polytope to the disjoint union of its facets, and therefore $d:\,\mathcal P_{n}\to\mathcal P_{n-1}$. Note that $d^{2}(P^n)\neq 0$ for $n>1$. 
\end{definition}

\begin{remark}
One can similarly consider a bigraded ring 
$\mathcal{P}^O=\oplus\mathcal P^{O}_{n}$, where $\mathcal P^{O}_{n}=\oplus\mathcal P^{O}_{n,k}$, generated by the set of all oriented polytopes. An operator $d^O$ on such a ring maps an oriented polytope $P$ to the disjoint union of its facets $\{F_i\}$ with induced orientations. Recall that for an oriented polytope $P$ normal vectors to facets $F_i$ are considered as pointing inwards the polytope (cf. Definition~\ref{Simplepolytopes}) and, therefore, an induced orientation of a facet is defined correctly. 
In this case, $(d^O)^2 = 0$ and 
\[
d^O(P_1\times P_2) = (d^OP_1)\times P_2 \cup (-1)^{\dim P_1}P_1\times d^OP_2.
\] 
\end{remark}

\begin{prob}
Compute the homology $H_{*}[\mathcal{P}^O, d^O]$.
\end{prob}

Throughout this paper we deal only with the differential ring $(\mathcal P,d)$. To avoid consideration of families with trivial combinatorial structure, which are not interesting for us here, we introduce the following notions.

\begin{definition}
Let us call a linear family $\mathcal F$ \emph{n-reducible}, if each polytope $P\in\mathcal F_k$ with $k\geq n$ is a product of polytopes of positive dimensions. If $\mathcal F$ is not $n$-reducible for any $n>1$, then we call such a family $\mathcal F$ an \emph{irreducible} family.
\end{definition}

Among all polytopes $P$ the simple polytopes are characterized by the following formula, which is valid only for them:
$$
F(dP) = \frac{\partial}{\partial t} F(P),
$$
where $F(P) = \alpha^n+f_{n-1}\alpha^{n-1}t+\ldots+f_0t^n$.
Here we denote $n=\dim P$ and $f_k$ is a number of $k$-dimensional faces.

In the theory of simple polytopes, alongside with $F$-polynomials, it is useful to consider $H$-polynomials: $H(s,t)=F(s-t,t)$.
For any simple polytope $P$ the classical Dehn--Sommerville theorem on $f$-vector of a simple polytope $P$ is equivalent to the identity $H(\alpha,t)=H(t,\alpha)$.

The following relation on $H$-polynomials holds in the ring $(\mathcal P,d)$:
$$
H(dP)=\partial H(P),
$$
where $\partial=\frac{\partial}{\partial s} + \frac{\partial}{\partial t}$.
Therefore, $\partial$ is a linear operator, which sends symmetric polynomials to symmetric polynomials. Moreover, $F$ and $H$ induce homomorphisms of differential rings:
$$
F:\,\mathcal P\to\mathbb{Z}[\alpha,t]
$$
and
$$
H:\,\mathcal P\to\mathbb{Z}[s,t],
$$
that send differential $d$ to the operator $\frac{\partial}{\partial t}$ in the case of $F$, and send differential $d$ to the operator $\partial$ in the case of $H$.

Finally, what we get is that the formulae for the values of the differential $d$ on families of polytopes provide us with partial differential equations.
In the series of works~\cite{B, BK, BV2} it was shown that for certain important families of polytopes these differential equations have analytical solutions (see also~\cite[\S1.7, 1.8]{TT}). Note that the work~\cite{B} was motivated by the work of Buchstaber and Koritskaya~\cite{BK}, in which it was proved that the generating series for $H$-polynomials in the case of associahedra $As$ satisfies the classical quasilinear E.Hopf equation.

Among the new notions introduced in our work, the next two are the main ones. We recall that we denote by $m(P^n)$ the number of facets of an $n$-dimensional polytope $P$; when the choice of $P$ is clear from the context, we write $m(n)$, or simply, $m$, instead of $m(P^n)$, to simplify our notation.

\begin{definition}\label{ADFP}
A linear family of polytopes $\mathcal F=\{P^n|\,n\geq 0\}$ is called an \emph{algebraic direct family of polytopes} (in brief, ADFP), if the following condition holds:

For any $r\geq 0$ and $n>r$ there exists an $r$-dimensional face $F=F(r,n)$ of a polytope $P^n$ such that $F$ is combinatorially equivalent to $P^r$, and $\{P^n, i_{r}^{n}\}$ is a direct system of polytopes with respect to face embeddings $i_{r}^{n}:\,P^{r}\hookrightarrow\partial P^n\subset P^n$ (here we identified $F$ with $P^r$).
\end{definition}

\begin{definition}\label{GDFP}
An algebraic direct family of polytopes $\mathcal F=\{P^n|\,n\geq 0\}$
is called a \emph{geometric direct family of polytopes} (in brief, GDFP), if the following condition holds:

For any $r\geq 0$ and $n>r$ there exists a set of vertices $J=J(r,n)\subset [m(n)]=\{1,2,\ldots,m(n)\}$ of cardinality $m(r)$ such that the full subcomplex in $K_{P^n}$ on $J$ is combinatorially equivalent to $K_{P^r}$, and $\{K_{P^n}, j_{r}^{n}\}$ is a direct system of triangulated spheres with respect to simplicial embeddings $j_{r}^{n}:\,K_{P^{r}}\hookrightarrow K_{P^n}$, induced by embeddings $J\subset [m(n)]$ (here we identified $(K_{P^n})_J$ with $K_{P^r}$).
\end{definition}

\begin{exam}
1) The family of simplices $\Delta=\{\Delta^n\,|\,n\geq 0\}$ is an ADFP, but not a GDFP. 2) The family of cubes $\mathcal I=\{I^n\,|\,n\geq 0\}$ is a GDFP.
\end{exam}

Observe that each ADFP $\mathcal F$, consisting of flag polytopes, is a GDFP due to Lemma~\ref{FlagCommuteEmbed}. The condition from the definition of GDFP implies that $\{H^*(\mathcal Z_{P^n}), (i_{r}^{n})^{*}\}$ is an inverse system of rings with respect to split epimorphisms of rings $(i_{r}^{n})^*$, cf. Corollary~\ref{splitepi}. Recall that each map $i_{r}^n$ also induces a map of moment-angle manifolds $\hat{\phi}_{r}^{n}:\,\mathcal Z_{P^r}\rightarrow\mathcal Z_{P^n}$, see Construction~\ref{mapmfds}. 

It is easy to see thaht for any GDFP the induced maps $(\hat{\phi}_{r}^{n})^{*}$ and $(j_{r}^{n})^{*}$ are equivalent homomorphisms.
It follows that $(\hat{\phi}_{r}^{n})^{*}:\,H^*(\mathcal Z_{P^n})\rightarrow H^*(\mathcal Z_{P^r})$ are also split ring epimorphisms. 

One of main goals of this paper is to construct examples of direct families of polytopes (not necessarily flag ones) with rich combinatorial structure. This problem will be solved by means of the theory of nestohedra.

\begin{definition}\label{DFPM}
A direct family of polytopes $\mathcal F=\{P^n|\,n\geq 0\}$ will be called a \emph{direct family with nontrivial Massey products} (in brief, DFNM), if there exist a nontrivial Massey product of order $k$ in $H^*(\mathcal Z_{P^n})$, and $k\to\infty$ when $n\to\infty$.
\end{definition}

Note that if $\mathcal F$ is a GDFP, then the next property follows immediately from Corollary~\ref{splitepi} and Theorem~\ref{zkcoh}. If $\langle\alpha_{1}^{r},\ldots,\alpha_{k}^{r}\rangle$ is a defined Massey product in $H^*(\mathcal Z_{P^r})$, then for any $n>r$ there exists a defined Massey product $\langle\alpha_{1}^{n},\ldots,\alpha_{k}^{n}\rangle$ in $H^*(\mathcal Z_{P^n})$ such that $(\hat{\phi}_{r}^{n})^{*}(\alpha_{t}^{n})=\alpha_{t}^{r}$ for any $1\leq t\leq k$ and
$$
(\hat{\phi}_{r}^{n})^{*}\langle\alpha_{1}^{n},\ldots,\alpha_{k}^{n}\rangle=\langle\alpha_{1}^{r},\ldots,\alpha_{k}^{r}\rangle.
$$

\begin{definition}\label{special}
We call a direct family $\mathcal F=\{P^n|\,n\geq 0\}$ of polytopes with nontrivial Massey products a {\emph{special family}}, if for any $n\geq 2$ there exists a strictly defined nontrivial Massey product $\langle\alpha_{1},\ldots,\alpha_{k}\rangle$ in $H^*(\mathcal Z_{P^n})$ for each $2\leq k\leq n$.
\end{definition}

Our main result, see Theorem~\ref{DirFamNontrivial}, states that the family $\mathcal P_{Mas}$, mentioned in the Introduction, is a special geometric direct family with nontrivial Massey products, with an additional condition of $\dim\alpha_{i}=3$ taking place in Definition~\ref{special}. 


To introduce and study the family $\mathcal P_{Mas}$, we must turn now to certain important definitions and results from the theory of nestohedra. Following the works by Feichtner and Sturmfels~\cite{FS} and Postnikov~\cite{Post05}, we now introduce the notions of a building set and a nestohedron, which play fundamental role in our work.

\begin{definition}\label{nest}
A set $B\subseteq 2^{[n+1]}$ is called a {\emph{building set}} on the set of vertices $[n+1]$, if the following two conditions hold:\\
1) For any $i\in [n+1]$: $\{i\}\in B$;\\
2) If $S_{1}, S_{2}\in B$ and $S_{1}\cap S_{2}\neq\varnothing$, then $S_{1}\cup S_{2}\in B$.

A building set $B$ is called {\emph{connected}}, if $[n+1]\in B$.
A \emph{nestohedron} $P_B$ on a (connected) building set $B$ on $[n+1]$ is a  $n$-dimensional simple convex polytope, defined as a Minkowski sum of simplices in the Euclidean space $\mathbb{R}^{n+1}$ with a standard basis $\{e_{1},\ldots,e_{n+1}\}$:
$$
P_{B}=\sum\limits_{S\in B\backslash [n+1]}\,\conv(e_{i}|\,i\in S).
$$ 
\end{definition}

The nest statement is well known and is a direct consequence of the above definition.

\begin{coro}\label{BuildSetProd}
Let $B$ be a union of building sets $B_{1}$ on $[n_{1}+1]$ and $B_{2}$ on $[n_{2}+1]$, $[n_{1}+1]\cap [n_{2}+1]=\varnothing$ (such a $B$ is, ibviously, disconnected). Then $P_{B}=P_{B_1}\times P_{B_2}$. 
\end{coro}

\begin{exam}
1) Consider the set $B_{\Delta}=\{\{i\},[n+1]|1\leq i\leq n+1\}$. Then $P_{B_\Delta}=\Delta^n$;
2) Consider the set $B_{\Box}=\{\{1,\ldots,i\},\{i\}|1\leq i\leq n+1\}$. Then $P_{B_{\Box}}=I^n$. Due to Corollary~\ref{BuildSetProd}, $\hat{B}=\{\{i\},\{2j-1,2j\}|\,1\leq i\leq n+1,1\leq j\leq n\}$ is also a building set for $I^n$. 
\end{exam}

\begin{remark}
The above mentioned examples show that there might be different building sets $B_{1}$ and $B_{2}$ such that $P_{B_1}=P_{B_2}$. 
\end{remark}

In~\cite{PRW06} it was shown that if $P_B$ is a flag nestohedron, then there exists a building set $B_{0}\subseteq B$ such that $P_{B_0}$ is a combinatorial cube and $\dim P_{B_0}=\dim P_{B}$. Thus, a class of polytopes, arising from a cube by face truncations, contains all flag nestohedra and, in particular, all graph-associahedra. A stronger result was obtained in~\cite{BV}, that is: a nestohedron is flag of and only if it can be realized as a \emph{2-truncated cube}, i.e., a polytope, obtained  from a cube by a sequence of truncations of faces of codimension two. Observe that the resulting polytope of a face truncation applied to a cube is flag if and only if the truncated face has codimension two.

The next family of polytopes, introduced by Carr and Devadoss~\cite{CD} in a paper devoted to the study of Coxeter complexes, consists of flag nestohedra, and therefore, due to the above mentioned result of Buchstaber and Volodin~\cite{BV}, these polytopes can be realized as 2-truncated cubes.  

\begin{definition}\label{GA}
A \textit{graphical building set} $B(\Gamma)$ of a simple graph $\Gamma$ on the vertex set $[n+1]$ consists of such subsets $S$ that the induced subgraph $\Gamma_{S}$ on the vertex set $S\subset [n+1]$ is connected.\\
Then the nestohedron $P_{\Gamma}=P_{B(\Gamma)}$ is called a \emph{graph-associahedron}.
\end{definition}

\begin{exam}
The next families of graph-associahedra, playing an important role in our paper, found numerous applications in convex geometry, combinatorics, and representation theory: 
\begin{itemize}
\item $\Gamma$ is a complete graph on $[n+1]$.\\
Then $P_{\Gamma}=Pe^n$ is a \textit{permutohedron}, see Figure~\ref{permfig}.
\begin{figure}[h]
\includegraphics[scale=0.5]{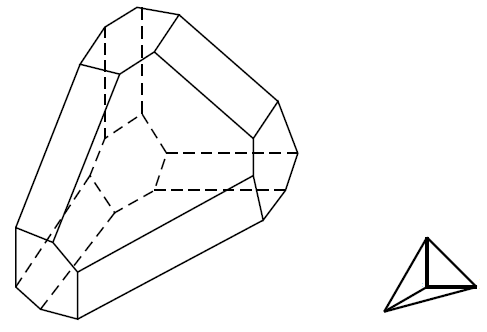}
\caption{3-dimensional permutohedron and the corresponding graph.}\label{permfig}
\end{figure}

\item $\Gamma$ is a star graph on $[n+1]$.\\
Then $P_{\Gamma}=St^n$ is a \textit{stellahedron}, see Figure~\ref{stelfig}.
\begin{figure}[h]
\includegraphics[scale=0.5]{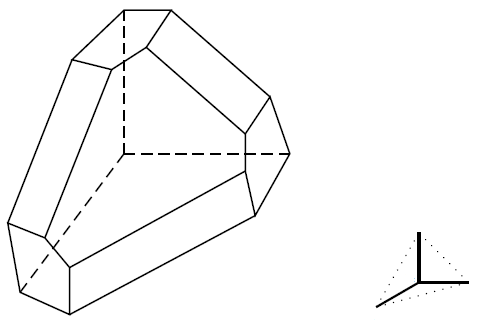}
\caption{3-dimensional stellahedron and the corresponding graph.}\label{stelfig}
\end{figure}

\item $\Gamma$ is a cycle on $[n+1]$.\\
Then $P_{\Gamma}=Cy^n$ is a \textit{cyclohedron} (or, Bott--Taubes polytope~\cite{BT}), see Figure~\ref{cyclfig}.
\begin{figure}[h]
\includegraphics[scale=0.5]{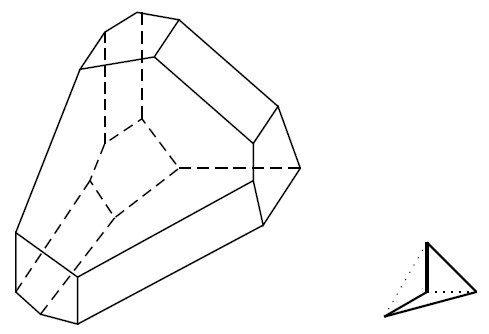}
\caption{3-dimensional cyclohedron and the corresponding graph.}\label{cyclfig}
\end{figure}

\item $\Gamma$ is a chain graph on $[n+1]$.\\
Then $P_{\Gamma}=As^n$ is an \textit{associahedron} (or, Stasheff polytope~\cite{S}), see Figure~\ref{assfig}.
\begin{figure}[h]
\includegraphics[scale=0.5]{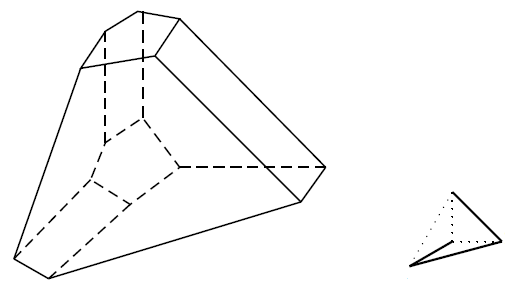}
\caption{3-dimensional associahedron and the corresponding graph.}\label{assfig}
\end{figure}

\end{itemize}
\end{exam}

\begin{exam}
The families of simplices $\Delta=\{\Delta^n|\,n\geq 0\}$, cubes $\mathcal I=\{I^n|\,n\geq 0\}$, associahedra $As=\{As^n|\,n\geq 0\}$, permutohedra $Pe=\{Pe^n|\,n\geq 0\}$, cyclohedra $Cy=\{Cy^n|\,n\geq 0\}$, and stellahedra $St=\{St^n|\,n\geq 0\}$ are linear families. Apart from the cube family, which is 2-reducible, all other families are irreducible.
\end{exam}


Now, suppose that $P=P_B$ is a nestohedron on a connected building set $B\subseteq 2^{[n+1]}$, $n\geq 2$ (cf. Definition~\ref{nest}). 

The next result takes place, which will be particularly important in what follows.

\begin{lemma}[\cite{FS}]\label{NestBound}
The operator $d$ acts on nestohedra by formula:
$$
dP_{B}^n=\sum\limits_{S\in B\backslash [n+1]}\,P_{B|_{S}}\times P_{B/S},\eqno (2)
$$
where \emph{restriction} of a building set $B|_{S}$ is defined to be a building set, induced on the subset of vertices $S\in B$, and a \emph{contraction} of $S$ in $B$ is defined to be a building set $B/S=\{T\subset [n+1]\backslash S\,|\,T\in B\text{ or }T\sqcup S\in B\}$.
\end{lemma}

\begin{definition}
Let $\mathcal F$ be a family of polytopes. Then set $\mathcal P_\mathcal F=\oplus_{n\geq 0}\mathcal {P}_{F,n}$ to be graded subring of the differential (bi)graded ring $\mathcal P$, generated over $\mathbb Z$ by all $P^n\in F_{n}$. Then $\mathcal F$ will be called a {\emph{differential family}} (or, in brief, a {\emph{$d$-family}}), if $\mathcal P_{\mathcal F}$ is itself a differential ring, that is, if $d:\,\mathcal P_{F,n}\rightarrow\mathcal P_{F,n-1}$.
A linear family which is a $d$-family will be called a {\emph{d-linear family}}.
\end{definition}

The families of simplices, cubes, associahedra, permutohedra are $d$-families. In fact, they are $d$-linear families. However, in general, a (linear) family is not necessarily a $d$-family. This remark motivates us to introduce the next notion.

\begin{definition}
By a \emph{$d$-closure} of a family $\mathcal F$ we mean a minimal (with respect to inclusion) extension of a family $\mathcal F$ in the ring $\mathcal P$ to a  $d$-family, if such extension exists.
\end{definition}

The next statement follows immediately from Lemma~\ref{NestBound}.

\begin{coro}\label{ClosurePolytopes}
A $d$-closure of any family of nestohedra contains all nestohedra of the type $P_{B_{n}|_S}$ and $P_{B_{n}/S}$, alongside with every nestohedron $P_{B_n}$, for all elements $S\in B_{n}\backslash [n+1]$. 
\end{coro}

\begin{exam}
It is well known that, cf.~\cite[\S1.6]{TT}, a $d$-closure of the family $Cy$ of cyclohedra equals $As\cup Cy$, and a $d$-closure of the family $St$ of stellahedra equals $St\cup Pe$.
\end{exam}

The above mentioned examples allow us to give the next definition.

\begin{definition}\label{complex}
We say that a family $\mathcal F$ has {\emph{complexity}} $k\geq 1$, if its $d$-closure exists and can be represented as a union of $k$ linear systems.
\end{definition}

In particular, the above examples show taht the families of associahedra and permutohedra have complexity 1 (in other words, they are $d$-linear families themselves), and the families of cyclohedra and stellahedra both have complexity 2.

It can be easily seen from Definition~\ref{nest} that the combinatorial type of an $n$-dimensional nestohedron $P_{B_n}$ is determined by the structure of its building set $B_n$ on $[n+1]$. Thus, any linear family $\mathcal F$ of nestohedra is given by a sequence of building sets $\{B_n, n\geq 1\}$.

The following question is among the important problems of combinatorics of families

\begin{prob}
Construct and describe families of polytopes with arbitrarily given complexity.
\end{prob}

We mentioned in the introduction the new family of flag nestohedra $\mathcal P_{Mas}$, which will be constructed in Definition~\ref{familyP} below. One of our main results is the construction of this family and the proof of the fact that the family $\mathcal P_{Mas}$ has complexity 4, see Theorem~\ref{boundary}. 

In our subsequent publications we shall continue the study of polytope families with finite complexity, having the DFNM property. In particular, we shall apply Corollary~\ref{ClosurePolytopes} to solve the above problem for the nestohedra family. We shall also construct new families of polytopes with nontrivial Massey products, having finite complexity.


\section{Generating series of nestohedra families and differential equations for them}

We start this section with the following construction due to Erokhovets. It already found various applications, among which is~\cite[Proposition 1.5.23]{TT}. We apply it in the next section in order to find nestohedra families with nontrivial Massey products.

\begin{constr}[{cf.~\cite[Construction 1.5.19]{TT}}]\label{buildsetsubst}
Let $B$ be a connected building set on $[n+1]$ and $B_{i}$ be connected building set on the vertex set $[k_i]$ for $1\leq i\leq n+1$. Then a {\emph{substitution of building sets}} is a connected building set $B(B_{1},\ldots,B_{n+1})$ on the set $[k_{1}]\sqcup\ldots\sqcup [k_{n+1}]=[k_{1}+\ldots+k_{n+1}]$, consisting of the elements $S_{i}\in B_{i}$, as well as $\sqcup_{i\in S}\,[k_i]$, where $S\in B$.
\end{constr}

This construction allows one for any building set $B$ to construct a connected building set $B'$ such that $P_{B}=P_{B'}$, see~\cite[\S1.5]{TT}.
Next, we introduce a new operation on connected building sets.

\begin{constr}
Suppose $B_{1}$ and $B_{2}$ are connected building sets on $[n+1]$ and $B_{1}\cap B_{2}=B_{\Delta}$. If a subset $B_{1}\cup B_{2}$ of the set $2^{[n+1]}$ is a building set, then we say that the {\emph{sum of building sets}} $B_{1}+B_{2}$ is defined and, by definition, is equal to $B_{1}\cup B_{2}$ as subsets of $2^{[n+1]}$.
\end{constr}

\begin{propos}\label{sumbuildset}
Let $B_{1}, B_{2}$ be connected building sets on $[n+1]$ and the following condition holds: for any $S_{i}\in B_{i},i=1,2$ with $S_{1}\cap S_{2}\neq\varnothing$ their union $S_{1}\cup S_{2}\in B_{1}\cup B_{2}$. Then the sum of building sets $B_{1}+B_{2}$ is defined.
\end{propos}

Now we turn to a definition of the flag nestohedra family $\mathcal P_{Mas}$. Consider the following sets of subsets in $[n+1]$ for $n\geq 2$:
$$
B_{1}(P,n)=\{\{i\}|\,1\leq i\leq n+1\},
$$
$$
B_{2}(P,n)=\{\{1,2,i_{1},\ldots,i_{k}\}|\,3\leq i_{1}<\ldots<i_{k}\leq n+1,0\leq k\leq n-1\},
$$
$$
B_{3}(P,n)=\{\{1,j_{1},\ldots,j_{p}\}|\,3\leq j_{1}<\ldots<j_{p}\leq n,1\leq p\leq n-2\}.
$$

\begin{propos}\label{Pbuildset}
The sets $B_{1}(P,n)\cup B_{2}(P,n),B_{1}(P,n)\cup B_{3}(P,n)\cup \{[n+1]\}$ are connected building sets on $[n+1]$ and their sum
$$
B(P,n)=(B_{1}(P,n)\cup B_{2}(P,n))+(B_{1}(P,n)\cup B_{3}(P,n)\cup\{[n+1]\})
$$
is defined. Furthermore, a nestohedron $P_{B(P,n)}$ is flag.
\end{propos}
\begin{proof}
The first statement follows immediately from the definition of a building set. To prove the second one, it suffices to observe that
$$
(B_{1}(P,n)\cup B_{2}(P,n))\cap(B_{1}(P,n)\cup B_{3}(P,n)\cup\{[n+1]\})=B_{\Delta}
$$ 
and if $S_{1}\in B_{1}(P,n)\cup B_{2}(P,n)$, $S_{2}\in B_{1}(P,n)\cup B_{3}(P,n)\cup \{[n+1]\}$ for $S_{1}\cap S_{2}\neq\varnothing$, then $S_{1},S_{2}\in B_{\Delta}$. Therefore, due to Proposition~\ref{sumbuildset}, the second statement holds. It is easy to see that each element from $B_{2}(P,n)\cup B_{3}(P,n)$ can be represented in the form $S_{1}\sqcup S_{2}$ with $S_{1},S_{2}\in B(P,n)$, and so the polytope $P_{B(P,n)}$ is flag, which finishes the proof. 
\end{proof}

\begin{definition}\label{familyP}
Denote by $\mathcal P_{Mas}=\{P^{n}_{Mas},n\geq 0\}$ a family of flag nestohedra, where $P^0_{Mas}$ is a point, $P^1_{Mas}$ is a segment, and $P^{n}_{Mas}=P_{B(P,n)}$ for $n\geq 2$. 
\end{definition}

\begin{figure}[h]
\includegraphics[scale=0.75]{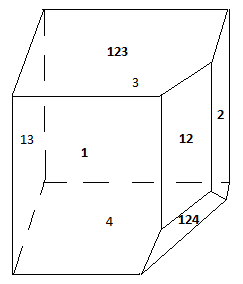}
\caption{3-dimensional polytope from the family $\mathcal P_{Mas}$}
\label{3dim2cube}
\end{figure}

A family f 2-truncated cubes $\mathcal Q$, introduced in~\cite{L1}, (see also Definition~\ref{2truncMassey} below) motivated us to construct the family of flag nestohedra $\mathcal P_{Mas}$. Note that $P^2_{Mas}=Q^2$ is a square. The three-dimensional flag nestohedron $P^3_{Mas}$ is shown in Figure~\ref{3dim2cube} (the visible facets are marked in bold).

Next, we shall apply the theory of differential ring of polytopes $(\mathcal P,d)$ to the family $\mathcal P_{Mas}$. Using Lemma~\ref{NestBound} for nestohedra we can compute the value of the boundary operator $d$ on polytopes from the family $\mathcal P_{Mas}$, as well as on the families of flag nestohedra, which elements arise in $dP$, where $P\in\mathcal P_{Mas}$. 

Consider the following four families of flag nestohedra: $\mathcal P_{Mas}$, permutohedra $Pe$, stellahedra $St$, and also a certain new family of graph-associahedra $\mathcal{P}_{\Gamma}$. We are going to define the latter family now.

Take the following sets of subsets in $[n+1]$:
$$
B_{1}(\Gamma,n)=\{\{i\}|\,1\leq i\leq n+1\},
$$
$$
B_{2}(\Gamma,n)=\{\{i_{1},\ldots,i_{k}\}|\,1\leq i_{1}<\ldots<i_{k}\leq n,2\leq k\leq n\},
$$
$$
B_{3}(\Gamma,n)=\{\{1,j_{1},\ldots,j_{p},n+1\}|\,2\leq j_{1}<\ldots<j_{p}\leq n,0\leq p\leq n-1\}.
$$

\begin{propos}
The sets $B_{1}(\Gamma,n)\cup B_{2}(\Gamma,n)\cup \{[n+1]\}$ and $B_{1}(\Gamma,n)\cup B_{3}(\Gamma,n)$ are connected building sets on the vertex set $[n+1]$, and their sum 
$$
B(\Gamma,n)=(B_{1}(\Gamma,n)\cup B_{2}(\Gamma,n)\cup \{[n+1]\})+(B_{1}(\Gamma,n)\cup B_{3}(\Gamma,n))
$$
is defined. Furthermore, a nestohedron $P^{n}=P_{B(\Gamma,n)}$ is flag.
\end{propos}
\begin{proof}
The proof goes analogously to that of Proposition~\ref{Pbuildset}.
\end{proof}

\begin{definition}\label{polytopesPG}
Denote by $\mathcal P_{\Gamma}=\{P_{\Gamma}^{n},n\geq 0\}$ the set of flag nestohedra, where $P_{\Gamma}^0$ is a point, $P_{\Gamma}^1$ is a segment, and $P_{\Gamma}^{n}=P_{B(\Gamma,n)}$ for $n\geq 2$.
\end{definition}

\begin{remark}
It is not hard to see that $B(\Gamma,n)$ is isomorphic to a graphical building set $B(\Gamma_n)$, where $\Gamma_n$ is a simple graph on $n+1$ vertices, consisting of a complete graph $K_n$ on the set $[n]$ and an edge joining the vertices $\{n\}$ and $\{n+1\}$. 
\end{remark}

To compute the value of the boundary operator on our families of flag nestohedra we need the following formulae for the action of the operator $d$ on the families $Pe$ and $St$, cf.~\cite{B}:
$$
dPe^n=\sum\limits_{s=0}^{n-1}\binom{n+1}{s+1}Pe^{s}\times Pe^{n-s-1},
$$
and
$$
dSt^n=nSt^{n-1}+\sum\limits_{s=0}^{n-1}\binom{n}{s}St^{s}\times Pe^{n-s-1}.
$$

\begin{lemma}\label{lemmboundary}
The operator $d$ acts on the family of graph-associahedra $\mathcal P_{\Gamma}$ by formula:
$$
dP^{n}_{\Gamma}=Pe^{n-1}+\sum\limits_{s=0}^{n-2}\binom{n-1}{s+1}Pe^{s}\times P^{n-s-1}_{\Gamma}+\sum\limits_{s=0}^{n-1}\binom{n-1}{s}Pe^{s}\times Pe^{n-s-1}+
$$
$$
+\sum\limits_{s=0}^{n-2}\binom{n-1}{s}P^{s+1}_{\Gamma}\times Pe^{n-s-2},
$$
In particular, $\mathcal P_{\Gamma}$ has complexity 2 (cf. Defintion~\ref{complex}).
\end{lemma}
\begin{proof}
It follows from the formula for the value of the boundary operator on a nestohedron, given in Lemma~\ref{NestBound}, that for the family $\mathcal{P}_{\Gamma}$ the next identities take place:
$$
B(\Gamma,n)/\{1,i_{1},\ldots,i_k,n+1\}=B(Pe^{n-k-2}),\quad 2\leq i_{1}<\ldots<i_{k}\leq n,
$$
$$
B(\Gamma,n)/\{1,i_{1},\ldots,i_{k}\}=B(Pe^{n-k-1}),\quad 2\leq i_{1}<\ldots<i_{k}\leq n,
$$
$$
B(\Gamma,n)/\{j_{1},\ldots,j_{p}\}=B(\Gamma,n-p),\quad 2\leq j_{1}<\ldots<j_{p}\leq n.
$$
For the corresponding induced building sets we have:
$$
B(\Gamma,n)|_{\{1,i_{1},\ldots,i_{k},n+1\}}=B(\Gamma,k+1),
$$
$$
B(\Gamma,n)|_{\{1,i_{1},\ldots,i_{k}\}}=B(Pe^k),
$$
$$
B(\Gamma,n)|_{\{j_{1},\ldots,j_{p}\}}=B(Pe^{p-1}).
$$
Finally, applying the definition of contraction operation, we get: $B(\Gamma,n)/\{n+1\}=B(Pe^{n-1})$.

Using the formula for the value of the boundary operator on a nestohedron and the above formulae, we conclude that:
$$
dP^{n}_{\Gamma}=Pe^{n-1}+\sum\limits_{s=0}^{n-2}\binom{n-1}{s+1}Pe^{s}\times P^{n-s-1}_{\Gamma}+\sum\limits_{s=0}^{n-1}\binom{n-1}{s}Pe^{s}\times Pe^{n-s-1}+
$$
$$
+\sum\limits_{s=0}^{n-2}\binom{n-1}{s}P^{s+1}_{\Gamma}\times Pe^{n-s-2}.
$$
\end{proof}

\begin{theorem}\label{boundary}
The operator $d$ acts on the family $\mathcal P_{Mas}$ by formula:
$$
dP^n_{Mas}=2St^{n-1}+(n-2)P^{n-1}_{Mas}+\sum\limits_{s=0}^{n-2}\binom{n-2}{s}St^{s}\times P^{n-s-1}_{\Gamma}+
$$
$$
+\sum\limits_{s=0}^{n-2}\binom{n-2}{s}St^{s+1}\times Pe^{n-s-2}+\sum\limits_{s=0}^{n-3}\binom{n-2}{s}P^{s+2}_{Mas}\times Pe^{n-s-3}.
$$
In particular, $\mathcal P_{Mas}$ has complexity 4 (cf. Definition~\ref{complex}).
\end{theorem}
\begin{proof}
It is easy to check that the following identities take place:
$$
B(P,n)/\{1\}=B(\Gamma,n-1), B(P,n)/\{2\}=B(P,n)/\{n+1\}=B(St^{n-1}), 
$$
$$
B(P,n)/\{3\}=\ldots=B(P,n)/\{n\}=B(P,n-1),
$$
$$
B(P,n)/\{1,2,S\}=B(Pe^{n-2-|S|}), B(P,n)|_{\{1,2,S\}}=B(St^{|S|+1}),
$$
$$
B(P,n)/\{1,2,n+1,S\}=B(Pe^{n-3-|S|}), B(P,n)|_{\{1,2,n+1,S\}}=B(P,|S|+2),
$$
$$
B(P,n)/\{1,S\}=B(\Gamma,n-1-|S|), B(P,n)|_{\{1,S\}}=B(St^{|S|}),
$$
for $S\subset\{3,4,\ldots,n\}$.

The proof is finished by applying Lemma~\ref{NestBound}.
\end{proof}

Due to Lemma~\ref{FlagCommuteEmbed} and the previous theorem, for each $n\geq 1$ there exists a retraction $r_{n}\colon\mathcal Z_{P_{Mas}^n}\to\mathcal Z_{P^{n-1}_{Mas}}$. 

\begin{propos}\label{nofib}
For $n\geq 3$ there exists no smooth retraction $\mathcal Z_{P_{Mas}^n}\to\mathcal Z_{P^{n-1}_{Mas}}$.
\end{propos}
\begin{proof}
Denote the moment-angle manifold $\mathcal Z_{P^{n}_{Mas}}$ simply by $M_n$. Suppose that there exists a smooth retraction $r_{n}$. Since $r_{n}i_{n}=id$, where $i_{n}\colon M_{n-1}\to M_n$, the differential $dr_{n}\colon\tau M_{n}\to\tau M_{n-1}$ must be an epimorphism. Thus a retraction $\pi_{n}\colon M_{n}\to M_{n-1}$ is a smooth fibration, for which the map $i_{n}$ provides a section. Let $M$ be the fiber of $\pi_n$. From existence of the section it follows that Serre spectral sequence for the fibration $\pi_{n}$ degenerates, that is, $E_{\infty}=E_{2}=H^{*}(M_{n-1};\,H^{*}(M))$. In the latter case we get an additive isomorphism $H^*(M_{n};\Q)\cong H^{*}(M_{n-1};\Q)\otimes H^{*}(M;\Q)$. Thus, for their Hilbert polynomials ($H(M;t):=\sum\,b_{i}t^{i}$) one must have: 
$$
H(M_{n};t)=H(M_{n-1};t)\cdot H(M;t).
$$
Therefore, the polynomial $H(M_{n-1};t)$ must be a factor of the polynomial $H(M_{n};t)$. From calculations with bigraded Betti numbers of flag nestohedra of $\mathcal P_{Mas}$ it follows immediately that it is a contradiction for all $n\geq 3$.
\end{proof}

Let us turn now to the problem of generating series and differential equations for them, in the case of flag nestohedra. Set
$$
Pe(x)=\sum\limits_{s=0}^{\infty}Pe^{s}\frac{x^{s+1}}{(s+1)!},\quad St(x)=\sum\limits_{s=0}^{\infty}St^{s}\frac{x^s}{s!}
$$
and
$$
P_{\Gamma}(x)=\sum\limits_{s=0}^{\infty}P^{s+1}_{\Gamma}\frac{x^s}{s!},\quad
P_{Mas}(x)=\sum\limits_{s=0}^{\infty}P^{s+2}_{Mas}\frac{x^{s+2}}{s!}.
$$

Applying Theorem~\ref{boundary} we are going to get differential equations, which are satisfied by above defined generating series of our families.
In the work by Buchstaber~\cite{B} the next formulae were obtained:
$$
dPe(x)=Pe^2(x),\,
dSt(x)=(x+Pe(x))St(x).
$$

Studying generating series for classical polynomials in their work~\cite{B-Kh}, Buchstaber and Kholodov suggested a fruitful idea: to consider generating series with certain structure sequences of coefficients $\{a_n\}$. This allowed them, varying the coefficients, to obtain generating series for different classes of classical polynomials, and therefore, to study them simultaneously, by changing the structure constants $a_n$.

We shall use the same idea to study the generating series of our families of polytopes. For a generating series of the type $Q(x)=\sum\limits_{n=0}^{\infty}a_{n}Q^{n}x^{n+n_0}$ set, by definition, $dQ(x)=\sum\limits_{n=0}^{\infty}a_{n}d(Q^{n})x^{n+n_0}$. 

\begin{theorem}\label{1paramGenSer}
The following identities hold:
$$
dP_{\Gamma}(x)=2Pe(x)P_{\Gamma}(x)+\Bigl(1+\frac{d}{dx}Pe(x)\Bigr)\frac{d}{dx}Pe(x),
$$
$$
dP_{Mas}(x)=(x+Pe(x))P_{Mas}(x)+x^{2}(2\frac{d}{dx}St(x)+St(x)P_{\Gamma}(x)+
$$
$$
+\frac{d}{dx}St(x)\frac{d}{dx}Pe(x)).
$$
\end{theorem}
\begin{proof}
The proof goes by direct computation, using Lemma~\ref{lemmboundary} and Theorem~\ref{boundary}.
\end{proof}

Now, following~\cite{BE}, we introduce a new paprameter and study the 2-paprameter generating series for our families of flag nestohedra.

\begin{definition}\label{2paramSer}
For the polytope family $\mathcal P=\{P^{n}|\,n\geq 0\}$ consider its generating series $P(x)=\sum\limits_{n=0}^{\infty}a_{n}P^{n}x^{n+n_0}$ and the following formal series ($P^{n}\in\mathcal P$):
$$
Q(P^{n};q)=\sum\limits_{k=0}^{n}(d^{k}P)\frac{q^k}{k!}.
$$
Then we can define a {\emph{2-parameter extension}} of a generaing series for the family $\mathcal P$ by formula:
$$
P(q,x)=\sum\limits_{n=0}^{\infty}a_{n}Q(P^{n};q)x^{n+n_0}.
$$
Note that for any family $\mathcal P$: $Q(P^{n};0)=P^{n}$ and $P(0,x)=P(x)$.
\end{definition}

In~\cite{B}, as well as in the works by Buchstaber and Volodin~\cite{BV, BV2} partial differential equations were obtained, the solutions of which were the 2-parameter generating series of well known linear families of flag nestohedra $\{P^n|\,n\geq 0\}$, where $P^n$ denotes either an $n$-dimensional simplex $\Delta^n$, or a cube $I^n$, or an associahedron $As^n$, or a cyclohedron $Cy^n$, or a permutohedron $Pe^n$, or a stellahedron $St^n$. These equations carry important information about complexity of the combinatorial structure of the family of polytopes. 

In particular, let us mention that the 2-parameter generating series $Pe(q,x)$ and $St(q,x)$ give solutions to the following Cauchy problems:
\begin{itemize}
\item[(1)] For the 2-parameter generating series of the permutohedra family one has:
$$
\frac{\partial}{\partial q}Pe(q,x)=Pe^2(q,x), Pe(0,x)=Pe(x).
$$
The solution takes the form:
$$
Pe(q,x)=\frac{Pe(x)}{1-qPe(x)};
$$
\item[(2)] For the 2-parameter generating series of the stellahedra family one has:
$$
\frac{\partial}{\partial q}St(q,x)=(x+Pe(q,x))St(q,x), St(0,x)=St(x).
$$
The solution takes the form:
$$
St(q,x)=St(x)\frac{e^{qx}}{1-qPe(x)}.
$$
\end{itemize}
 
\begin{theorem}\label{2paramGenSer}
The next statements hold.
\begin{itemize}
\item[(1)] The 2-parameter generating series $P_{\Gamma}(q,x)$ is a solution of the following Cauchy problem:
$$
\frac{\partial}{\partial q}P_{\Gamma}(q,x)=2Pe(q,x)P_{\Gamma}(q,x)+\frac{\partial}{\partial x}Pe(q,x)(1+\frac{\partial}{\partial x}Pe(q,x)), P_{\Gamma}(0,x)=P_{\Gamma}(x);
$$
\item[(2)] The 2-parameter generating series $P_{Mas}(q,x)$ is a solution of the following Cauchy problem:
$$
\frac{\partial}{\partial q}P_{Mas}(q,x)=(x+Pe(q,x))P_{Mas}(q,x)+x^{2}[2\frac{\partial}{\partial x}St(q,x)+St(q,x)P_{\Gamma}(q,x)+
$$
$$
+\frac{\partial}{\partial x}St(q,x)\frac{\partial}{\partial x}Pe(q,x)], P_{Mas}(0,x)=P_{Mas}(x).
$$
\end{itemize}
\end{theorem}
\begin{proof}
Direct computation with the use of Theorem~\ref{1paramGenSer}. 
\end{proof}

Denote the right hand sides of the equations from the previous theorem by $AF+B$, where $F=F(q,x)$ is the desired solution. Let us find solutions of the above mentioned Cauchy problems in the form
$F(q,x)=F_{1}(q,x)\cdot F_{2}(q,x)$ under conditions $\frac{\partial}{\partial q}F_{1}=A\cdot F_{1}$ and $F_{1}\cdot\frac{\partial}{\partial q}F_{2}=B$. 

Then it is easy to get the next formulae for the functions $F_{1}(q,x)$ in each of the above two cases. Indeed, for $F(q,x)=P_{\Gamma}(q,x)$ we get:
$F_{1}(q,x)=F_{1}(0,x)\frac{1}{(1-qPe(x))^{2}}$. For $F(q,x)=P_{Mas}(q,x)$ we get: $F_{1}(q,x)=F_{1}(0,x)\frac{e^{qx}}{1-qPe(x)}$.

In our subsequent publications we plan to give complete solutions to the above Cauchy problems, to analyze the behaviour of the solutions, and interpret the corollaries from those solutions in terms of combinatorics of the corresponding families.


\section{Flag nestohedra and Massey products}

We recall the results on Massey products in cohomology of moment-angle manifolds over 2-truncated cubes, which will play an important role when moving to the case of flag nestohedra.

The first examples of moment-angle manifolds $\zp$ with nontrivial higher Massey roducts  $\langle\alpha_{1},\ldots,\alpha_{k}\rangle$ for any  $k\geq 4$ in $H^*(\zp)$ were constructed by Limonchenko~\cite{L1,L2}. Namely, he introduced a family $\mathcal Q=\{Q^n|\,n\geq 0\}$ such that there exists a nontrivial Massey product $\langle\alpha_{1},\ldots,\alpha_{n}\rangle$ of order $n$ with $\dim\alpha_{i}=3,1\leq i\leq n$ in $H^*(\mathcal Z_{Q^n})$ for any $n\geq 2$.

We start with the construction of the linear family of 2-truncated cubes $\mathcal Q=\{Q^n|\,n\geq 0\}$, for which $\mathcal Z_{Q^{n}}$ has, by~\cite[Theorem 3.6]{L3}, a strictly defined nontrivial $n$-fold Massey product in cohomology for each $n\geq 2$.

\begin{definition}[\cite{L1,L2}]\label{2truncMassey}
Let $Q^0$ be a point and $Q^1\subset\R^1$ be a segment $[0,1]$. Denote by $I^{n}=[0,1]^n, n\geq 2$ an $n$-dimensional cube with facets $F_{1},\ldots,F_{2n}$ in such a way that $F_{i},1\leq i\leq n$ contains the origin 0, $F_{i}$ and $F_{n+i}$ are parallel for all $1\leq i\leq n$. 
Then the face ring of the cube has the form:
$$
\ko[I^n]=\ko[v_{1},\ldots,v_{n},v_{n+1},\ldots,v_{2n}]/I_{I^n},
$$
where $I_{I^n}=(v_{1}v_{n+1},\ldots,v_{n}v_{2n})$.

Consider a polynomial ring
$$
\ko[v_{1},\ldots,v_{2n},v_{k',n+k'+i'}|\,1\leq i'\leq n-2, 1\leq k'\leq n-i']
$$
and its monomial ideal, generated by square free monomials:
$$
I=(v_{k}v_{n+k+i},v_{k',n+k'+i'}v_{n+k'+l},v_{k',n+k'+i'}v_{p},v_{k',n+k'+i'}v_{k'',n+k''+i''}),
$$
where $v_{j}$ corresponds to $F_{j}$ for all $1\leq j\leq 2n$, and 
$$
0\leq i\leq n-2, 1\leq k\leq n-i, 1\leq i',i''\leq n-2, 1\leq k'\leq n-i', 
$$
$$
1\leq k''\leq n-i'', 1\leq p\neq k'\leq k'+i', 0\leq l\neq i'\leq n-2, 
$$
$$
k'+i'=k''\,\text{or }k''+i''=k'.
$$

Let us define $Q^n\subset\R^n$ to be a simple polytope such that $I_{Q^n}=I$. Note that $Q^n$ has a natural realization as a 2-truncated cube, and moreover, its combinatorial type does not depend on the order in which the faces of the cube $I^n$ are truncated (the generators $v_{i,j}$ correspond to the truncated faces $F_{i}\cap F_{j}$ of $I^n$).
\end{definition}

The next result on (higher) Massey products in cohomology of moment-angle manifolds holds.

\begin{theorem}[{\cite{L2,L3}}]\label{mainMassey}
Let $\alpha_i\in H^{3}(\mathcal Z_{Q^n})$ be represented by 3-dimensional cocycles of the type $v_{i}u_{n+i}\in R^{-1,4}(Q^n)$ for all $1\leq i\leq n$ and $n\geq 2$. Then all the Massey products of the consecutive elements from the ordered set $(\alpha_{1},\ldots,\alpha_{n})$ are defined and the $n$-fold Massey product $\langle\alpha_{1},\ldots,\alpha_{n}\rangle$ is strictly defined and nontrivial.
\end{theorem}

\begin{theorem}\label{Qdfpnm}
For each $Q^n\in\mathcal Q$ and any $0\leq r<n$ there exists a face $F^{r}\subset Q^n$ such that $F^r$ is combinatorially equivalent to $Q^r$.
\end{theorem}
\begin{proof}
To prove the statement it suffices to show that the polytope $Q^n$ has a facet equal to $Q^{n-1}$ for all $n\geq 1$.

Indeed, consider the facet $F_{n-1}$ of the polytope $Q^n$. Let us show that $F_{n-1}=Q^{n-1}$. By definition of $Q^n$, $F_{n-1}$ is obtained from the facet $G_{n-1}=I^{n-1}$ of the cube $I^{n}$ with facets $G_{i},1\leq i\leq 2n$ (recall that we set $G_{i}\cap G_{n+i}=\varnothing$ above, for all $1\leq i\leq n$, cf. Definition~\ref{2truncMassey}) by a set of consecutive truncations of faces of codimensions two, and moreover, is also a 2-truncated cube itself (more generally, each face of a 2-truncated cube is a 2-truncated cube), see~\cite[\S1.6]{TT}. The truncated faces of codimension two, in the case of $F_{n-1}$, have the form (here we use the fact that truncation of a facet does not effect the combinatorial type of a polytope and that $G_{n-1}\cap G_{2n-1}=\varnothing$):
$$
G_{n-1}\cap G_{i}\cap G_{j},
$$
where $1\leq i\leq n-2$, $n+2\leq j\leq 2n$, $j\neq 2n-1$. It suffices to observe that the sets of minimal nonfaces $\MF(F_{n-1})$ and $\MF(Q^{n-1})$ coincide (up to changing the notation), thus $F_{n-1}$ is combinatorially equivalent to $Q^{n-1}$. It can be checked similarly that $F_{2n-1}=Q^{n-1}$.

Furthermore, the proof shows that $F_{r}\cap F_{r+1}\cap\ldots\cap F_{n-1}=Q^{r}$.
\end{proof}

Now, Theorem~\ref{FlagCriterion} and the theorem just proved, taken together, immediately imply the next result.

\begin{coro}\label{AllProducts}
There exists a nontrivial strictly defined $k$-fold Massey product in $H^*(\mathcal Z_{Q^n})$ for any $k, 2\leq k\leq n$.
\end{coro}

\begin{remark}
Note that it was proved in~\cite[Theorem 4.1]{L3} that a full subcomplex in $K=K_{Q^n}$ on the vertex set $\{1,\ldots,r-1,n,n+1,\ldots,n+r-1,2n\}$ is combinatorially equivalent to the full subcomplex in $K_{Q^r}$ on the vertex set $[2r]$ for each $r, 2\leq r\leq n$. This gives another proof of the above Corollary~\ref{AllProducts}.
\end{remark}

\begin{propos}\label{Qdfpnm}
The family $\mathcal Q$ is a special geometric direct family with nontrivial Massey products.
\end{propos}
\begin{proof}
Due to Corollary~\ref{AllProducts}, it suffices to prove that $\mathcal Q$ is an algebraic direct family, cf. Definition~\ref{ADFP}. The latter follows from Theorem~\ref{Qdfpnm}.
\end{proof}

Starting with an arbitrary irreducible flag polytope $F^r$, we can construct a family of irreducible flag polytopes $\{P^n|\,n\geq r\}$ having $P^{r}=F^{r}$ and such that $P^{k}$ is a face of $P^{l}$ for any $l>k\geq r$.

\begin{constr}\label{familyFlag}
If given a flag polytope $F^r$, one can define a sequence of flag polytopes $\mathcal P(F)=\{P^n|\,n\geq r\}$ in the following way. Set $P^{r}=F^{r}$ and if $P^n$ is already defined we can construct $P^{n+1}$ as a result of a truncation from the polytope $P^{n}\times I$ of a certain face of codimension two of the type $F_{i(n)}\times\{1\}\subset P^{n}\times\{1\}\subset P^{n}\times I$ (here,  $F_{i(n)}$ denotes a facet of the polytope $P^n$). Observe that the resulting polytope is flag again and has $P^n$ as its facet $P^{n}\times\{0\}\subset P^{n+1}$ for each $n\geq r$. Obviously, the combinatorial type of $P^n$ for $n>r$, in general, depends not only on the original polytope $F^r$, but also from the choice of the truncated facet $F_{i(n)}$ in $P^{n}$.
\end{constr}

In the above construction we determined a new operation $\fc$ (``face cut'') on the set of flag polytopes; one has: $P^{n+1}=\fc(P^n)$ for all $n\geq r$. Let us denote by $Q=\fc^{k}(P)$ the polytope, obtained from $P$ by consecutive application of $k$ operations, described in Construction~\ref{familyFlag}; note that $P=\fc^{0}(P)$ and $\dim Q=\dim P+k$. 

As an application of the construction described above, we immediately get the following result.

\begin{coro}\label{MasseySeq}
Let $F^r$ be a flag polytope and suppose that there exists a nontrivial $k$-fold Massey product in $H^*(\mathcal Z_{F^r})$. Then there exists a sequence of polytopes $\mathcal P=\{P^n|\,n\geq r\}$ such that there exists a nontrivial $k$-fold Massey product in $H^*(\mathcal Z_{P^n})$ for all $n\geq r$. 
\end{coro}
\begin{proof}
The proof goes by applying Theorem~\ref{FlagCriterion} to a sequence of flag polytopes $\mathcal P(F^r)=\{P^n|\,n\geq r\}$, defined in the Construction~\ref{familyFlag}.
\end{proof}

\begin{definition}
A sequence of irreducible flag polytopes $\mathcal P=\{P^n\}$ such that there exists a nontrivial $k$-fold Massey product in $H^*(\mathcal Z_{P^n})$ with $k\to\infty$ as $n\to\infty$ and, moreover, existence of a nontrivial $k$-fold Massey product in $H^*(\mathcal Z_{P^n})$ implies existence of a nontrivial $k$-fold Massey product in $H^*(\mathcal Z_{P^l})$ for all $l>n$, we will call a \emph{sequence of polytopes with strongly related Massey products}.
\end{definition}

\begin{remark}
Each geometric direct family of polytopes with nontrivial Massey products is an example of a sequence of polytopes with strongly related Massey products, cf. Definition~\ref{GDFP} and Definition~\ref{DFPM}.
\end{remark}

\begin{definition}
We say that two sequences of polytopes $\mathcal P_{1}=\{P_{1}^{n}\}$ and $\mathcal P_{2}=\{P_{2}^{n}\}$ are \emph{combinatorially different}, if for any $N\geq 0$ there exists $n>N$ such that $P_{1}^{n}$ and $P_{2}^{n}$ are not combinatorially equivalent.
\end{definition}

\begin{propos}\label{MasseySeqInfinity}
There exists an infinite set $\mathcal S=\{\mathcal P_{\alpha}|\,\alpha\in I\}$ of pairwisely combinatorially different sequences of polytopes with strongly related Massey products.
\end{propos}
\begin{proof}
Observe that if $P^{k}, P^{l}\in\mathcal P(F)$ ($k<l$) for a certain flag polytope $F$, then $P^{k}$ is a proper face of $P^l$ and, moreover, $\fc(P^k)=P^{k+1}$ is a face of $P^l\subset\fc(P^l)$. For any sequence $s\in\{0,1\}^{\infty}$, not stabilizing at zero, we can define a sequence of flag polytopes $\mathcal P_{s}=\{P_{s}^{n}|\,n\geq 4\}$ such that $P_{s}^{n}=Q^{n}$, if $s_{n}=1$, and $P_{s}^{n}=\fc^{n-k(n)}(Q^{k(n)})$, if $s_{n}=0$. Here, $k(n)=1$, if $s_{1}=\ldots=s_{n}=0$, and $k(n)=\max\,\{m|\,m<n, s_{m}=1\}$, otherwise. Then from Theorem~\ref{FlagCriterion} it follows that $\mathcal P_{s}$ is a sequence of polytopes with strongly related Massey products. Furthermore, the sequence $\mathcal P_{s}$ is combinatorially different with $\mathcal P_{s'}$, 
if $|s-s'|$ is not stabilizing at zero, since $m(\fc(P))=m(P)+3$ for any polytope $P$ and $m(Q^n)=\frac{n(n+3)}{2}-1$. 
\end{proof}

Another way to get different sequences of flag polytopes with nontrivial (higher) Massey products in cohomology of their moment-angle manifolds is given in the next statement.

\begin{theorem}\label{mainMasseysequence}
There exist infinitely many sequences $\mathcal P_{k}=\{P_{k}^n\}, k\geq 1$ of simple flag irreducible polytopes such that: 
\begin{itemize}
\item[(a)] If $P\in\mathcal P_{i}$ and $Q\in\mathcal P_{j}$ for $i\neq j$, then $P$ and $Q$ are not combinatorially equivalent;
\item[(b)] For any $k\geq 1$ and $r\geq 2$ there exists $N=N(k,r)$ such that there exists a strictly defined nontrivial $l$-fold Massey product in $H^*(\mathcal Z_{P_{k}^n})$, for all $2\leq l\leq r, n\geq N$.
\end{itemize}
\end{theorem}
\begin{proof}
Consider the following sequences of flag polytopes: $\mathcal P_{k}=\{\fc^{k-1}(Q^{n})|\,n\geq 3\}$, $k\geq 1$. 

To prove (a) suppose that the converse is true; then the following identity for the numbers of facets of combinatorially equivalent polytopes of the same dimension holds:
$$
m(\fc^{k}(Q^l))=m(\fc^{l}(Q^k)).
$$
Since $m(Q^n)=\frac{n(n+3)}{2}-1$ for all $n\geq 2$, and $m(\fc(P))=m(P)+3$, the above formula implies $(l-k)(l+k-3)=0$. For $l=k$ both polytopes belong to the same sequence $\mathcal P_{l+1}$, and for $l+k=3$ one of the dimensions, $l$ or $k$, will be greater than 3. In both cases we get a contradiction with the definition of sequences $\mathcal P_k$, which finishes the proof of statement (a).

Applying Construction~\ref{familyFlag}, we obtain: $Q^n$ is a face of each polytope of dimension greater or equal to $n+k-1$ in the sequence $\mathcal P_k$. Then statement (b) (with $N(k,r)=r+k-1$) follow from Theorem~\ref{FlagCriterion} and Corollary~\ref{AllProducts}. This finishes the proof.
\end{proof}


Now, we turn to studying the family of flag nestohedra $\mathcal P_{Mas}$ constructed above. 
We shall prove that $\mathcal P_{Mas}$ is a linear geometric direct family of flag nestohedra with nontrivial Massey products, which motivates its notation.

\begin{lemma}\label{FlagNestMassey}
For each $P^{n}_{Mas}\in\mathcal P_{Mas}$ with $n\geq 2$ there exists a strictly defined nontrivial $n$-fold Massey product in $H^*(\mathcal Z_{P^n_{Mas}})$.
\end{lemma}
\begin{proof}
Let us find a realization of the polytope $P^{n}_{Mas}\in\mathcal P_{Mas}$ for $n\geq 2$ as a 2-truncated cube, by means of the iterated procedure of codimension two face truncations from $I^n$, described in~\cite{BV}. In terms of the Definition~\ref{2truncMassey}, we identify $F_{i}$ with the set $\{1,\ldots,i\}$ for $1\leq i\leq n$ and we identify $F_{i}$ with the set $\{i-n+1\}$ for $n+1\leq i\leq 2n$. Then we consecutively truncate the following faces of codimension two from the cube $I^n$:
$$
\{1\}\sqcup\{3\}, \{1,2\}\sqcup\{4\},\ldots,\{1,\ldots,n-1\}\sqcup\{n+1\}\\
$$
$$
\cdots\\
$$
$$
\{1\}\sqcup\{n\}, \{1,2\}\sqcup\{n+1\}.
$$
in the order opposite to inclusion. It is easy to see that the building set $B=B(P,n)$ can be represented as a union of the connected building set $B_0$ of the cube $I^n$, the set $B_{1}$ of the above described subsets of $[n+1]$, and the set of all subsets of $[n+1]$, which are the unions of nontrivially intersecting (that is, intersecting and not being subsets of one another) elements from $B$: if $S_{1},S_{2}\in B$ and $S_{1}\cap S_{2}\neq\varnothing,S_{1},S_{2}$, then we must have: $S_{1}\cup S_{2}\in B$ by definition of a building set. We get that $S=S_{1}\sqcup S_{2}\in B(P,n)$ with $S_{1},S_{2}\in B_0$ if and only if $S\in B_1$, and therefore, full subcomplexes: on the vertex set $[2n]$ in $K_{Q^n}=\partial Q^*$ (here $Q^n$ denotes the polytope from Theorem~\ref{mainMassey}) and on the vertex set $\{\{i\},\{1,2,\ldots,k\}|1\leq i\leq n+1, 2\leq k\leq n\}$ in $K_{P^n_{Mas}}$ are isomorphic. To finish the proof it suffices to apply Theorem~\ref{mainMassey} and Theorem~\ref{zkcoh}.
\end{proof}

\begin{remark}
Observe that the polytopes $Q^n$ and $P^n_{Mas}$ are combinatorially different for all $n\geq 4$, since $f_{0}(P^n_{Mas})=3\times 2^{n-2}+n-1>f_{0}(Q^n)=\frac{n(n+3)}{2}-1$ for $n\geq 4$. Note also that for any $n>2$ the polytope $P^n_{Mas}\in\mathcal P_{Mas}$ is a flag nestohedron, but not a graph-associahedron. The latter follows from Theorem~\ref{boundary}, the fact that $P^{3}_{Mas}=Q^{3}$ is not a graph-associahedron ($m(Q^3)<m(As^3)$), and Lemma~\ref{NestBound}, since each facet of a graph-associahedron is a product of graph-associahedra.
\end{remark}


\begin{definition}(\cite{FS})
Suppose $B$ is a building set. The nerve complex $K_P=\partial P^*$ of the nestohedron $P=P_B$ is called a {\emph{nested set complex}} and is denoted by $N_B$.
\end{definition}

The poset of elements of the nested set complex $N_B$ can be described in terms of the structure of the building set $B$ in the following way.
 
\begin{propos}[{\cite[Theorem 1.5.13]{TT}}]\label{Fposet}
The vertices of $N_{B}$ are in one-to-one correspondence with nonmaximal elements of the building set $B$.\\ 
Furthermore, the set of vertices, corresponding to such elements: $S=\{S_{i_1},\ldots, S_{i_k}\}$, is a simplex in $N_B$ if and only if:
\begin{itemize}
\item[(1)] For any two elements $S_{i_p},S_{i_q}$ with $1\leq p,q\leq k$, one has: $S_{i_p}\cap S_{i_q}=\varnothing$, or $=S_{i_p}$, or $=S_{i_q}$;

\item[(2)] If elements of the subset $\{S_{i_{t_1}},\ldots,S_{i_{t_l}}\}\subseteq S,l\geq 2$ are pairwisely disjoint, then $S_{i_{t_1}}\sqcup\ldots\sqcup S_{i_{t_l}}\notin B$.
\end{itemize}
\end{propos}

\begin{lemma}\label{NestSetRestrict}
Let $P=P_B$ be a nestohedron on a connected building set $B$ on $[n+1]$. Then for any element $S\in B$ a full subcomplex of a nested set complex $N_{B}$ on the vertex set $B|_S$ is combinatorially equivalent to a nested set complex $N_{B|_S}$. 
\end{lemma}
\begin{proof}
Due to Proposition~\ref{Fposet}, it suffices to observe that for pairwisely disjoint subsets $S_{i_{t_1}},\ldots,S_{i_{t_l}}\in B_S$ we have: $S_{i_{t_1}},\ldots,S_{i_{t_l}}\subseteq S$ and $S_{i_{t_1}}\sqcup\ldots\sqcup S_{i_{t_l}}\subseteq S$, and so $S_{i_{t_1}}\sqcup\ldots\sqcup S_{i_{t_l}}\notin B$ holds if and only if $S_{i_{t_1}}\sqcup\ldots\sqcup S_{i_{t_l}}\notin B|_S$.
\end{proof}

\begin{exam}
1. If $P=Pe^n$, then for any $S\in B(P)$ one has: $P_{B|_S}=Pe^{|S|-1}$ and $N_{B|_S}$ is a full subcomplex in $N_{B}$ (both nested set complexes are barycentric subdivisions of boundaries of simplices of corresponding dimensions).
2. If $P=As^n$, then for any chain graph $S$ on the vertex set $[i],i\leq n+1$ one has: $P_{B|_S}=As^{|S|-1}$ and $N_{B|_S}$ is a full subcomplex in $N_B=\partial (As^{n})^*$.
\end{exam}

\begin{remark}
Note that Lemma~\ref{NestSetRestrict} does not hold, when the restriction $B|_S$ is replace by contraction $(B/S)\cap B$. Indeed, consider $B=B(P,4)$, that is, a building set for $P^{4}\in\mathcal P_{Mas}$, $S_{1}=\{2\}$, $S_{2}=\{4\}$ and $S=\{1,3\}$. Then $S\in B$, $S_{1},S_{2}\in B/S$, $(B/S)\cap B=\{\{2\},\{4\}\}$ and $S_{1}\sqcup S_{2}\notin B$. However, $S_{1}\sqcup S_{2}\in B/S$, since $S_{1}\sqcup S_{2}\sqcup S=\{1,2,3,4\}\in B$. We immediately get that the full subcomplex on the vertex set $(B/S)\cap B$ in $N_{B}$ is a one-dimensional simplex (segment), but $N_{B/S}$ consists of two points.
\end{remark}

Consider a connected building set $B$ on $[n+1]$ and an element $S\in B\backslash [n+1]$. By Construction~\ref{mapmfds}, we have an induced embedding of moment-angle manifolds $\hat{\phi}_{S}:\,\mathcal Z_{P_{B|_S}}\to\mathcal Z_{P_B}$, and due to Lemma~\ref{NestSetRestrict} and Construction~\ref{mapmcxs}, we have an induced embedding of moment-angle-complexes $j_{S_{*}}:\,\mathcal Z_{N_{B|_S}}\to\mathcal Z_{N_B}$, having a retraction, where we denoted $N_B=K_{P_B}$ and $N_{B|_S}=K_{(P_{B|_S})}$. Combinatorial properties of nesohedra allow us to obtain the next result.

\begin{propos}\label{phiPsi}
Suppose $P_B$ is a nestohedron. Then the following diagram commutes
$$\begin{CD}
  \mathcal Z_{P_{B|_S}} @>\hat{\phi}_{S}>> \mathcal Z_{P_{B}}\\
  @VVh_{1} V\hspace{-0.2em} @VVh_{2} V @.\\
  \mathcal Z_{N_{B|_S}} @>j_{S_*}>>\mathcal Z_{N_B},
\end{CD}\eqno 
$$
where $h_{1}$ and $h_{2}$ are homeomorphisms. In particular, $\hat{\phi}_{S}$ induces split epimorphism in cohomology for any element $S\in B\backslash [n+1]$.
\end{propos}
\begin{proof}
This follows directly from Lemma~\ref{NestSetRestrict}.
\end{proof}

Note that Proposition~\ref{phiPsi} implies that the induced embedding of moment-angle manifolds $\hat{\phi}_{S}$ has a retraction (and therefore, it induces split epimorphism of cohomology rings) for any nestohedron $P_B$, not necessarily a flag one.

Let $B$ be a building set and $S\in B$. Denote by $P_S$ the nestohedron on a restriction of the building set $P_{B|_S}$.

\begin{propos}\label{CohomEpi}
The following statements hold.
\begin{itemize}
\item[(1)] If $j_{S}:\,N_{B|_S}\hookrightarrow N_B$, then $j^{*}_{S}:\,H^*(\mathcal Z_{P_B})\rightarrow H^*(\mathcal Z_{P_S})$ is a natural split epimorphism of cohomology rings;
\item[(2)] If $\langle\alpha_{1},\ldots,\alpha_{k}\rangle\in H^*(\mathcal Z_{P_S})$ is a defined Massey product, then there exist cohomology classes $\beta_{t}\in H^*(\mathcal Z_{P_B})$ for all $1\leq t\leq k$ such that $j^*(\beta_{t})=\alpha_t$, a Massey product $\langle\beta_{1},\ldots,\beta_{k}\rangle$ is defined, and $j^*\langle\beta_{1},\ldots,\beta_{k}\rangle=\langle\alpha_{1},\ldots,\alpha_{k}\rangle$. Furthermore, if $\langle\alpha_{1},\ldots,\alpha_{k}\rangle$ is nontrivial (resp. strictly defined), then $\langle\beta_{1},\ldots,\beta_{k}\rangle$ is nontrivial (resp. strictly defined).
\end{itemize}
\end{propos}
\begin{proof}
Statement (1) is a direct consequence of Proposition~\ref{phiPsi} and Corollary~\ref{splitepi}.
To prove statement (2) we firstly observe that $N_{B|_S}$ is a full subcomplex in the nested set complex $N_B$ by Lemma~\ref{NestSetRestrict}, and therefore, we can apply statement (1) to the cohomology homomorphism $j^{*}_{S}$, which is induced by embedding of building sets $B|_S\subset B$; the latter map also induces an embedding of simplicial complexes $j_{S}:\,N_{B|_S}\hookrightarrow N_B$. The proof is finished by applying Theorem~\ref{zkcoh}. 
\end{proof}

Using Proposition~\ref{CohomEpi}, we are able to get the following statement on a linear family  $\mathcal F$ of nestohedra, which provides us with a sufficient condition for $\mathcal F$ to be a GDFP. 

\begin{propos}\label{NestDirect}
Let $\mathcal F=\{P_{B(n)}|\,n\geq 0\}$ be a linear family of nestohedra on connected building sets $B(n)$ on $[n+1]$ for $n\geq 0$. If for any $n>2$ there exists an element $S(n)\in B(n)$ with $|S(n)|=n$ such that $P_{B|_{S(n)}}=P_{B(n-1)}$, then the family $\mathcal F$ is a geometric direct family of polytopes.\\
In particular, the following nestohedra families are GDFP: cubes $\mathcal I$, permutohedra $Pe$, stellahdra $St$, cyclohedra $Cy$, and associahedra $As$.
\end{propos}
\begin{proof}
By Lemma~\ref{NestSetRestrict}, $N_{B|_{S(n)}}$ is a full subcomplex in the nested set complex $N_{B(n)}$ for each $n>2$. On the other hand, a formula for the value of the boundary operator $d$ on the family $\mathcal F$ (cf. Lemma~\ref{NestBound}) shows that $P_{B|_{S(n)}}=P_{B(n-1)}$ is combinatorially equivalent to a facet of $P_{B(n)}$. According to Proposition~\ref{CohomEpi} (1), an embedding $j_{S(n)}:\,N_{B|_{S(n)}}\hookrightarrow N_{B(n)}$ of a full subcomplex into a nested set complex induces a split epimorphism in cohomology of their moment-angle-complexes. The proof of the first part of the statement now finished by applying Proposition~\ref{phiPsi} and induction on $n$. The rest follows from explicit formulae for the boundary operator $d$ action on the nestohedra families mentioned above, see~\cite{B}, and the first part of the statement.
\end{proof}

\begin{exam}
Note that for an algebraic direct family of simplices $\Delta$ a condition from Proposition~\ref{NestDirect} does not hold, since for any $n>2$ and each $S\in B(n)\backslash [n+1]$ one has: $P_{B|_S}$ is a point, and an arbitrary full subcomplex $(K_P)_J$ is either $K_P$, or a simplex, for any $P\in\Delta$.
\end{exam}

\begin{theorem}\label{MasAllOrder}
For $P^{n}_{Mas}\in\mathcal P_{Mas}$ there exist a strictly defined nontrivial $k$-fold Massey product $\langle\alpha_{1}^{n},\ldots,\alpha_{k}^{n}\rangle$ in $H^*(\mathcal Z_{P^n_{Mas}})$ for all $2\leq k\leq n$.
Furthermore, for each $2\leq r\leq s$ there exists a natural embedding $j_{r}^{s}:\,N_{B(P,r)}\hookrightarrow N_{B(P,s)}$ such that $(j_{r}^{s})^{*}\langle\alpha_{1}^{s},\ldots,\alpha_{k}^{s}\rangle=\langle\alpha_{1}^{r},\ldots,\alpha_{k}^{r}\rangle$ for all $2\leq k\leq r$.
\end{theorem}
\begin{proof}
Set $B=B(P,n)$. By Theorem~\ref{boundary} for $S_{t}=\{1,2,n+1,I\}$, where we denoted $I=\{3,4,\ldots,t\}$ and $2\leq t\leq n-1$ (for $t=2$ we set $S_{t}=\{1,2,n+1\}$ and $I=\varnothing$), we obtain: $P_{B|_{S_t}}=P^{t}_{Mas}\in\mathcal P_{Mas}$. Consider a natural embedding of building sets $B|_{S_{t}}\subset B$. Due to Lemma~\ref{NestSetRestrict}, the latter map induced an embedding of nested set complexes $j_{t}^{n}:\,N_{B|_{S_t}}\hookrightarrow N_{B}$, that is, there exists a full subcomplex in $N_{B(P,n)}=K_{P^n_{Mas}}$ isomorpic to $N_{B(P,k)}$ for all $2\leq k\leq n-1$. The proof finishes now by applying Lemma~\ref{FlagNestMassey} and Proposition~\ref{CohomEpi}. 
\end{proof}

\begin{remark}
An alternative way to prove the above theorem goes by applying~\cite[Corollary 3.7]{L3} to the full subcomplex in the nested set complex $N_{B(P,n)}$  on the vertex set $\{\{i\},\{1,2,\ldots,k\}|1\leq i\leq n+1, 2\leq k\leq n\}$ (cf. the proof of Lemma~\ref{FlagNestMassey}). 
\end{remark}

\begin{theorem}\label{DirFamNontrivial}
The family $\mathcal P_{Mas}$ is a special geometric direct family with nontrivial Massey products.
\end{theorem}
\begin{proof}
From Proposition~\ref{NestDirect} and Theorem~\ref{boundary} it follows that  $\mathcal P_{Mas}$ is a geometric direct family of polytopes.
Then, using Lemma~\ref{FlagNestMassey}, we immediately deduce that $\mathcal P_{Mas}$ satisfies the condition concerning the nontrivial Massey products from Definition~\ref{DFPM}. Finally, the family $\mathcal P_{Mas}$ is special due to Theorem~\ref{MasAllOrder}, which finishes the proof.
\end{proof}


\begin{definition}\label{length}
We say that a space $X$ has \emph{length} $l(X)\geq k$ with respect to a given spectral sequence for the path loop fibration $\Omega X\to PX\to X$, if its differential $d_k$ is nontrivial. 
\end{definition}

We are going to use the well known statements about the relation between Massey operations in $H^*(X)$ and the values of higher differentials in Eilenberg--Moore spectral sequence for $X$, in order to get lower bounds for  $l_{EM}(\zp)$ for moment-angle manifolds $\zp$ of polytopes from special geometric direct families of polytopes with nontrivial Massey products. 

\begin{theorem}\label{MasseyEM}
Suppose $\mathcal F$ is a special GDFP, $P\in\mathcal F_n$ for $n\geq 2$.
Then all differentials $d_{r}$ for $r\leq n-1$ in Eilenberg--Moore spectral sequence for the path loop fibration for the moment-angle manifold $\zp$ are trivial and its Eilenberg--Moore length $l_{EM}(\zp)\geq n-1$.
\end{theorem}
\begin{proof}
This follows from the well known result, see~\cite[Theorem 8.31]{McCleary}, which states that if there exists a defined Massey product $\langle\alpha_{1},\ldots,\alpha_{n}\rangle$ in $H^*(X)$, then it belongs to the kernel of cohomology suspension homomorphism $\Omega^*$, since it equals (up to sign) the boundary $d_{n-1}(\alpha_{1},\ldots,\alpha_n)$ of an element from the $E_{n-1}$ term of Eilenberg--Moore spectral sequence (converging to $H^*(\Omega X)$), corresponding to the element $[a_{1}|\ldots|a_{n}]$ ($[a_{i}]=\alpha_i$) of the cobar complex. The existence of a strictly defined nontrivial  $k$-fold Massey product $\langle\alpha_{1},\ldots,\alpha_{k}\rangle$ in $H^*(\zp)$ for each $2\leq k\leq n$ finishes the proof.   
\end{proof}

\begin{propos}\label{MasseyNestohedra}
The following statements hold.
\begin{itemize}
\item[(1)] Let $P_{B}$ be a nestohedron on a connected building set $B$ on the vertex set $[n+1]$. Then if $P_{B|_S}$ for some element $S\in B\backslash [n+1]$ has a nontrivial strictly defined Massey product of order $k$ in cohomology of its moment-angle manifold, then the same holds for the polytope $P_B$;
\item[(2)] Let $P_{B}$ be a flag nestohedron on a connected building set $B$ on the vertex set $[n+1]$. Then if $P_{B/S}$ for some element $S\in B\backslash [n+1]$ has a nontrivial strictly defined Massey product of order $k$ in cohomology of its moment-angle manifold, then the same holds for the polytope $P_B$.
\end{itemize}
\end{propos}
\begin{proof}
Statement (1) is a direct consequence of Proposition~\ref{CohomEpi} (2). Statement (2) follows from Theorem~\ref{FlagCriterion} and Theorem~\ref{zkcoh}, since, due to Lemma~\ref{NestBound}, the nestohedron $P_{B/S}$ is a face of the nestohedron $P_B$.
\end{proof}

\begin{coro}\label{GAseqMas}
Let $\mathcal P$ be one of the classical families of graph-associahedra: associahedra, permutohedra, cyclohedra, or stellahedra. Then $\mathcal P=\{P^n|\,n\geq 0\}$ is a sequence of flag nestohedra such that there exists a nontrivial strictly defined Massey product of order $k$ in $H^*(\mathcal Z_{P^n})$, where $n\geq 2$ for $k=2,3$. Furthermore, existence of a nontrivial $r$-fold Massey product in $H^*(\mathcal Z_{P^n})$ implies existence of a nontrivial $r$-fold Massey product in $H^*(\mathcal Z_{P^l})$ for each $r\geq 3$ and $l>n$.
\end{coro}
\begin{proof}
In~\cite{L2,L3} it was proved that in cohomology of moment-angle manifolds over $As^3$, $Pe^3$, $Cy^3$, and $St^3$ there exists a nontrivial strictly defined Massey product. For any $P^{n}=P_B\in\mathcal P, n\geq 4$ there exists an element $S\in B$ on 4 vertices such that $P_{B|_S}$ coincides with one of the 3-dimensional graph-associahedra, considered above: $As^3$, $Pe^3$, $Cy^3$, or $St^3$. Thus, Proposition~\ref{MasseyNestohedra} implies our statement.
\end{proof}

\begin{prob}\label{probMasseyFamilies}
Can a family $\mathcal P$, considered in the previous statement, be a direct family of polytopes with nontrivial Massey products?
\end{prob}

Next, let us introduce the following notion.

\begin{definition}
A sequence of connected graphs 
$$
\Gamma=\{\Gamma_{n}|\,n\geq 0,\Gamma_{n}\,\text{ -- graph on the vertex set }[n+1]=\{1,2,\ldots,n+1\}\}
$$ 
will be called a \emph{direct family of graphs} (DFG), if for any $r$ and any $n>r$ there exists a subset of vertices $S\subset [n+1]$ such that the induced subgraph $\Gamma_{n}|_{S}$ is isomorphic to $\Gamma_{r}$, and the data $\{\Gamma_{n}, k_{r}^{n}\}$ yields a direct family of graphs with respect to graph embeddings $k_{r}^{n}:\,\Gamma_{r}\hookrightarrow\Gamma_{n}$, induced by embeddings $S\subset [n+1]$ (here, we identified $\Gamma_{n}|_S$ with $\Gamma_{r}$).
\end{definition}

Note that the sequences of graphs for the families of associahedra, stellahedra, cyclohedra, and permutohedra are direct families of graphs. It follows from the previous definitions that $\mathcal P_{\Gamma}=\{P_{\Gamma_n}|\,n\geq 0\}$ is a geometric direct family of polytopes whenever $\Gamma$ is a direct family of graphs. Therefore, it is natural to state the following problem.

\begin{prob}
Does there exist a family of graphs $\{\Gamma_{n}|\,n\geq 0\}$ such that there is a defined nontrivial Massey product of order $k$ in $H^*(\mathcal Z_{\Gamma_n})$ with $k\to\infty$ as $n\to\infty$?
\end{prob}

Polytopes $P^n_{Mas}\in\mathcal P_{Mas}$ are not graph-associahedra for any $n\geq 3$. Therefore, one can ask the next question which is a generalization of the previous problem.

\begin{prob}
Is there a DFG $\Gamma$ such that $\mathcal P_{\Gamma}$ is a geometric direct family of polytopes with nontrivial Massey products?
\end{prob}


We end this section with a generalization of Lemma~\ref{FlagNestMassey}.

\begin{definition}\label{FamilyF}
Consider a family of flag nestohedra on connected building sets that can be obtained from the polytopes of the $\mathcal P_{Mas}$ family by application of substitution of building sets operation. Let us denote the resulting nestohedra family by $\mathcal F_{Mas}$.
\end{definition}

\begin{remark}
Note that the set $\mathcal F_{Mas}$ has finitely many elements (pairwisely combinatorially different) in each dimension, and therefore, it is indeed a family.
\end{remark}

The notation from the previous definition is justified by the next statement.

\begin{theorem}\label{NestMassey}
For any $l\geq 2, r\geq 1$ and any finite set of positive integers $S=\{n_{i}\geq 2|\,1\leq i\leq r\}$ there exists an element $P(S)=P_{B(S)}\in\mathcal F_{Mas}$ and an $l$-connected moment-angle manifold $M(l,S)$ over a multiwedge of the polytope $P(S)$ such that there exists a strictly defined nontrivial $n$-fold Massey product in $H^*(M(l,S))$ for all $n\in S$.
\end{theorem}
\begin{proof}
Consider the case $l=2$. By Lemma~\ref{FlagNestMassey}, for connected building sets $B(i)=B(P,n_{i}),1\leq i\leq r$ with $n_{i}\in S$ we have: there exists a strictly defined nontrivial Massey product of order $n_i$ in $H^*(\mathcal Z_{P_{B(i)}})$. Then we can define a connected building set $B_{1}(S)=B(B(1),\ldots,B(r))$ as a result of a substitution of building sets (cf. Construction~\ref{buildsetsubst}), where $B=B(P,r-1)$. We get that $P_{B_{1}(S)}$ is a flag nestohedron, combinatorially equivalent to a product of nestohedra $P^{r-1}_{Mas}\times P^{n_1}_{Mas}\times\ldots\times P^{n_r}_{Mas}$ (cf.~\cite[Lemma 1.5.20]{TT}) and, moreover, $P_{B_{1}(S)}\in\mathcal F$ with $M(2,S)=\mathcal Z_{P_{B_{1}(S)}}$ (note that every moment-angle manifold is 2-connected).

Now, suppose $l\geq 3$. Following~\cite[Definition 3.5]{L3}, let us denote by $P(n,s)=P^{n}_{Mas}(J_{n,s})$ the corresponding multiwedge, see Construction~\ref{simpmultwedge}. Then for $s=[\frac{l+1}{2}]$ the nerve complex of the polytope $P(n,s)$ is $[\frac{l+1}{2}]$-connected, and therefore, a generalized moment-angle manifold $M(l,n)=\mathcal Z_{P^{n}_{Mas}}^{J_{n,s}}\cong\mathcal Z_{P(n,s)}$ is $l$-connected (see~\cite[Proposition 7.34.2]{BP04}). Due to~\cite[Theorem 3.6]{L3}, there exists a strictly defined nontrivial $n$-fold Massey product in $H^*(M(l,n))$.

Consider a product of the above defined building sets in the ring of building sets, cf.~\cite[\S 1.7]{TT}
$$
B_{2}(S)=B(1)\cdot\ldots\cdot B(r),
$$
and an ordered set of positive integers:
$$
J(l,S)=(J_{n_{1},s},\ldots,J_{n_{r},s}).
$$
Then set
$$
M(l,S)=\mathcal Z_{P_{B_{2}(S)}}^{J(l,S)}\cong\prod\limits_{i=1}^{r}\mathcal Z_{P(n_{i},s)}.
$$
The latter generalized moment-angle manifold is $l$-connected as a product of $l$-conncted manifolds and has strictly defined nontrivial Massey product of order $n$ in $H^*(M(l,S))$ for all $n\in S$.

Finally, we have: $B_{1}(S)\in\mathcal F_{Mas}$ for all nonempty finite sets of subsets $S$ from $\mathbb{N}$ with elements greater than one. Applying~\cite[Proposition 1.5.23]{TT}, we can construct a {\emph{connected}} building set $B(S)$, which is a result of iterated application of a substitution of building sets operation to polytopes from the $\mathcal P_{Mas}$ family, thus $P_{B_{2}(S)}$ is combinatorially equivalent to $P_{B(S)}$, which finishes the proof.
\end{proof}


\section{Lusternik--Schnirelmann category and Milnor spectral sequence}

In this section, following~\cite[\S5]{B1969}, we give an exposition of the theory of Milnor spectral sequence and prove the main theorem on the structure of the $E_1$ and $E_2$ terms structure, as well as strong convergence of this spectral sequence. For convenience of the further discussion, we start with the well known Milnor's results on the theory of fiber bundles~\cite{Milnor1956}. For more details on the constructions we introduce below and the proofs of the structure theorems about Milnor filtration and spectral sequence, we refer the reader to the work of Buchstaber~\cite{B1969}, and the papers by Milnor~\cite{Milnor1956}, Dold and Lashof~\cite{DoLa}. 

We denote by $W$ the category of pointed cellular spaces, whose $n$-dimensional skeleta are finite complexes.
Recall that a $\infty$-universal bundle is a principal fibration with a total space having trivial homotopy groups.
 
\begin{theorem}
For any $X\in W$ there exists a topological group $G_{X}\in W$ with cellular multiplication and such that $X$ is a base space of a $\infty$-universal $G_X$-fibration $p_{X}\colon E_{X}\to X$.
\end{theorem}

Any continuous mapping $f\colon X\to Y$ is associated to a group homomorphism $\tilde{f}\colon G_{X}\to G_{Y}$.

\begin{remark}
It can be shown that the group $G_X$ is homotopy equivalent to the loop space $\Omega X$.
\end{remark}

\begin{theorem}
Any principal bundle $E\to X$ with a structure group $G$ is induced by a $\infty$-universal bundle $E_{X}\to X$ with a group $G_X$ by means of a certain continuous homomorphism $h\colon G_{X}\to G$.
\end{theorem}

Starting with the group $G_X$ one can construct a $\infty$-universal $G_X$-bundle $\tilde{E}_{X}\to\tilde{X}$, see~\cite{Milnor1956}, such that $\tilde{E}_X$ and $\tilde{X}$ are both cellular spaces and $\tilde{X}$ is homotopy equivalent to the space $X$. In the space $\tilde{X}$ there exists a canonical filtration such that any continuous mapping $f\colon X\to Y$ is associated to a mapping $\bar{f}\colon\tilde{X}\to\tilde{Y}$, consistent with the filtrations.

This allows one to define 
a functorial spectral sequence
on the category $W$, the so called \emph{Milnor spectral sequence}.

Now, we give a construction of a $\infty$-universal bundle $\tilde{E}_{X}\to\tilde{X}$, following Dold and Lashof~\cite{DoLa}. A small change that we make in their construction allows us to give a simple proof of multiplicativity of Milnor spectral sequence.

A universal bundle $\tilde{E}_{X}\to\tilde{X}$ is defined to be a direct limit of a sequence of principal $G_X$-bundles $\{p_{n}\colon E_{n}\to X_{n}\}_{n\geq 0}^{\infty}$, $\tilde{E}_{X}=\varinjlim\,E_{n}$, $\tilde{X}=\varinjlim\,X_{n}$.

The bundles $p_{n}\colon E_{n}\to X_{n}$ are constructed consecutively. For $n=0$ we set $E_{0}=G_{X}$, $X_{0}=*$ is a point. In $E_0$ we can take $e$ as a based point, that is, the unit of the group $G_X$. Suppose, a principal bundle $p_{n-1}\colon E_{n-1}\to X_{n-1}$ is already constructed. Let us denote by $\mu_{n-1}\colon E_{n-1}\times G_{X}\to E_{n-1}$ the structure group $G_X$ action on the total space $E_{n-1}$ of the principal bundle. By inductive assumption, we may suppose that $\mu_{n-1}$ is a cellular map. Consider the diagram
$$
\begin{CD}
  E_{n-1} @<\mu_{n-1}<< E_{n-1}\times G_{X} @>>> CE_{n-1}\times G_{X}\\
  @VVp_{n-1} V\hspace{-0.2em} @VVq V\hspace{-0.2em} @VVq V\\
  X_{n-1} @<p_{n-1}<<E_{n-1} @>>> CE_{n-1},
\end{CD}\eqno 
$$
where $CE_{n-1}$ is a based cone over $E_{n-1}$ and $q$ is a natural projection map. Gluing the space $CE_{n-1}\times G_X$ to $E_{n-1}$ by the map $\mu_{n-1}$ and gluing the space $CE_{n-1}$ to $X_{n-1}$ by the map $p_{n-1}$, we obtain the spaces $E_{n}$ and $B_{n}$, respectively, alongside with a projection map $p_{n}\colon E_{n}\to B_{n}$. The space $E_n$ is the space of triples of the type $(y,t,g)$, where $y\in E_{n-1}$, $t\in [0,1]$, $g\in G_X$, in which the following identifications are made:
$$
(y_{1},0,g)\sim (y_{2},0,g);\,(y,1,g)\sim (yg,1,e);\,(*,t_{1},g)\sim (*,t_{2},g)
$$
for all $y,y_{1},y_{2}\in E_{n-1}$, $g\in G_X$, and $t_{1},t_{2}\in [0,1]$, where we denote $yg=\mu_{n-1}(y,g)$. For the base point of the space $E_{n}$ we can take the point, which appears after identifying all the points of the type $(*,t,e)$. The action $\mu_{n}\colon E_{n}\times G_{X}\to E_{n}$ of the group $G_X$ on such a represented space $E_{n}$ can be described by formula 
$$
\mu_{n}(y,g)=\mu_{n}((y_{1},t,g_{1}),g)=(y_{1},t,g_{1}g).
$$
It is easy to check that $p_{n}\colon E_{n}\to B_{n}$ is a principal $G_X$-bundle.

In this way one can construct a direct sequence of principal $G_X$-bundles
$$
\begin{CD}
  E_{0} @>i_{0}>> E_{1} @>i_{1}>> \cdots @>>> E_{n} @>>>\cdots\\
  @VVV\hspace{-0.2em} @VVV\hspace{-0.2em} @VVV\hspace{-0.2em} @VVV\hspace{-0.2em} @VVV\\
  X_{0} @>\bar{i}_{0}>> X_{1} @>\bar{i}_{1}>>\cdots @>>>X_{n} @>>>\cdots,
\end{CD}\eqno 
$$
where the mapping $i_{k}\colon E_{k}\to E_{k+1}$ can be described by formula $i_{k}(y)=(y,1,e)$. The mappings $i_{k}$ commute with the $G_{X}$-action, thus, turning to direct limits, we get a principal $G_X$-bundle $\varinjlim\,E_{n}\to\varinjlim\,X_{n}$. (At this point it is essential that for each $n$ the action $\mu_{n}\colon E_{n}\times G_{X}\to E_{n}$ is cellular, cf.~\cite{Milnor1956}). Set $\tilde{E}_{X}=\lim\,E_{n}$, $\tilde{X}=\lim\,X_n$. The mapping $i_{k}\colon E_{k}\to E_{k+1}$ can be extended to a mapping of the reduced cone $J_{k}\colon CE_{k}\to E_{k+1}$, $J_{k}(y,t)=(y,t,e)$, that is, the subspace $i_{k}(E_k)$ is contracting by $E_{k+1}$ to a point, and therefore, the space $\tilde{E}$ is aspherical.
Thus, the fibration $\tilde{E}_{X}\to\tilde{X}$ is a $\infty$-universal bundle for the group $G_X$.

\begin{definition}
{\emph{Milnor filtration}} of a space $X\in W$ is the above constructed filtration $\{X_n\}$ in the space $\tilde{X}$, which is homotopy equivalent to the space $X$.
\end{definition}

Let us list the main properties of the Milnor filtration that follow directly from its construction; all homeomorpisms (notation: $=$) and homotopy equivalences (notation: $\simeq$) below are functorial with respect to the maps, induced by continuous mappings $f\colon X\to Y$ (here, $SX$ denotes a reduced suspension over $X$):
\begin{itemize}
\item[(a)] $X_{1}=SG_{X}$;
\item[(b)] $E_{n}/E_{n-1}=(CE_{n-1}\times G_{X})/(E_{n-1}\times G_{X})$, $X_{n}/X_{n-1}=SE_{n-1}$;
\item[(c)] $E_{n}\simeq G_{X}\wedge SE_{n-1}\simeq G_{X}\wedge \underbrace{SG_{X}\wedge\ldots\wedge SG_{X}}_{n\text{ times}}$;
\item[(d)] $X_{n+1}/X_{n}\simeq\underbrace{SG_{X}\wedge\ldots\wedge SG_{X}}_{n+1\text{ times}}$.
\end{itemize}

\begin{theorem}\label{mainMilnor}
For any space $X\in W$ there exist the mappings $\varphi\colon\tilde{E}_{X}\to\tilde{E}_{X}\times\tilde{E}_{X}$ and $\bar{\varphi}\colon\tilde{X}\to\tilde{X}\times\tilde{X}$ such that:
\begin{itemize}
\item[(1)] The mappings $\varphi$ and $\bar{\varphi}$, alongside with the diagonal homomorphism $G_{X}\to G_{X}\times G_{X}$, yield a map of a principal $G_X$-bundle $(\tilde{E}_{X},p,\tilde{X})$ to a principal $G_{X}\times G_{X}$-bundle $(\tilde{E}_{X}\times\tilde{E}_{X},p\times p,\tilde{X}\times\tilde{X})$;
\item[(2)] The mapping $\bar{\varphi}$ is homotopy equivalent to the diagonal map $\tilde{X}\to\tilde{X}\times\tilde{X}$;
\item[(3)] The mappings $\varphi$ and $\bar{\varphi}$ are functorial with respect to the maps $f\colon X\to Y$;
\item[(4)] $\varphi(E_s)\subset\cup_{s^{'}+s^{''}=s}\,E_{s^{'}}\times E_{s^{''}}$; $\varphi(X_s)\subset\cup_{s^{'}+s^{''}=s}\,X_{s^{'}}\times X_{s^{''}}$.
\end{itemize}
\end{theorem}
\begin{proof}
Denote by $f_{i}(\tau),i=1,2$ the continuous maps, given by the following formulae
$$
f_{1}(\tau)=\begin{cases}
0,&\text{if $0\leq\tau\leq\frac{1}{2}$;}\\
2\tau-1,&\text{if $\frac{1}{2}\leq\tau\leq 1$,}
\end{cases}
$$
$$
f_{2}(\tau)=\begin{cases}
2\tau,&\text{if $0\leq\tau\leq\frac{1}{2}$;}\\
1,&\text{if $\frac{1}{2}\leq\tau\leq 1$.}
\end{cases}
$$
Let us define a mapping $\varepsilon_{i}\colon\tilde{E}_{X}\to\tilde{E}_{X}$ for $i=1,2$ in the following way. On the subspace $E_{0}=G_{X}\subset\tilde{E}_{X}$ the mapping $\varepsilon_{i}$ is identical; for any point of the type $(y,t,g)\in E_{s}\subset\tilde{E}_{X}$, where $y\in E_{s-1}$, $t\in [0,1]$, $g\in G_X$, we set $\varepsilon_{i}(y,t,g)=(\varepsilon_{i}(y),f_{i}(t),g)$. The mapping $\varepsilon_i$ is then defined correctly, it is continuous and equivariant with respect to the group $G_X$ action. Therefore, it determines a continuous mapping $\bar{\varepsilon}_{i}\colon\tilde{X}\to\tilde{X}$. Set 
$$  
\varphi(x)=(\varepsilon_{1}(x),\varepsilon_{2}(x)),\,\bar{\varphi}(x)=(\bar{\varepsilon}_{1}(x),\bar{\varepsilon}_{2}(x)).
$$
The statements (1)-(3) follow directly from the construction of the mappings $\varphi$ and $\bar{\varphi}$. 

Let us show that $\varphi(E_s)\subset\cup_{s^{'}+s^{''}=s}\,E_{s^{'}}\times E_{s^{''}}$, that is, for any point $x\in E_s$ there exist nonnegative integers $s^{'}$ and $s^{''}$ such that $s^{'}+s^{''}=s$ and that $\varepsilon_{1}(x)\in E_{s^{'}}$,  $\varepsilon_{2}(x)\in E_{s^{''}}$. For $s=0$ this is obvious. Suppose that we already proved this statement for certain $s=s_0$. Let $x\in E_1$ be an arbitrary point. Then $x=(y,t,g)$, where $y\in E_{s_0}$, and therefore, $\varepsilon_{1}(y)\in E_{s^{'}_0}$ and $\varepsilon_{2}(y)\in E_{s^{''}_0}$, where $s^{'}_{0}+s^{''}_{0}=s_{0}$.

We have $\varepsilon_{1}(x)=(\varepsilon_{1}(y),f_{1}(t),g)$ and $\varepsilon_{2}(x)=(\varepsilon_{2}(y),f_{2}(t),g)$. For any $t\in [0,1]$ one has: either $f_{1}(t)=0$, or $f_{2}(t)=1$. In the latter case we get $\varepsilon_{1}(x)\in E_{s^{'}_{0}+1}$ and $\varepsilon_{2}(x)=(\varepsilon_{2}(y),1,g)=(\varepsilon_{2}(y)g,1,e)\in E_{s^{''}_{0}+1}$, and the statement is proved.
In the former case we get $\varepsilon_{2}(x)\in E_{s^{''}_{0}+1}$ and $\varepsilon_{1}(x)=(\varepsilon_{1}(y),0,g)=(*,0,g)=(*,1,g)=(*g,1,e)\in E_{0}\subset E_{s^{'}_{0}}$, and the statement is proved. 

Since the mapping $\varphi$ is equivariant, we obtain $\bar{\varphi}(X_s)\subset\cup_{s^{'}+s^{''}=s}\,X_{s^{'}}\times X_{s^{''}}$, which finishes the proof of the theorem.
\end{proof}

Note that the proof of statement (4) of Theorem~\ref{mainMilnor} is the only place where we needed the identification $(*,t_{1},g)\sim (*,t_{2},g)$ in the space $E_s$, which is missed in the original construction of Dold and Lashof~\cite{DoLa}.

Let $\begin{CD} p\colon E @>F>> X\end{CD}$ be any locally trivial fibration. Milnor filtration $\{X_n\}$ in the space $X$ determines a filtration $\{Y_n\}$ in the space $E$, where $Y_{n}=p^{-1}(X_n)$. It follows from property (c) of Milnor filtration that $Y_{n}/Y_{n-1}=((X_{n}/X_{n-1})\times F)/F$.

\begin{definition}
We call a \emph{Milnor spectral sequence} for the fibration $\begin{CD} p\colon E @>F>> X\end{CD}$ a spectral sequence of the filtration $\{Y_{n}=p^{-1}(X_n)\}_{n\geq 0}^{\infty}$, where $\{X_n\}$ is a Milnor filtration of the space $X$.
\end{definition}

Applying Theorem~\ref{mainMilnor}, we get the next theorem (see~\cite[\S15]{Dold1966}) in a standard way.

\begin{theorem}\label{multspectMilnor}
For any multiplicative cohomology theory $h^*$ the Milnor spectral sequence for a locally trivial fibration $\begin{CD} p\colon E @>F>> X\end{CD}$ is multiplicative, that is, its terms $E_r$ are graded skew-commutative rings, differentials $d_{r}$ are derivations. Furthermore, if the spectral sequence strongly converges, then multiplication in the $E_{\infty}$ term is induced by multiplication in the ring $h^*(E)$.
\end{theorem}

It was shown in~\cite[Example 7.43]{BP04} that the polyhedral product operation, applied to a direct sequence of simplicial complexes, allows one to obtain a construction of a filtration in the universal principal $G$-bundle $EG\to BG$, where $G$ is a topological group, which generalizes the Milnor flitration construction that we introduced above. The latter one arises from a sequence of simplices. We give this construction below. 

\begin{constr}
Suppose we have a direct sequence of simplicial complexes
$$
K_{0}\subset K_{1}\subset\ldots\subset K_{n}\subset\ldots
$$
in which a complex $K_{n}$ is $n$-neighbourly. Then the direct limit of the filtration of the corresponding polyhedral products $M=\varinjlim\mathcal\,(\cone(G),G)^{K_{n}}$ is $n$-connected for all $n$, and therefore, is a contractible cellular space, on which the topological group $G$ as acting freely. Therefore, we get a filtration in the universal principal $G$-bundle $EG\to BG$:
$$
\begin{CD}
  (\cone(G),G)^{\Delta^{0}} @>>>(\cone(G),G)^{\Delta^{1}} @>>>\ldots @>>>(\cone(G),G)^{\Delta^{n}}\\
  @VVV\hspace{-0.2em} @VVV\hspace{-0.2em} @VVV\hspace{-0.2em} @VVV\hspace{-0.2em} @.\\
  (\cone(G),G)^{\Delta^{0}}/G @>>>(\cone(G),G)^{\Delta^{1}}/G @>>>\ldots @>>>(\cone(G),G)^{\Delta^{n}}/G, 
\end{CD}\eqno 
$$ 
\end{constr}

For the geometric direct families of polytopes $\mathcal P=\{P^n\,|\,n\geq 0\}$, which we studied earlier in this work, we also obtain a direct sequence of their moment-angle manifolds 
$$
\ast\subset S^{3}\to\ldots\to\mathcal Z_{P^n}\to\mathcal Z_{P^{n+1}}\to\ldots
$$
with a free action of the diagonal circle $G=S^1$, when $n>0$. Note that all the manifolds $\mathcal Z_{P^n}$ as well as their direct limit $M=\varinjlim\mathcal Z_{P^n}$ are 2-connected. Observe that $H^{3}(M;\Z)\neq 0$, and, in particular, the space $M$ is not contractible. We are going to study the topology of the quotient spaces $\mathcal Z_{P^{n+1}}/\mathcal Z_{P^n}$ and the topology of the space $M$ in our subsequent publications.

 
Now, we turn to a more detailed study of the Milnor spectral sequence. For our purposes it is sufficient to assume in what follows that $h^*$ is a multiplicative cohomology theory with coefficient ring $\ko=h^{*}(\ast)$ being a module over a field. 

Our results are valid also in the case when the considered spaces from the category $W$ do not have torsion in ordinary cohomology, since in the latter case the cohomology ring of such a space is a free module over the coefficient ring $\ko$.

For the cohomology theory and the spaces to which we restrict ourselves here, the multiplication homomorphism 
$$
\mu\colon\tilde{h}^{*}(X)\hat{\otimes}_{\ko}\tilde{h}^{*}(Y)\to\tilde{h}^{*}(X\wedge Y)
$$
is an isomorphism. Here, for infitely dimensional cellular spaces $X=\cup_{n=0}^{\infty}\,X^{n}$ and $Y=\cup_{m=0}^{\infty}\,Y^{m}\in W$, by definition, we set: $\tilde{h}^{*}(X)=\varprojlim\limits_{n}\,\tilde{h}^{*}(X^{n})$ and $\tilde{h}^{*}(X\wedge Y)=\varprojlim\limits_{n,m}\,\tilde{h}^{*}(X^{n}\wedge Y^{m})$. Note that in the adjunctioned tensor product $\hat{\otimes}$, unlike the case of the ordinary tensor product $\otimes$, the infinite series are involved.

In the applications of Milnor spectral sequence one can take as $h^*$ the classical cohomology theorie, or $K$-theory, or complex cobordism.

In Theorem~\ref{propertiesspecMilnor} for cohomology theory $h^*$ we describe the terms $E_1$ and $E_2$. We also prove strong convergence of Milnor spectral sequence of a locally trivial fibration.

Consider a locally trivial fibration $\begin{CD} p\colon E @>F>> X\end{CD}$, associated to the universal bundle $\begin{CD} p_X\colon E_X @>G_X>>\tilde{X}\end{CD}$ by means of a certain action $\mu\colon G_{X}\times F\to F$ of the group $G_X$ on the fiber $F$. Since, as we noted above, for our cohomology theory $\tilde{h}^{*}$ and our spaces $X,Y\in W$ the multiplication homomorphism 
$$
\mu\colon\tilde{h}^{*}(X)\hat{\otimes}_{\ko}\tilde{h}^{*}(Y)\to\tilde{h}^{*}(X\wedge Y)
$$
is an isomorphism, it can be shown that the ring ${h}^{*}(G_X)$
is a graded commutative Hopf algebra over $\ko$, and the ring ${h}^{*}(F)$
is a graded comodule over the Hopf algebra $h^{*}(G_X)$. Let us set ${h}^{*}(G_X)=A$ and ${h}^{*}(F)=M$. Denote by $\Delta\colon A\to A\hat{\otimes}A$ the comultiplication in the algebra $A$, and denote by $\Delta^{'}\colon M\to M\hat{\otimes}A$ the coaction in the comodule $M$. Then the group $M\hat{\otimes}M$ is a comodule over the Hopf algebra $A\hat{\otimes}A$. In our case the multiplication homomorpism $\varphi\colon M\hat{\otimes}M\to M$ in the ring $M$ is a homomorphism of a $A\hat{\otimes}A$-comodule in a $A$-comodule, and therefore, it induces a homomorphism
$$
\varphi^{*}\colon Cotor_{A}(M,\ko)\hat{\otimes}Cotor_{A}(M,\ko)\to Cotor_{A}(M,\ko),
$$ 
and thus the homology group $Cotor_{A}(M,\ko)$ of the comodule $M$ over $A$ is, in fact, a ring. (For the theory of the functor $Cotor$ we refer the reader to~\cite{EiMoore}).

\begin{theorem}\label{propertiesspecMilnor}
Milnor spectral sequence for a locally trivial fibration $\begin{CD} p\colon E @>F>> X\end{CD}$ has the following properties:
\begin{itemize}
\item[(1)] The term $E_{1}=\sum\,E_{1}^{s,t}$ can be described by formula
$$
E_{1}^{s,t}=M\hat{\otimes}\underbrace{\bar{A}\hat{\otimes}\cdots\hat{\otimes}\bar{A}}_{s\text{ times}},
$$
where $\bar{A}$ is a kernel of the natural augmentation $A\to\ko$;
\item[(2)] The product of elements $b_{0}\otimes a_{1}\otimes\ldots\otimes a_{s}\in E_{1}^{s}$ and $b_{0}^{'}\otimes a_{1}^{'}\otimes\ldots\otimes a_{s^{'}}^{'}\in E_{1}^{s^{'}}$ is equal to 
$$
(-1)^{s\cdot\deg\,b_{0}^{'}}b_{0}\cdot b_{0}^{'}\otimes a_{1}\otimes\ldots\otimes a_{s}\otimes a_{1}^{'}\otimes\ldots\otimes a_{s^{'}}^{'}\in E_{1}^{s+s^{'}};
$$
\item[(3)] Differential $d_{1}^{s}$ is given by
$$
d_{1}^{s}(b_{0}\otimes a_{1}\otimes\ldots\otimes a_{s})=\bar{\Delta}^{'}b_{0}\otimes a_{1}\otimes\ldots\otimes a_{s}+\sum\limits_{i=1}^{s}\,b_{0}\otimes a_{1}\otimes\ldots\otimes\tilde{\Delta}a_{i}\otimes\ldots\otimes a_{s},
$$
where $\bar{\Delta}^{'}(b_0)=\Delta^{'}b_{0}-b_{0}\otimes 1$, $\tilde{\Delta}a_{i}=\Delta\,a_{i}-a_{i}\otimes 1-1\otimes a_{i}$;
\item[(4)] The ring $E_{2}^{*,*}$ is naturally isomorphic to the ring $Cotor_{A}(M,\ko)$, where $A=h^*(G_X;\ko)$ and $M=h^*(F;\ko)$;
\item[(5)] The spectral sequence is strongly convergent, and therefore, multiplication in the term $E_{\infty}^{*,*}$ is induced by multiplication in the ring $h^*(E;\ko)$.
\end{itemize}
\end{theorem}
\begin{proof}
Statement (1) of the theorem follows directly from properties (b) and (c) of Milnor filtration in the base space $X$ and the fact that $p^{-1}(X_n)/p^{-1}(X_{n-1})=((X_{n}/X_{n-1})\times F)/F$. Statement (2) follows from the explicit formula for the mapping $\bar{\varphi}$, Theorem~\ref{mainMilnor}, and the definition of multiplication in the ring $E_{1}^{*,*}$. The term $E_{1}^{s,*}$ is multiplicatively generated by elements of the group $h^*(F)=E_{1}^{0,*}$ and $1\otimes\tilde{h}^{*}(G_X)\subset h^*(F)\hat{\otimes}\tilde{h}^{*}(G_X)=E_{1}^{1,*}$. Since the differential $d_{1}^{s}$ satisfies the ``Leibniz rule'', it suffices to check the formula from statement (3) on the generators. For the generators it is a formal consequence of the definition of the differential $d_1$. The term $(E_{1}^{s,*},d_{1}^{s})$ alongside with multiplication are precisely the standard cobar construction of the comodule $M$ over the coalgebra $A$, and therefore, statement (4) is a well known algebraic fact about coincidence of the two multiplications in the ring $Cotor_{A}(M,\ko)$, where the first multiplication is defined via the cobar construction, and the second one is defined via the homomorphism of comodules
$\varphi\colon M\hat{\otimes}M\to M$, cf.~\cite[\S2.2]{Adams1961}. 

Finally, statement (5) follows from the following general result on strong convergence of a spectral sequence. Suppose $X\in W$ and $\{X_n\}_{n=0}^{\infty}$ is a filtration in the space $X=\varinjlim\,X_n$. A spectral sequence $\{E_{r}^{s,t},d_{r}^{s,t}\}$ of the filtration $\{X_n\}$ strongly converges in the cohomology theory $h^*$. This finishes the proof. 
\end{proof}


Our next goal is to obtain a result about vanishing of certain differentials in Milnor spectral sequence for moment-angle manifolds. 
To do this, we shall use a classical notion of Lusternik--Schnirelmann category of a topological space. For more details on Lusternik--Schnirelmann category we refer, for example, to the monograph~\cite{C-L-O-T}, which is dedicated to it. 

We are going to compute this important homotopy invariant for moment-angle manifolds over 2-truncated cubes from the $\mathcal Q$ family, cf. Definition~\ref{2truncMassey}, and then link this result with degeneration of Milnor spectral sequences for those moment-angle manifolds.

\begin{definition}
A covering of a topological space $X$ is called {\emph{categorical}}, if each element of the covering is a set, which is open and contractible in $X$, that is, the embedding of such a set into $X$ is homotopy trivial.\\
{\emph{Lusternik--Schnirelmann category}} (or, simply, a {\emph{category}}) $cat(X)$ of a space $X$ is the least integer $k$ such that $X$ has a categorical covering, which consists of exactly $k+1$ sets.
\end{definition}

Below we list the basic properties of the category.

\begin{itemize}
\item $cat(X)=0$ $\Longleftrightarrow$ $X$ is contractible;
\item $cat(X)\leq 1$ $\Longleftrightarrow$ $X$ is a co-$H$-space;
\item A smooth closed manifold $M$ having $cat(M)=1$ is homeomorphic to a sphere;
\item $cup(X)\leq cat(X)$, where \emph{cohomology length} $cup(X)$ of a space $X$ is defined to be the maximal number of positive degree elements from $H^*(X;\Q)$, having a nonzero product.
\end{itemize}

We recall that, cf. Definition~\ref{length}, a topological space $X$ has \emph{length} $l(X)\geq k$ with respect to a given spectral sequence for the path loop fibration $\Omega X\to PX\to X$, if differential $d_k$ is nontrivial. 

\begin{theorem}\label{lengthMilnor}
Consider the $n$-dimensional 2-truncated cube $Q\in\mathcal Q_n$ for $n\geq 2$.
Then $cup(\mathcal Z_{Q})=cat(\mathcal Z_{Q})=n$.
Milnor spectral sequence for the moment-angle manifold $\mathcal Z_{Q}$ degenerates in the $E_{n+1}$ term and $l_{M}(\mathcal Z_{Q})\leq n$.
\end{theorem}
\begin{proof}
Due to the properties of the category, listed above, to prove the first part of the statement it suffices to show that the first and the last of the inequalities below hold:  
$$
n\leq cup(\mathcal Z_{Q})\leq cat(\mathcal Z_{Q})\leq n.
$$
To prove the left inequality, let us consider the following cohomology classes from $H^*(\mathcal Z_{Q})$ (here, we use the notation from Definition~\ref{2truncMassey}):
$$
\alpha_{1}=[v_{1}u_{n+1}],\alpha_{2}=[v_{2}u_{n+2}u_{n+3}\cdots u_{2n}u_{1,n+2}],
$$
$$
\alpha_{k}=[v_{k}u_{k-1,n+k}u_{k-2,n+k}\cdots u_{1,n+k}],\,3\leq k\leq n-1,
$$
$$
\alpha_{n}=[v_{n}u_{n-1,2n}u_{n-2,2n}\cdots u_{2,2n}],
$$
where (see Theorem~\ref{zkcoh}) $du_{i,j}=v_{i,j}$ in the differential graded algebra $R^*(Q)$.
Observe that, by Theorem~\ref{zkcoh}, the product $\alpha_{1}\alpha_{2}\cdots\alpha_{n}$ is equal (up to sign) the generator of the top degree cohomology group $H^{m(n)+n}(\mathcal Z_{Q})$. Since $H^{m(n)+n}(\mathcal Z_{Q})\cong\mathbb Z$, the desired inequality follows.

The right inequality above holds, since $cat(\zk)\leq n$ for any simplicial complex $K$ of dimension $n-1$ due to~\cite[Lemma 3.3]{BG}. Thus, $cup(\mathcal Z_{Q})=cat(\mathcal Z_{Q})=n$ for all $n\geq 2$, which finishes the proof of the first part of the statement.

Due to the result of Ginsburg~\cite{Gin} (for modification of one of the steps of the proof see the work of Ganea~\cite{Gan}), all differentials $d_{r}$ for $r>n$ in Milnor spectral sequence, cf. Theorem~\ref{propertiesspecMilnor}, are trivial for spaces of category $n$. 
The theorem is proved.
\end{proof}


\section{Discussion}

There is a categorical notion of a bracket. It is defined in any category, where for every two objects there is an object called their product, for every single object there is an object called its boundary, and for every triple $A$, $B$, $C$ of objects with morphisms $C\to A$ and $C\to B$ there exists an object $D$ called a union of $A$ and $B$ along $C$.

In cohomology theory, we obtain the notion of a triple Massey product and in homotopy groups of spheres we get a triple Toda bracket~\cite{Toda62}.
In cobordism theory, there exists a geometric operation, which gives a closed manifold $M$ if given three closed manifolds, $M_{1},M_{2},M_{3}$, two manifolds $W_{1}$ and $W_{2}$ with boundaries such that $\partial W_{1}=M_{1}\times M_{2}$, $\partial W_{2}=M_{2}\times M_{3}$ and a fixed homeomorphism $\partial W_{1}\cong\partial W_{2}$. The resulting closed manifold $M$ is a union of $W_{1}\times M_{3}$ and $M_{1}\times W_{2}$ along their common boundary.

\begin{prob}
Find a construction of a higher order bracket in general categorical terms due to Eilenberg and Mac Lane.
\end{prob} 

One of the first applications of Massey products was made in the paper of Massey and Uehara~\cite{M-U} (1957) and consisted in the proof of Jacobi identity for Whitehead products in $\pi_{*}(X)$. This appeared to be one of the manifestations of the famous Eckmann--Hilton duality between homotopy and cohomology, be which a triple Massey product is associated to a (triple) Toda bracket, see~\cite{Toda62}. The higher Massey products, introduced shortly in the work of Massey~\cite{Mass}, by this duality are associated to the (higher) Toda brackets, introduced by Cohen in his note~\cite{Cohen67}, in which it was stated that elements of stable homotopy groups of spheres are generated by (higher) Toda brackets of Hopf classes. 

In the literature, see, for example~\cite{McCleary,BaTa}, the following statements can be found, for which, however, no rigorous proof has ever been provided:
\begin{itemize}
\item[(1)] Higher differentials in Eilenberg--Moore spectral sequence for the path loop fibration $\Omega X\to PX\to X$ of a connected and simply connected space $X$ are described by (matric) Massey products. An accurate statement, which illustrates this thesis, is our Theorem~\ref{MasseyEM};
\item[(2)] The kernel of the cohomology suspension homomorphism $\Omega^*: H^{*}(X)\to H^{*-1}(\Omega X)$ is generated by matric (alongside with the ordinary) Massey products.
\end{itemize}

Both of these statements and the definition of the matric Massey product were first announced in the note by May~\cite{May1968}.

In this paper we constructed a direct sequence of moment-angle manifolds $\{\mathcal Z_{P^n}\}$, which determines nontrivial strictly defined Massey products in the ring $\varprojlim H^*(\mathcal Z_{P^n})$, see Section 5.
In relation with this, we want to mention here the inverse sequence
$$
M_{0}\leftarrow M_{1}\leftarrow\ldots\leftarrow M_{n}\leftarrow M_{n+1}\leftarrow\ldots
$$ 
of nilmanifolds, where $\pi_{q}\colon M_{n+1}\to M_{n}$ is a principal $S^1$-bundle, see~\cite{B1999}. 

In the beginning of 1970s the first author formulated a conjecture (based on the results of the works~\cite{B1970,Gon1,Gon2,B-Sh}), that each element of the ring $\varinjlim H^*(M_n;\Q)$ can be realized as a nontrivial Massey product of the generators of the group $H^{1}(\varprojlim M_n;\Q)=\Q\oplus\Q$. A partial answer to that conjecture was given in the work of Artelnykh~\cite{Art}, and the full solution was given in the work of~\cite{Million}.  

Finally, let us return back to the problem of relation between homotopy invariants of Lusternik--Schnirelmann type and the differentials in Milnor spectral sequence.

There is a significant number of publications about Lusternik--Schnirelmann category and its generalizations. In algebraic topology there is a well known notion of the genus of a space, introduced by Schwartz~\cite{Sh}, and its applications. 

By the use of homological algebra methods, Halperin and Lemaire~\cite{HL} introduced a new homotopy invariant, $Mcat(X)$. This invariant was defined in terms of the notion of a free tensor model $(T(V),d)$, introduced by the same authors, which is quasiisomorphic to the singular cochain complex $(C^*(X;\ko),d)$ with coefficients in a field $\ko$. It was shown in~\cite{HL} that $cup(X)\leq Mcat(X)\leq cat(X)$ holds for any connected and simply connected space $X$.
In~\cite{HL} one can find the following example: let $X=Sp(2)$ be the symplectic group, then $cup(X)=Mcat(X)=2$, but $cat(X)=3$.

Ginsburg--Ganea theorem, which was used in the previous section, follows from a stronger result, which was obtained in the work of Jessup~\cite{Jess}. A theorem of Jessup states that if $Mcat(X)=n$, where $X$ has a homotopy type of a simply connected cellular space of a finite $\ko$-type, then Milnor spectral sequence for the space $X$ degenerates in the $E_{n+1}$ term. 

Applying~\cite[Proposition 1.6]{HL} that links Milnor filtration for $X$ to degree filtration in the free tensor model $\{(T^{\geq k}(X),d)\}$, Jessup proved that for degeneration of Milnor spectral sequence in $E_{n+1}$ the following condition is sufficient: for any nonnegative $k$, if $a\in T(V)$ is such that $d(a)\in T^{\geq n+k}(V)$, then there exists an element $\tilde{a}\in T^{\geq k}(V)$ with $d(a-\tilde{a})=0$. Jessup gives an explicit algorithm of how to find such an element $\tilde{a}$ for the given element $a$, using the existence of a certain morphism of tensor models, which is equivalent (due to definition of $Mcat(\cdot)$) to the condition $Mcat(X)=n$, cf.~\cite{HL}.

Since $cup(X)\leq Mcat(X)\leq cat(X)$, it follows immediately from Theorem~\ref{lengthMilnor} that $cup(\zp)=Mcat(\zp)=cat(\zp)=n$, where $P=Q^{n}\in\mathcal Q_n$ is the $n$-dimensional 2-truncated cube from the special geometric direct family $\mathcal Q$ of polytopes with nontrivial Massey products. Thus, we immediately obtain the upper bound for the Milnor length $l_{M}(\zp)\leq n$ (Theorem~\ref{lengthMilnor}). 

On the other hand, due to Theorem~\ref{mainMassey}, there exist strictly defined nontrivial Massey products of all the orders from 2 to $n$ in $H^*(\zp)$.
We call a {\emph{Massey length}} $Mas(X)$ of a space $X$ the maximal number $k$, for which there exists a nontrivial $k$-fold Massey product in $H^*(X;\mathbb Q)$. We assume, as before, that a 2-fold Massey product coincides with the usual product in cohomology of a space, up to sign. Therefore, following the above definition, we get $Mas(X)=cup(X)$ for spaces with trivial triple and higher Massey products. Our results on lengths of manifolds $\zp$ with respect to Eilenberg--Moore and Milnor spectral sequences naturally lead us to the next question.  

\begin{prob}\label{ineq}
Let $X$ be a moment-angle manifold. Is it true that $Mas(X)\geq cat(X)$? 
\end{prob}
In general, from the inequality $Mas(X)\geq cat(X)$ the chain of inequalities follows: 
$$
l_{EM}(X)+1\geq Mas(X)\geq cat(X)\geq Mcat(X)\geq l_{M}(X).
$$
Observe that for connected and simply connected CW-complexes of finite type over $\mathbb Q$, which are not moment-angle manifolds, the inequality from Problem~\ref{ineq} may not hold. Indeed, in the work by Iriye and Yano~\cite{IY} a simplicial complex $K$ was constructed such that $\mathbb Q[K]$ is a Golod ring, but $\zk$ is not a co-H-space, and therefore, $cat(\zk)>Mas(\zk)=1$. The above inequalities also do not hold for the manifold $Sp(2)$, see the discussion above.




\begin{thebibliography}{99}

\bibitem{Adams1961}
J. F. Adams, \emph{On the non-existence of elements of Hopf invariant one}, Ann. of Math., {\bf{72}}:1 (1960), 20--104.

\bibitem{And}
E. M. Andreev, \emph{On convex polyhedra in Lobachevskii spaces}, Math. USSR-Sb., {\bf{10}}:3 (1970), 413--440.

\bibitem{Art}
I. V. Artelnykh, \emph{Massey Products and the Bukhshtaber spectral sequence}, Russian Math. Surveys, {\bf{55}}:3 (2000), 559--561.

\bibitem{BaTa} 
I. K. Babenko, I. A. Taimanov, \emph{Massey products in symplectic manifolds}, Sb. Math., {\bf{191}}:8 (2000), 1107--1146.

\bibitem{BBCG10} 
A. Bahri, M. Bendersky, F. R. Cohen, and S. Gitler, \emph{The polyhedral product functor: a method of decomposition for moment-angle complexes, arrangements and related spaces,} Adv. Math. {\bf{225}} (2010), no. 3, 1634--1668. MR 2673742 (2012b:13053)

\bibitem{BBCG14} 
A. Bahri, M. Bendersky, F. R. Cohen, and S. Gitler, \emph{On free loop spaces of toric spaces}, preprint (2014); arXiv:1410.6458v2.

\bibitem{BBCG15} 
A. Bahri, M. Bendersky, F. R. Cohen, and S. Gitler, \emph{Operations on polyhedral products and a new topological construction of infinite families of toric manifolds,} Homology Homotopy Appl. 17 (2015), no. 2, 137--160. MR 3426378; arXiv:1011.0094.

\bibitem{BaskM} 
I. V. Baskakov, \emph{Massey triple products in the cohomology of moment-angle complexes}, Russian Math. Surveys, {\bf{58}}:5 (2003), 1039--1041.

\bibitem{BBP}
I. V. Baskakov, V. M. Buchstaber, T. E. Panov, \emph{Cellular cochain algebras and torus actions}, Russian Math. Surveys, {\bf{59}}:3 (2004), 562--563.

\bibitem{BG} 
P. Beben and J. Grbi\'c, \emph{LS-category of moment-angle manifolds, Massey products, and a generalisation of the Golod property,} Transactions of the American Mathematical Society, 1-21 (Submitted) (2016); arXiv:1601.06297v2.

\bibitem{BJ} 
Alexander Berglund and Michael J\"ollenbeck, \emph{On the Golod property of Stanley--Reisner rings}, J. Algebra \textbf{315}:1 (2007), 249--273.

\bibitem{Bg}
F. A. Bogomolov, \emph{The actions of a circumference on smooth odd-dimensional spheres}, Math. Notes, {\bf{5}}:4 (1969), 241--244.

\bibitem{bo-me06}
Fr\'ed\'eric Bosio and Laurent Meersseman, \emph{Real quadrics in
$\mathbb C^n$, complex manifolds and convex polytopes,} Acta
Math.~\textbf{197} (2006), no.~1, 53--127.

\bibitem{BT}
Raul Bott and Clifford Taubes, \emph{On the self-linking of knots. Topology and physics}. J.Math.Phys. {\bf{35}} (1994), no. 10, 5247--5287.

\bibitem{br-gu96}
Winfried Bruns and Joseph Gubeladze, \emph{Combinatorial
invariance of Stanley--Reisner rings,} Georgian Math.
J.~\textbf{3} (1996), no.~4, 315--318.

\bibitem{B1969}
V. M. Buchstaber. \emph{Cohomology operations in generalized cohomology theories,} PhD thesis (in Russia), Moscow State University, M., 1969.

\bibitem{B1970}
V. M. Buchstaber, \emph{The Chern–Dold character in cobordisms. I}, Math. USSR-Sb., {\bf{12}}:4 (1970), 573--594.

\bibitem{B1999}
V. M. Buchstaber, \emph{Groups of polynomial transformations of a real line, non-formal symplectic manifolds, and the Landweber--Novikov algebra}, Russian Math. Surveys, {\bf{54}}:4 (1999), 837--838.

\bibitem{B} 
V. M. Buchstaber, \emph{Ring of simple polytopes and differential equations}, Proc. Steklov Inst. Math., {\bf{263}} (2008), 13--37.

\bibitem{BE} 
V. M. Buchstaber, N. Yu. Erokhovets, \emph{Polytopes, Fibonacci numbers, Hopf algebras, and quasi-symmetric functions}, Russian Math. Surveys, {\bf{66}}:2 (2011), 271--367.

\bibitem{BE2017}
V. M. Buchstaber, N. Yu. Erokhovets, \emph{Constructions of families of three-dimensional polytopes, characteristic patches of fullerenes, and Pogorelov polytopes}, Izv. Math., {\bf{81}}:5 (2017), 901--972.

\bibitem{BEMPP}
V. M. Buchstaber, N. Yu. Erokhovets, M. Masuda, T. E. Panov, S. Park, \emph{Cohomological rigidity of manifolds defined by 3-dimensional polytopes}, Russian Math. Surveys, {\bf{72}}:2 (2017), 199--256.

\bibitem{B-Kh} 
V. M. Buchstaber, A. N. Kholodov, \emph{Boas-Buck structures on sequences of polynomials}, Funct. Anal. Appl., {\bf{23}}:4 (1989), 266--276.

\bibitem{B-Sh}
V. M. Buchstaber, A. V. Shokurov, \emph{The Landweber—Novikov algebra and formal vector fields on the line}, Funct. Anal. Appl., {\bf{12}}:3 (1978), 159--168.

\bibitem{BK} 
V. M. Buchstaber, E. V. Koritskaya, \emph{Quasilinear Burgers-Hopf equation and Stasheff polytopes}, Funct. Anal. Appl., {\bf{41}}:3 (2007), 196--207.

\bibitem{BL}
V. M. Buchstaber, I. Yu. Limonchenko, \emph{Massey products, toric topology and combinatorics of polytopes,} Izv. Math., {\bf{83}}:6 (2019), 1081–1136.

\bibitem{BL1} 
V. M. Buchstaber, and I. Yu. Limonchenko, \emph{Embeddings of moment-angle manifolds and sequences of Massey products}, preprint (2018); arXiv:1808.08851v1.

\bibitem{BL2} 
V. M. Buchstaber, and I. Yu. Limonchenko, \emph{Direct families of polytopes with nontrivial Massey products}, preprint (2018); arXiv:1811.02221v1.

\bibitem{bu-pa00-2}
V. M. Buchstaber, T. E. Panov, \emph{Torus actions, combinatorial topology, and homological algebra}, Russian Math. Surveys, {\bf{55}}:5 (2000), 825--921.

\bibitem{BP04}
V. M. Buchstaber, T. E. Panov. \emph{Torus actions in topology and combinatorics} (in Russian), M., 2004.

\bibitem{TT} 
V. M. Buchstaber, T. E. Panov. \emph{Toric Topology}, Mathematical Surveys and Monographs, 204, American Mathematical Society, Providence, RI, 2015.

\bibitem{BR} 
V. M. Buchstaber, N. Ray, \emph{Tangential structures on toric manifolds, and connected sums of polytopes}, Internat. Math. Res. Notices, 2001, no. 4, 193--219.

\bibitem{BV} 
V. M. Buchstaber, V. D. Volodin, \emph{Sharp upper and lower bounds for nestohedra}, Izv. Math., {\bf{75}}:6 (2011), 1107--1133; arXiv: 1005.1631v2.

\bibitem{BV2} 
Victor M. Buchstaber, and Vadim D. Volodin, \emph{Combinatorial 2-truncated cubes and applications}, Associahedra, Tamari Lattices, and Related Structures, Tamari Memorial Festschrift, Progress in Mathematics, {\bf{299}}, Birkh\"auser, Basel, 2012, 161--186.

\bibitem{CD} 
Michael P. Carr and Satyan L. Devadoss, \emph{Coxeter complexes and graph-associahedra}, Topology Appl. {\bf{153}} (2006), no. 12, 2155--2168.

\bibitem{Cohen67}
J. M. Cohen, \emph{The decomposition of stable homotopy}, Ann. Math., Second series, {\bf{87}} (2), (1968), 305--320.

\bibitem{C-L-O-T}
O. Cornea, G. Lupton, J. Oprea, and D. Tanr\'e. \emph{Lusternik--Schnirelmann Category}, Mathematical Surveys and Monographs, 103, American Mathematical Society, Providence, RI, 2003.

\bibitem{DJ}
Michael W. Davis and Tadeusz Januszkiewicz, \emph{Convex
polytopes, Coxeter orbifolds and torus actions}. Duke Math.
J.~{\bf62} (1991), no.~2, 417--451.

\bibitem{D-G-M-S}  
P. Deligne, Ph. A. Griffiths, J. W. Morgan, D. Sullivan, \emph{Real homotopy theory of K\"ahler manifolds}, Invent. Math. {\bf{29}} (1975), 245--274.

\bibitem{DS} 
Graham Denham and Alexander I. Suciu, \emph{Moment-angle complexes, monomial ideals, and Massey products}, Pure and Applied Mathematics Quarterly, 3(1) (2007), (Robert MacPherson special issue, part 3), 25--60.

\bibitem{Dold1966}
A. Dold. \emph{Halbexakte Homotopiefunktoren}, Berlin, Springer Verlag, 1966.

\bibitem{DoLa}
A. Dold, R. Lashof, \emph{Principal quasifibrations and fibre homotopy equivalence of bundles}, Illinois J. Math., {\bf{3}} (1959), 285--305.

\bibitem{EiMoore}
S. Eilenberg, J. C. Moore, \emph{Homology and Fibrations I, II, III}, Comm. Math. Helv. {\bf{40}} (1966), 199--236.

\bibitem{FS} 
Eva Maria Feichtner and Bernd Sturmfels, \emph{Matroid polytopes, nested sets and Bergman fans}, Port. Math. (N.S.) {\bf 62} (2005), no. 4, 437--468.

\bibitem{Gan}
T. Ganea, \emph{Lusternik-Schnirelmann category and strong category}, Illinois. J. Math. {\bf{11}} (1967), 417--427.

\bibitem{Gin}
M. Ginsburg, \emph{On the Lusternik-Schnirelmann category}, Ann. of Math. {\bf{77}} (1963), 538--551.

\bibitem{G-LdM} 
Samuel Gitler and Santiago Lopez de Medrano, \emph{Intersections of quadrics, moment-angle manifolds and connected sums}, Geom. Topol. {\bf{17}} (2013), No. 3, 1497--1534.

\bibitem{Go} 
E.S. Golod, \emph{Homologies of some local rings,} Dokl. Akad. Nauk SSSR {\bf{144}} (1962), 479--482.

\bibitem{Gon1}
L. V. Goncharova, \emph{The cohomologies of Lie algebras of formal vector fields on the line}, Funct. Anal. Appl., {\bf{7}}:2 (1973), 91--97.

\bibitem{Gon2}
L. V. Goncharova, \emph{Cohomologies of Lie algebras of formal vector fields on the straight line}, Funct. Anal. Appl., {\bf{7}}:3 (1973), 194--203.

\bibitem{GL1} 
J. Grbic and A.H.M. Linton, \emph{Lowest-degree triple Massey products in moment-angle complexes}, preprint (2019); arXiv:1908.02222
 
\bibitem{GL2} 
J. Grbic and A.H.M. Linton, \emph{Non-trivial higher Massey products in moment-angle complexes}, preprint (2019); arXiv:1911.07083

\bibitem{GT} 
Jelena Grbi\'c and Stephen Theriault, \emph{The homotopy type of the complement of a coordinate subspace arrangement}, Topology {\bf{46}} (2007), no. 4, 357--396.

\bibitem{G-T13} 
Jelena Grbi\'c and Stephen Theriault, \emph{The homotopy type of the polyhedral product for shifted complexes,} Adv. Math. {\bf{245}} (2013), 690--715.

\bibitem{GT2} 
Jelena Grbi\'c and Stephen Theriault, \emph{Homotopy theory in toric topology}, Russian Math. Surveys, {\bf{71}}:2 (2016), 185--251.

\bibitem{G-P-T-W} 
Jelena Grbi\'c, Taras Panov, Stephen Theriault, and Jie Wu, \emph{Homotopy types of moment-angle complexes for flag complexes}, Trans. Amer. Math. Soc., {\bf{368}} (2016), no.9, 6663--6682; arXiv:1211.0873.

\bibitem{GL}
T. H. Gulliksen and G. Levin. \emph{Homology of local rings}, Queen’s Papers in Pure and Applied Mathematics. V. 20, Queen’s University, Kingston, Ontario, 1969.

\bibitem{HL}
S. Halperin and J.-M. Lemaire, \emph{Notions of category in differential algebra}, Lecture Notes in Math., vol. 1318, Springer-Verlag, Berlin, 1988, 138--154.

\bibitem{HH} 
J. Herzog and T. Hibi. \emph{Monomial Ideals}, Graduate Texts in Mathematics {\bf{260}}, Springer, 2011. 

\bibitem{Hoch} 
M. Hochster, \emph{Cohen-Macaulay rings, combinatorics, and simplicial complexes}, in Ring theory, II (Proc. Second
Conf.,Univ. Oklahoma, Norman, Okla., 1975),  Lecture Notes in Pure
and Appl. Math., V. 26, 171--223, Dekker, New York, 1977.

\bibitem{IK} 
K. Iriye and D. Kishimoto, \emph{Decompositions of polyhedral products for shifted complexes}, Adv. Math. {\bf{245}} (2013), 716--736.

\bibitem{IK2014}
K. Iriye and D. Kishimoto, \emph{Fat-wedge filtration and decomposition of polyhedral products,} Kyoto J. Math. {\bf{59}}:1 (2019), 1--51; arXiv: 1412.4866v4.

\bibitem{IY}
K. Iriye and T. Yano, \emph{A Golod complex with non-suspension moment-angle
complex,} Topology Appl. {\bf{225}} (2017), 145--163; arXiv: 1601.03472v4.

\bibitem{Jess}
B. Jessup, \emph{New estimates for the collapse of the Milnor--Moore spectral sequence over a field}, Proc. of the AMS, Vol. 114, No. 4 (1992), 1115--1117.

\bibitem{Kat} 
L. Katth\"an, \emph{A  non-Golod  ring  with  a  trivial  product  on  its  Koszul  homology},  J.  Algebra {\bf{479}} (2017),  244–-262.

\bibitem{Kr} 
D. Kraines, \emph{Massey higher products}, Trans. Amer. Math. Soc., {\bf{124}} (1966), 431--449. 

\bibitem{L2014} 
I. Yu. Limonchenko, \emph{Stanley–Reisner rings of generalized truncation polytopes and their moment–angle manifolds}, Proc. Steklov Inst. Math., {\bf{286}} (2014), 188--197; arXiv:1401.2124.

\bibitem{L2015} 
I. Yu. Limonchenko, \emph{Families of minimally non-Golod complexes and polyhedral products}, Far Eastern Math. J., {\bf{15}}:2 (2015),  222--237; arXiv:1509.04302.

\bibitem{L1}
I. Yu. Limonchenko, \emph{Massey products in cohomology of moment-angle manifolds for 2-truncated cubes}, Russian Math. Surveys, {\bf{71}}:2 (2016), 376--378.

\bibitem{L2}
I. Yu. Limonchenko, \emph{Topology of moment-angle-manifolds arising from flag nestohedra,} Ch. Ann. Math., {\bf{38B}}(6) (2017), 1287--1302; arXiv:1510.07778.

\bibitem{L3}
I. Yu. Limonchenko, \emph{On higher Massey products and rational formality
for moment-angle-manifolds over multiwedges}, Proc. of the Steklov Inst. Math., 2019, {\bf{305}}, 161--181; arXiv:1711.00461v2.

\bibitem{McCleary}
J. McCleary. \emph{User's guide to spectral sequences}, Cambridge University Press, 2001.

\bibitem{LdM} 
Santiago Lopez de Medrano, \emph{Topology of the intersection of quadrics in $\mathbb{R}^n$}, Lecture Notes in Mathematics {\bf{1370}} (1989), 280--292.

\bibitem{M-U} 
W. S. Massey, H. Uehara,\emph{The Jacobi identity for Whitehead products,} Algebraic geometry and topology. A symposium in honor of S. Lefschetz, Princeton, N. J.: Princeton University Press, (1957), 361--377.

\bibitem{Mass} 
W. S. Massey, \emph{Some higher order cohomology operations}, International symposium on algebraic topology, Mexico City: Universidad Nacional Aut\'{o}noma de M\'{e}xico and UNESCO, (1958), 145--154.

\bibitem{May1968}
J. P. May, \emph{The cohomology of principal bundles, homogeneous spaces, and two-stage Postnikov systems}, Bull. AMS {\bf{74}} (1968), 334--339.

\bibitem{M} 
J. P. May, \emph{Matric massey products}, J. Algebra, {\bf{12}} (1969), 533--568.

\bibitem{Million}
D. V. Millionshchikov, \emph{Algebra of formal vector fields on the line and Buchstaber's conjecture}, Funct. Anal. Appl., {\bf{43}}:4 (2009), 264--278. 

\bibitem{Milnor1956}
J. W. Milnor, \emph{Construction of universal bundles I, II}, Ann. Math. {\bf{63}} (1956), 272--284, 430--436.

\bibitem{N} 
E. J. O'Neill, \emph{On Massey products}, Pacific J. of Math., Vol. 76, No. 1 (1978), 123--127.

\bibitem{P} 
Taras E. Panov, \emph{Cohomology of face rings, and torus actions},
in ``Surveys in Contemporary Mathematics''. London Math. Soc.
Lecture Note Series, vol.~\textbf{347}, Cambridge, U.K., 2008, 165--201; arXiv:math.AT/0506526.

\bibitem{P-R}
Taras Panov and Nigel Ray, \emph{Categorical aspects of toric
topology}. In: \emph{Toric Topology}, M.~Harada \emph{et al.},
eds. Contemp. Math.,~460. Amer. Math. Soc., Providence,
RI, 2008, 293--322.

\bibitem{p-r-v04}
Taras Panov, Nigel Ray and Rainer Vogt, \emph{Colimits,
Stanley--Reiner algebras, and loop spaces.} In: \emph{Categorical
decomposition techniques in algebraic topology (Isle of Skye,
2001)}. Progress in Math.,~215, Birkh\"auser, Basel, 2004, pp.
261--291.

\bibitem{PV} 
T. E. Panov, Ya. A. Veryovkin, \emph{Polyhedral products and commutator subgroups of right-angled Artin and Coxeter groups}, Sb. Math., {\bf{207}}:11 (2016), 1582--1600; arXiv:1603.06902.

\bibitem{Pog}
A. V. Pogorelov, \emph{A regular partition of Lobachevskian space}, Math. Notes, {\bf{1}}:1 (1967), 3--5.

\bibitem{Post05} 
A.~Postnikov, \emph{Permutohedra, associahedra, and beyond}, Int. Math. Res. Not. (2009), no. 6, 1026--1106; arXiv: math.CO/0507163.

\bibitem{PRW06} 
A.~Postnikov, V.~Reiner, L.~Williams, \emph{Faces of generalized permutohedra},  Documenta Mathematica, {\bf{13}} (2008), 207--273; arXiv:math/0609184v2.

\bibitem{RuTr}
Yu. Rudyak and A. Tralle, \emph{On Thom spaces, Massey products, and nonformal
symplectic manifolds,} Internat. Math. Res. Notices {\bf{2000}}:10, 495--513; Preprint MPI 1999-71, 1999.

\bibitem{Sh}
A. S. Schwarz, \emph{The genus of fiber space}, Dokl. Akad. Nauk SSSR, {\bf{119}}:2 (1958), 219--222.

\bibitem{S}
James D. Stasheff, \emph{Homotopy associativity of H-spaces. I.}, Trans. Amer. Math. Soc.~\textbf{108} (1963), 275--292.

\bibitem{sull77}
D. Sullivan, \emph{Infinitesimal computations in topology},
Publ. I.H.E.S.~\textbf{47} (1977), 269--331.

\bibitem{Toda62}
H. Toda, \emph{Composition Methods in homotopy groups of spheres,} Princeton University Press, (1962).

\bibitem{TrOp}
A. Tralle, J. Oprea. \emph{Symplectic manifolds with no K\"ahler structure}. Berlin: Springer-Verlag, 1997. (Lecture Notes in Math. V. 1661.)

\bibitem{Zh} 
E. G. Zhuravleva, \emph{Massey products in the
cohomology of the moment-angle manifolds corresponding to Pogorelov polytopes},
Math. Notes, {\bf{105}}:4 (2019), 519--527; arXiv: 1707.07365v2.

\end{thebibliography}
\end{document}